\documentclass [11pt,a4paper]{article}
\usepackage{amsmath}
\usepackage{amsthm}
\usepackage{graphicx}
\usepackage{latexsym}
\usepackage{amssymb}
\usepackage{enumerate}
\usepackage[latin2]{inputenc}
\usepackage{color}
\usepackage{calc}
\usepackage{subfigure}

\numberwithin{equation}{section}
\newtheorem  {theorem}       {Theorem}
\newtheorem  {lemma}         {Lemma}
\newtheorem  {corollary}     {Corollary}
\newtheorem  {proposition}   {Proposition}

\newtheorem* {theorem*}      {Theorem}
\newtheorem* {lemma*}        {Lemma}
\newtheorem* {corollary*}    {Corollary}
\newtheorem* {proposition*}  {Proposition}
\newtheorem* {definition*}   {Definition}
\newtheorem* {remark*}       {Remark}
\newtheorem* {remarks*}      {Remarks}

\newcounter{aux}

\def \N {\mathbb N}
\def \Z {\mathbb Z}

\def \R {\mathbb R}

\def \T {\mathbb T}

\def \Ordo {{\cal O}}
\def \ind{1\!\!1}

\def \dsum {\displaystyle\sum}

\def \alp   {\alpha}
\def \bet   {\beta}

\def \zet   {\zeta}

\def \kap   {\kappa}
\def \lam   {\lambda}
\def \ome   {\omega}
\def \om {\omega}

\def \Ome   {\Omega}
\def \Om   {\Omega}

\def\uome{\underline{\omega}}
\def\uo{\uome}

\def\vareps{\varepsilon}
\def\wt{\widetilde}

\def\ulr{\underline{r}}
\def\olr{\overline{r}}

\def\beqs{\begin{eqnarray*}}
\def\eeqs{\end{eqnarray*}}
\def\beq{\begin{eqnarray}}
\def\eeq{\end{eqnarray}}
\def\beas{\begin{eqnarray*}}
\def\eeas{\end{eqnarray*}}
\def\bea{\begin{eqnarray}}
\def\eea{\end{eqnarray}}

\def \prob        {\ensuremath{\mathbf{P}}}
\def \expect      {\ensuremath{\mathbf{E}}}
\def \var         {\ensuremath{\mathbf{Var}}}
\def \cov         {\ensuremath{\mathbf{Cov}}}

\def\pt{\partial_t}
\def\px{\partial_x}
\newcommand{\abs}[1]{\left|{#1}\right|}

\def \b {\beta}
\def \d {\delta}
\def \Omn {\Omega^{n}}

\def \Tn {\T^{n}}
\def \Ln {L^{n}}
\def \Kn {K^{n}}
\def \Gn {G^{n}}

\def \grn {\nabla^{n}}
\def \gradn {\grn}

\def \Xn {{\cal X}^{n}}

\def \mun  {\mu^{n}}
\def \munt {\mu^{n}_t}

\def \nun  {\nu^{n}}
\def \nunt {\nu^{n}_t}

\def \pin  {\pi^{n}}
\def \pil  {\pi^{l}}

\def\sn {s^n}

\def\hn   {h^{n}}
\def\fnt  {f^{n}_t}

\def\taun {\tau^{n}}
\def\thetan {\theta^{n}}
\def\Psin {\Psi^{n}}
\def\Phin {\Phi^{n}}

\def\rn{\wih{\rho}}
\def\un{\wih{u}}
\def\rnt{\widetilde{\rho}^n}
\def\unt{\widetilde{u}^n}

\newcommand{\vct}[1]{ \text{\boldmath${#1}$} }

\def\vzeta {\vct{\zeta}}

\def\vx  {\vct{x}}
\def\vlam{\vct{\lam}}

\newcommand{\wih}[1] {\widehat{#1}^{n}}

\def \dom    {{\cal D}}

\def \Jn {J^{n}}
\def \In {I^{n}}
\def \Sn {S^{n}}
\def \Fn {F^{n}}
\def \I {\In}
\def \J {\Jn}
\def \S {\Sn}
\def \F {\Fn}

\def \rruu {(\rho,u)}

\def \g {\gamma}

\def\pp{\partial}

\def\ulr{\underline{r}}
\def\olr{\overline{r}}
\def\ulw{\underline{w}}
\def\olw{\overline{w}}
\def\ulz{\underline{z}}
\def\olz{\overline{z}}
\def \a {\alpha}

\title{\textsc{Perturbation of singular
equilibria of hyperbolic two-component
  systems:\\ a universal hydrodynamic limit}}

\author {
{Bálint Tóth \qquad Benedek Valkó}
\\[8pt]
Institute of  Mathematics
\\
Technical University Budapest
}

\begin{document}

\setlength{\baselineskip}{1.23\baselineskip}

\maketitle
\begin{abstract}
We consider one-dimensional, locally finite interacting particle
systems with two conservation laws which under Eulerian
hydrodynamic limit lead to two-by-two systems of conservation
laws:
\beq
\notag
\left\{
\begin{array}{l}
  \pt \rho +\px \Psi(\rho , u)=0
\\
  \pt u+\px \Phi(\rho,u)=0,
\end{array}
\right.
\eeq
with $(\rho,u)\in{\cal D}\subset\R^2$, where ${\cal D}$ is a convex
compact polygon in  $\R^2$. The system  is \emph{typically} strictly
hyperbolic in the interior of ${\cal D}$ with possible non-hyperbolic
degeneracies  on the boundary $\partial {\cal D}$.  We consider the
case of isolated singular (i.e. non hyperbolic)  point on the interior
of one of the edges
of ${\cal D}$, call it $(\rho_0,u_0)=(0,0)$ and assume ${\cal
  D}\subset\{\rho\ge0\}$. This can be achieved by a linear
transformation of the conserved quantities.
We investigate the propagation of \emph{small nonequilibrium
  perturbations} of the steady state
of the microscopic interacting particle  system,
corresponding to the densities
$(\rho_0,u_0)$ of the conserved quantities. We prove that
for a very rich  class of systems, under proper hydrodynamic limit the
propagation of these small perturbations are \emph{universally} driven
by the two-by-two system
\beq%
\notag
\left\{
\begin{array}{l}
\pt\rho + \px\big(\rho u\big)=0
\\
\pt u   + \px\big(\rho + \gamma u^2\big) =0
\end{array}
\right.
\eeq
where the parameter $\gamma:=\frac12 \Phi_{uu}(\rho_0,u_0)$ (with a
proper choice of space and time scale) is the only trace of the
microscopic structure. The proof is valid for the cases with
$\gamma>1$.

The proof essentially relies on the relative entropy method
and thus, it is valid only in the regime of smooth solutions of the
pde. 
But there are  essentially new elements: in order to control
the fluctuations of the terms with Poissonian (rather than
Gaussian) decay coming from the low density approximations we have
to apply refined pde estimates. In particular Lax entropies of
these pde systems 
play a \emph{not merely technical} key role in the main part of
the proof.
\end{abstract}
\pagebreak

\tableofcontents

\section{Introduction}
\label{section:intro}

\subsection{The PDE to be derived and some facts about it}
\label{subs:pde}

We consider the pde
\begin{eqnarray}
\label{eq:pde}
\left\{
\begin{array}{l}
\pt\rho + \px\big(\rho u\big)=0
\\[5pt]
\pt u   + \px\big(\rho + \gamma u^2\big) =0
\end{array}
\right.
\end{eqnarray}
where $\rho=\rho(t,x)\in[0,\infty)$, $u=u(t,x)\in(-\infty,\infty)$
are density, respectively, velocity field and $\gamma\in\R$  is a
fixed parameter. For any fixed $\gamma$ this is a \emph{hyperbolic
system of conservation laws} in the domain
$(\rho,u)\in\R_+\times\R$.

Phenomenologically, the pde describes a deposition/domain growth
-- or, in biological term: chemotaxis -- mechanism: $\rho(t,x)$
is the density of population performing the deposition and
$h(t,x)$ is the height of the deposition. Let
\[
u(t,x):=-\px h(t,x).
\]
The physics of the phenomenon is contained in the following two rules:

\begin{enumerate}[(a)]

\item
The velocity field of the population is proportional to the
\emph{negative gradient of the height}  of the deposition. That is,
the  population is pushed towards the local decrease of the deposition
height. This rule, together with the  conservation of total mass of
the population leads to the continuity equation (the first
equation in our system).

\item
The deposition rate is
\begin{eqnarray*}
\pt h = \rho + \gamma \big(\px h\big)^2.
\end{eqnarray*}
The first term on the right hand side is just saying that
deposition is done additively by the population. The second term
is a self-generating deposition, introduced and phenomenologically
motivated by Kardar-Parisi-Zhang
\cite{kardarparisizhang} and commonly accepted in the
literature. Differentiating this last  equation with respect to the
space variable
$x$ results in the second equation of our system.

\end{enumerate}

\noindent The pde (\ref{eq:pde})  is invariant under the scaling:
\begin{eqnarray*}
\label{eq:scale}
\wt\rho(t,x)
:=
A^{2\beta}\rho(A^{1+\beta}t,Ax),
\qquad
\wt u(t,x)
:=
A^{\beta}u(A^{1+\beta}t,Ax),
\end{eqnarray*}
where $A>0$ and $\beta\in\R$ are arbitrarily fixed. The choice
$\beta=0$  gives the straightforward hyperbolic scale invariance,
valid for  any system of conservation laws. More interesting is
the $\beta=1/2$ case. This is the natural scale invariance of the
system, since the physical variables (density  and velocity
fields) change \emph{covariantly}  under this scaling. This is the
(presumed, but never rigorously  proved)  asymptotic scale
invariance of the Kardar-Parisi-Zhang deposition phenomena. The
nontrivial scale invariance of the pde (\ref{eq:pde}) suggests its
\emph{universality} in some sense. Our main result indeed states
its validity in a very wide context.

It is also clear that the pde is invariant under the left-right
reflection symmetry $x\mapsto -x$:

The parameter $\gamma$  of the pde (\ref{eq:pde}) is of crucial
importance: different values of $\gamma$ lead to completely
different  behavior. Here are listed some particular cases which
arose in the past in various contexts:
\\
---
The pde (\ref{eq:pde}) with $\gamma=0$ arose in the context of the
`true self-repelling motion' constructed by T\'oth and Werner in
\cite{tothwerner1}.  For a survey of this case see also
\cite{tothwerner2}. The same equation, with viscosity terms added,
appear in mathematical biology under the name of (negative)
chemotaxis equations (see e.g. \cite{rascle},
\cite{othmerstevens}, \cite{levinesleeman})
).
\\
---
Taking $\gamma=1/2$  we get the `shallow water equation'. See
\cite{evans}, \cite{leveque}.  This is the only value of the
parameter $\gamma$ when $m=\rho u$  is conserved and as a
consequence the pde (\ref{eq:pde}) can be interpreted as gas
dynamics equation.
\\
---
With $\gamma=1$ the pde is called `Leroux's equation' which is of
Temple class and for this reason much investigated. For many
details about this equation see \cite{serre}. In the recent paper
\cite{fritztoth} Leroux's system has been derived as hydrodynamic
limit under Eulerian scaling for a two-component lattice gas,
going even beyond the appearance of shocks.

The main facts about the pde (\ref{eq:pde}) are presented in
Section \ref{section:pdedetails}. Here we only mention that

\begin{enumerate} [1.]

\item
For any
$\gamma\in\R$ the system
(\ref{eq:pde}) is strictly \emph{hyperbolic} in
$(\rho,u)\in(0,\infty)\times\R$, with hyperbolicity marginally lost at
$(\rho,u)=(0,0)$ for $\gamma\not=1/2$ and at
$\rho=0$ for $\gamma=1/2$. This follows from straightforward
computations.

\item The \emph{Riemann invariants} (or characteristic
coordinates)  are explicitly computed in section
\ref{section:pdedetails}, for a first impression see Figure
\ref{fig:level} of the Appendix where the level lines of the
Riemann invariants are shown. It turns out that the picture
changes qualitatively at the critical values $\gamma=1/2$,
$\gamma=3/4$ and $\gamma=1$. It is of crucial importance for our
later problem that the level curves, expressed as
$u\mapsto\rho(u)$  are convex for $\gamma<1$, linear for
$\gamma=1$  and concave for $\gamma>1$.

\item
For any
$\gamma\ge0$ the system
(\ref{eq:pde}) is \emph{genuinely nonlinear} in
$(\rho,u)\in(0,\infty)\times\R$, with genuine nonlinearity  marginally
lost at
$(\rho,u)=(0,0)$ for $\gamma\not=0,1/2$ and at
$\rho=0$ for
$\gamma=0,1/2$. (For
$\gamma<0$ genuine nonlinearity is lost on the parabola
$\rho=-4\gamma(2\gamma-1)^2(\gamma+1)^{-2}u^2$.)

\item
The system is sufficiently rich in
\emph{Lax entropies}.

\item
For
$\gamma\ge0$ the system
(\ref{eq:pde})  satisfies the conditions of the
Lax-Chuey-Conley-Smoller
\emph{Maximum Principle}  (see
\cite{lax1}, \cite{lax2},
\cite{serre}).

\end{enumerate}

From the Maximum Principle a
very essential difference between the cases
$\gamma<1$,
$\gamma=1$ and
$\gamma>1$ follows, which is of crucial importance for our further
work. In the case
$\gamma<1$ all convex domains bounded by level curves of the Riemann
invariants are
\emph{unbounded   (non-compact)}   and thus there is no a priori bound
on the entropy  solutions. Even starting with smooth initial data with
compact support nothing prevents the solutions to blow up
indefinitely. On the other hand, if
$\gamma\ge1$ any bounded subset of
$\R_+\times\R$ is contained in a compact convex domain bounded by
level  sets of the Riemann invariants, which fact yields a priori
bounds on the entropy solutions, given bounded initial data.

The goal of the present paper is to derive the two-by-two hyperbolic
system of conservation laws (\ref{eq:pde}) as decent hydrodynamic
limit of some systems of interacting particles with two conserved
quantities.

We consider one-dimensional, locally finite interacting particle
systems with two  conservation laws  which under \emph{Eulerian}
hydrodynamic limit lead to two-by-two systems of conservation laws
\beqs
\left\{
\begin{array}{l}
  \pt \rho +\px \Psi(\rho , u)=0
\\[5pt]
  \pt u+\px \Phi(\rho,u)=0,
\end{array}
\right.
\eeqs
with $(\rho,u)\in{\cal D}\subset\R^2$, where ${\cal D}$ is a convex
compact polygon in  $\R^2$. The system  is \emph{typically} strictly
hyperbolic in the interior of ${\cal D}$ with possible non-hyperbolic
degeneracies  on the boundary $\partial {\cal D}$.  We consider the
case of isolated singular (i.e. non hyperbolic)  point on the interior
of one of the edges
of ${\cal D}$, call it $(\rho_0,u_0)=(0,0)$ and assume ${\cal
  D}\subset\{\rho\ge0\}$ (otherwise we apply an appropriate
   linear transformation on the conserved quantities)
We investigate the propagation of
\emph{small nonequilibrium
  perturbations} of the steady state
of the microscopic interacting particle  system, corresponding to
the densities $(\rho_0,u_0)$ of the conserved quantities. We prove
that for a very rich  class of systems, under proper hydrodynamic
limit the propagation of these small perturbations are
\emph{universally} driven by the  system (\ref{eq:pde}), where the
parameter $\gamma:=\frac12 \Phi_{uu}(\rho_0,u_0)$ (with a proper
choice of space and time scale) is the only trace of the
microscopic structure. The proof is valid for the cases with
$\gamma>1$.

Actually, in order to simplify some of the arguments, we
impose the left-right reflection symmetry of the pde (\ref{eq:pde}) on
the systems of interacting particles on microscopic level, see
condition (\ref{cond:lrsymm}) in subsection \ref{subs:rates}.
But we note that the whole proof 
can be extended without this condition, just some arguments would
be longer.

The proof essentially relies on H-T. Yau's relative entropy method
and thus, it is valid only in the regime of smooth solutions of the
pde (\ref{eq:pde}).

We should emphasize here the essential new ideas of the proof.
Since we consider a \emph{low density} limit, the distribution of
particle numbers in blocks of mesoscopic size will have a
\emph{Poissonian} tail. The fluctuations of the other conserved
quantity will be Gaussian, as usual. It follows that when
controlling the fluctuations of the empirical block averages the
usual large deviation approach would lead us to the disastrous
estimate $\expect\big(\exp\{\vareps \,\, GAU\cdot
POI\}\big)=\infty$. It turns out that some very special cutoff
must be applied. Since the large fluctuations which are cut off
can not be
estimated by robust methods (i.e. by applying entropy inequality),
only some cancellation due to martingales can help. This is the
reason why the cutoff function must be chosen  in a very special
way, in terms of a particular Lax entropy of the Euler equation
(\ref{eq:euler}). In this way the proof becomes a mixture (in our
opinion rather interesting mixture) of probabilistic and pde
arguments. The fine properties of the limiting pde, in particular
the global behavior of Riemann invariants and some particular Lax
entropies, play
an essential role in the proof. The radical difference between the
$\gamma\ge1$ vs. $\gamma<1$ cases, in particular applicability vs.
non-applicability of the Lax-Chuey-Conley-Smoller maximum principle,
manifests itself on the microscopic, probabilistic level.

\subsection{The structure of the paper}
\label{subs:structure}

In Section \ref{section:models} we define the class of models to
which our main theorem applies: we formulate the conditions to be
satisfied by the interacting particle systems to be considered, we
compute the steady state measures and the fluxes corresponding to
the conserved quantities. At the end of this section we formulate
the Eulerian hydrodynamic limit, for later reference.

In Section \ref{section:intermediate} first we perform asymptotic
analysis of the Euler equations close to the singular point
considered, then we formulate our main result, Theorem \ref{thm:main},
and its immediate consequences.

Sections \ref{section:proofprep} to \ref{section:ingred} are devoted
to the proof of  Theorem \ref{thm:main}.

In Section  \ref{section:proofprep} we perform the necessary
preliminary computations for the proof. After introducing the minimum
necessary notation we apply some standard procedures in the context of
relative entropy method.
Empirical block averages are introduced, numerical error terms are
separated and estimated. In this first estimates only straightforward
numerical approximations (Taylor expansion bounds) and the most direct
entropy inequality is applied.

Section \ref{section:cutoff} is of crucial importance: here it is
shown why the traditional approach of the relative entropy method
fails to apply. Here it becomes apparent that in the fluctuation
bound (usually referred  to as \emph{large deviation estimate})
instead of the tame $\expect\big(\exp\{\vareps \,\,GAU^2\}\big)$
we would run into the wild $\expect\big(\exp\{\vareps \,\,GAU\cdot
POI\}\big)$ which is, of course, infinite. It is explained here
what kind of cutoff is applied: the large fluctuations cut off can
not be estimated by robust methods (i.e. by applying entropy
inequality). Only some cancellation due to martingales can help.
This is the reason why the cutoff function must be chosen  in a
very particular way, in terms of a particular Lax entropy of the
Euler equation. The cutoff function is constructed and its key
estimates are stated. Proofs of the lemmas formulated in this
section are postponed to Section \ref{section:cutoffproofs}. At
the end of this section the outline of the further steps is
presented.

In Section \ref{section:tools} all the necessary probabilistic
ingredients
of the forthcoming steps are gathered. These are: fixed time large
deviation bounds and fixed time fluctuation bounds, the time averaged
block replacement bounds (one block estimates) and the time averaged
gradient bounds (two block estimates). The proof of these last two rely
on Varadhan's large deviation bound cited  in that section and
on some  probability lemmas stated and proved in section
\ref{section:ingred}. We should mention here that these proofs of the
one- and two block  estimates, in particular the probability lemmas
involved also contain some new (and, 
we hope, instructive) elements.

Sections \ref{section:firstterm} and \ref{section:otherterms}
conclude the proof: the various terms arising in section
\ref{section:cutoff} are estimated using all the tools
(probabilistic and pde) developed in earlier sections. One can see
that these estimates  rely heavily on the fine properties of the
Lax entropy used in the cutoff procedure.

As we already mentioned sections  \ref{section:cutoffproofs} and
\ref{section:ingred} are devoted to proofs of various lemmas
stated in earlier parts. Section  \ref{section:cutoffproofs}
deals with the pde estimates while Section \ref{section:ingred}
is probabilistic.

Finally in the Appendix (Section \ref{section:pdedetails}) we give
some details about the pde (\ref{eq:pde}). This is included for sake
of completeness and in order to let the reader have some more
information about these, certainly interesting, pde-s. Strictly
technically speaking this Appendix is not used in the proof.


\section{Microscopic models}
\label{section:models}

Our  interacting particle systems to be defined in the present
section model on a microscopic level the same deposition phenomena
as the pde (\ref{eq:pde}).  There will be two conserved physical
quantities: the particle number $\eta_j\in\N$  and the (discrete)
negative gradient of the deposition height $\zeta_j\in\Z$.

\noindent The dynamical driving mechanism is of such nature that
\begin{enumerate}[(i)]
\item The deposition height growth is influenced by the local
particle density. Typically: growth is enhanced by higher particle
densities.
\item The particle motion is itself influenced by the
deposition profile. Typically: particles are pushed in the
direction of the negative gradient of the deposition height.
\end{enumerate}
The left-right reflection symmetry of the pde will be also
implemented on the microscopic level. Actually, this is not really
necessary in order to prove our main result, but without this
assumption some of the arguments would be somewhat longer.

\subsection{State space, conserved quantities}
\label{subs:statespace}

Throughout this paper we denote by
$\Tn$ the discrete tori
$\Z/n\Z$,
$n\in\N$,  and by
$\T$  the continuous torus
$\R/\Z$.  We will denote the local spin state by
$\Om$,  we only consider the case when
$\Om$ is finite. The state space of the interacting particle system
of size
$n$ is
\begin{equation*}
\Omn:=\Om^{\Tn}.
\end{equation*}
Configurations will be denoted
\begin{equation*}
\uo:=(\omega_j)_{j\in \Tn}\in\Omn,
\end{equation*}
For sake of simplicity we consider discrete (integer valued)
conserved quantities only. The two conserved quantities are
\begin{eqnarray}
\notag
&&
\eta: \Om\rightarrow \N,
\\[5pt]
\label{eq:zetamap}
&&
\zeta: \Om\rightarrow v_0\Z,
\text{ or }
\zeta: \Om\rightarrow v_0\big(\Z+1/2\big).
\end{eqnarray}
The trivial scaling factor $v_0$ will be conveniently chosen later
(see
(\ref{eq:choiceofv0})).
We also use the notations
$\eta_j=\eta(\omega_j),\,\zeta_j=\zeta(\omega_j).$  This means that
the sums
$\sum_j\eta_j$ and
$\sum_j\zeta_j$ are conserved by the dynamics. We assume that the
conserved quantities are different and non-trivial, i.e. the
functions
$\zeta ,\eta$ and the constant function 1 on
$\Om$ are linearly independent.

The left-right reflection symmetry of the model is implemented by an
involution
\[
R:\Om\to \Om, \qquad R\circ R = Id
\]
which acts on the conserved quantities  as follows:
\begin{eqnarray}
\label{eq:mirroretazeta} \eta(R\omega)=\eta(\omega), \qquad
\zeta(R\omega)=-\zeta(\omega).
\end{eqnarray}
%

\subsection{Rate functions, infinitesimal generators, stationary
  measures}
\label{subs:rates}

Consider a (fixed) probability measure
$\pi$  on
$\Om$, which is invariant under the action of the involution
$R$, i.e.
$\pi(R\ome)=\pi(\ome)$.  Since eventually we consider
\emph{low densities} of
$\eta$, in order
to exclude trivial cases we assume that
\beq
\label{eq:pinondeg}
\pi\big(\,\zeta=0\,\big|\,\eta=0\,)<1.
\eeq
The scaling factor $v_0$ in
(\ref{eq:zetamap}) is chosen so that
\beq
\label{eq:choiceofv0}
\var\big(\zeta\,|\,\eta=0)=1.
\eeq
This choice simplifies some formulas (fixing a recurring constant to
be equal to 1, see
(\ref{eq:phi1}))  but does not restrict generality.

For later use we introduce the notation
\beqs
&&
\rho^*:=\max\{\eta(\omega):\, \pi(\omega)>0\},
\\[5pt]
&&
u^*:=\max\{\zeta(\omega): \,\pi(\omega)>0\},
\\[5pt]
&&
u_*:=\max\{\zeta(\omega): \,\eta(\omega)=0, \,
\pi(\omega)>0\},
\eeqs
For
$\tau, \theta \in \R$ let
$G(\tau,\theta)$ be the moment generating function defined below:
\begin{eqnarray*}
G(\tau,\theta)
:=
\log \sum_{\omega\in \Om}
e^{\tau \eta(\omega)+\theta \zeta(\omega)}
\pi(\omega).
\end{eqnarray*}
In thermodynamic terms
$G(\tau,\theta)$  corresponds to the Gibbs free energy. We define the
probability measures
\begin{eqnarray}
\label{eq:gibbs1}
\pi_{\tau, \theta}(\omega)
:=
\pi(\omega)
\exp(\tau \eta(\omega)+ \theta \zeta(\omega)-G(\tau,\theta))
\end{eqnarray}
on
$\Om$. We are going to define dynamics which conserve the quantities
$\sum_j\eta_j$ and
$\sum_j\zeta_j$, posses no other (hidden) conserved quantities and for
which the product measures
\begin{eqnarray*}
\pin_{\tau,\theta}
:=
\prod_{j \in \Tn} \pi_{\tau,\theta}
\end{eqnarray*}
are stationary.

We need to separate a symmetric (reversible) part of the dynamics
which will be speeded up sufficiently in order to enhance
convergence to local equilibrium and thus help estimating some
error term in the hydrodynamic limiting procedure. So
 we consider two
\emph{rate functions} $r: \Om\times
\Om\times\Om\times\Om\rightarrow \R_+$ and $s:
\Om\times\Om\times\Om\times\Om\rightarrow \R_+$, $r$ will define
the \emph{asymmetric} component of the dynamics, while $s$ will
define the \emph{reversible} component.   The dynamics of the
system consists of elementary jumps affecting nearest neighbor
spins,
$(\omega_j,\omega_{j+1})\longrightarrow(\omega'_j,\omega'_{j+1})$,
performed with rate $\lambda
r(\omega_j,\omega_{j+1};\omega'_j,\omega'_{j+1}) +
 \kappa  s(\omega_j,\omega_{j+1};\omega'_j,\omega'_{j+1})$, where
$\lambda, \kappa  > 0$  are speed-up  factors, depending on
 the size of the system in the limiting procedure.

We require that the rate functions
$r$ and
$s$ satisfy the following conditions.

\begin{enumerate}[(A)]

\item
\label{cond:cons}
\emph{Conservation laws:}
If
$r(\omega_1,\omega_2;\omega'_1,\omega'_2)>0$
\emph{or}
$s(\omega_1,\omega_2;\omega'_1,\omega'_2)>0$ then
\begin{eqnarray*}
\begin{array}{rcl}
&&
\eta(\omega_1)+\eta(\omega_2)=
\eta(\omega'_1)+\eta(\omega'_2).
\\[5pt]
&&
\zeta(\omega_1)+\zeta(\omega_2)=
\zeta(\omega'_1)+\zeta(\omega'_2),
\end{array}
\end{eqnarray*}

\item
\label{cond:irred}
\emph{Irreducibility:}
For every
$N\in [0,n \rho^*]$,
$Z\in [-nu^*,nu^*]$  the set
\begin{eqnarray*}
\Omn_{N,Z}
:=
\left\{
\uo\in\Omn:
\sum_{j\in\Tn}\eta_j=N,
\sum_{j\in\Tn}\zeta_j=Z
\right\}
\end{eqnarray*}
is an irreducible component of
$\Omn$, i.e. if
$\uo, \uo'\in \Omn_{N,Z}$ then there exists a series of elementary
jumps with positive rates transforming
$\uo$ into
$\uo'$.

\item
\label{cond:lrsymm}
\emph{Left-right symmetry:}
The jump rates are invariant under left-right reflection
\emph{and} the action of the involution
$R$ (jointly):
\begin{eqnarray*}
&&
r(R\omega_2,R\omega_1;R\omega'_2,R\omega'_1)
=
r(\omega_1,\omega_2;\omega'_1,\omega'_2).
\\
&&
s(R\omega_2,R\omega_1;R\omega'_2,R\omega'_1)
=
s(\omega_1,\omega_2;\omega'_1,\omega'_2).
\end{eqnarray*}

\item
\label{cond:rstaci}
\emph{Stationarity of the asymmetric part:}
For any $\omega_1, \omega_2, \omega_3\in \Ome$
\beqs
Q(\omega_1,\omega_2)+Q(\omega_2,\omega_3)+Q(\omega_3,\omega_1)=0,
\eeqs
where
\beqs
Q(\omega_1,\omega_2)
:=
\sum_{\omega_1',\omega_2' \in \Ome}
\left\{
\frac{\pi(\omega_1')\pi(\omega_2')}{\pi(\omega_1)\pi(\omega_2)}
r(\omega_1',\omega_2';\omega_1,\omega_2)
-
r(\omega_1,\omega_2;\omega'_1,\omega'_2)
\right\}.
\eeqs

\item
\label{cond:srevers}
\emph{Reversibility of the symmetric part:}
For any
$\omega_1, \omega_2, \omega^\prime_1, \omega^\prime_2\in \Om$
\begin{eqnarray*}
\pi(\omega_1) \pi(\omega_2)
s(\omega_1,\omega_2;\omega^\prime_1,\omega^\prime_2)
=
\pi(\omega_1^\prime) \pi(\omega_2^\prime)
s(\omega_1^\prime,\omega_2^\prime;\omega_1,\omega_2).
\end{eqnarray*}

\setcounter{aux}{\value{enumi}}
\end{enumerate}

For a precise formulation of the infinitesimal generator on
$\Omn$ we first define the map
$\Theta_{j}^{\omega'\omega''} : \Omn \rightarrow \Omn$ for every
$\omega',\omega''\in \Om$,
$j \in \Tn$:
\begin{eqnarray*}
\left( \Theta_{j}^{\omega'\omega''} \uo \right)_{i} =
\left\{
\begin{array}{lcl}
\omega'\quad&\text{ if }& i=j
\\[4pt]
\omega''\quad&\text{ if }& i=j+1
\\[4pt]
\omega_i\quad&\text{ if }& i\not=j,j+1.
\end{array}
\right.
\end{eqnarray*}
The infinitesimal generators defined by these rates will be denoted:
\begin{eqnarray*}
&&
\Ln f(\uo)
=
\sum_{j\in \Tn} \sum_{\omega',\omega''\in \Om}
r(\omega_j,\omega_{j+1};\omega',\omega'')
\big(f(\Theta_{j}^{\omega'\omega''} \uo)-f(\uo)\big).
\\[5pt]
&&
\Kn f(\uo)
=
\sum_{j\in \Tn} \sum_{\omega',\omega''\in \Om}
s(\omega_j,\omega_{j+1};\omega',\omega'')
\big(f(\Theta_{j}^{\omega'\omega''} \uo)-f(\uo)\big).
\end{eqnarray*}
We denote by
$\Xn_t$  the Markov process on the state space
$\Omn$ with infinitesimal generator
$\Gn:=\lambda(n)\Ln + \kappa(n) \Kn$.   with speed-up factors
$\lambda(n)$ and
$\kappa(n)$ to be specified later

\medskip
\noindent
{\bf Remarks:}

\begin{enumerate}[(1)]

\item
Conditions
(\ref{cond:cons}) and
(\ref{cond:irred}) together imply   that
$\sum_j\eta_j$ and
$\sum_j\zeta_j$ are indeed the only conserved quantities of the
dynamics.

\item Condition (\ref{cond:lrsymm}) together with
(\ref{eq:mirroretazeta}) is implementation on a microscopic level
 of the left-right symmetry of the pde
(\ref{eq:pde}). Actually, our main result, Theorem \ref{thm:main}, is
valid without this assumption but some of the arguments would be more
technical.

\item
Condition
(\ref{cond:rstaci}) implies that the product measures
$\pin_{\tau,\theta}$ are indeed stationary for the dynamics defined by
the asymmetric rates
$r$. This is  seen by applying similar  computations to those of
\cite{balazs},
\cite{cocozza},
\cite{rezakhanlou} or
\cite{tothvalko2}. Mind that  this is
\emph{not} a detailed balance condition for the rates $r$.

\item Condition (\ref{cond:srevers}) is a straightforward detailed
balance condition. It implies that the product measures
$\pin_{\tau,\theta}$ are reversible for the dynamics defined by
the symmetric rates $s$.

\end{enumerate}

We will refer to  the measures
$\pin_{\tau,\theta}$ as the
\emph{canonical} measures. Since
$\sum_j \zeta_j$ and
$\sum_j \eta_j$  are conserved the  canonical measures on
$\Omn$ are
\emph{not} ergodic. The conditioned measures defined on
$\Omn_{N,Z}$ by:
\begin{eqnarray*}
\pin_{N,Z}(\uo)
:=
\pin_{\tau,\theta}
\big(\uo
\big|
\sum_{j\in\Tn} \eta_j=N,
\,
\sum_{j\in\Tn} \zeta_j=Z,
\,
\big) =
\frac { \pin_{\tau,\theta}(\uo) \ind\{\uo\in\Omn_{N,Z}\} }
      { \pin_{\tau,\theta}(\Omn_{N,Z}) }
\end{eqnarray*}
are also stationary and due to condition (\ref{cond:irred})
satisfied by the rate functions  they are ergodic. We shall call
these measures the \emph{microcanonical measures} of our system.
(It is easy to see that the measure $\pin_{N,Z}$ does not depend
on the choice of thew values of $\tau$ and $\theta$ in the
previous definition.)

The assumptions are by no means excessively restrictive.
Here follow some concrete examples of interacting particle systems
which belong to the class specified by conditions
(\ref{cond:cons})-(\ref{cond:srevers}) and also satisfy  the further
conditions
(\ref{cond:gradflux}),
(\ref{cond:psi1'}),
(\ref{cond:phi0'-psi1}),
(\ref{cond:lsi})
to be formulated later.

\begin{enumerate}[]

\item \underline{\emph{$\{-1,0,+1\}$-model}} The model is
described and analyzed in full detail in \cite{tothvalko2} and
\cite{fritztoth}. The one spin state space is $\Omega=\{-1,0,+1\}$
. The left-right reflection symmetry is implemented by $R:\Om\to
\Om$, $R\omega=-\omega$.  The dynamics consists of nearest
neighbor spin exchanges  and the two conserved quantities are
$\eta(\omega)=1-\abs{\omega}$ and $\zeta(\omega)=\omega$. The jump
rates are
\[
\begin{array}{rr}
 r(1,-1;-1,1)=0, & \qquad r(-1,1;1,-1)=2,
\\[5pt]
 r(0,-1;-1,0)=0, &\qquad r(-1,0;0,-1)=1,
\\[5pt]
 r(1,0,0,1)=0, & \qquad r(0,1,1,0)=1.
\end{array}
\]
and
\[
s(\omega_1,\omega_2;\omega'_1,\omega'_2)
= \left\{
\begin{array}{cl}
1
&
\text{ if } (\omega_1,\omega_2)=(\omega'_2,\omega'_1)
\text{ and } \omega_1\not=\omega_2
\\
0
&
\text{ otherwise}.
\end{array}
\right.
\]
The one dimensional marginals of the stationary measures are
\[
\pi_{\rho,u}(0)=\rho,
\quad
\pi_{\rho,u}(\pm1)=\frac{1-\rho\pm u}{2}
\]
with the domain of variables
${\cal D}=\{(\rho,u)\in\R_+\times\R: \rho+|u|\le1\}$.

\item \underline{\emph{Two-lane models}} The following family of
examples are  finite state space  versions of the  bricklayers
models introduced in \cite{tothwerner2}.  Let $\Om=\{0,1, \dots,
\bar n \}\times\{-\bar z,-\bar z+1,\dots,\bar z-1,\bar z\}$,
where $\bar n\in\N$ and $\bar z\in\{\frac12,1,\frac32,2\dots\}$.
The elements of $\Om$  will be denoted
$\omega:=\left(\substack{\eta\\\zeta}\right)$. Naturally enough,
$\sum_j\eta_j$ and $\sum_j\zeta_j$ will be the conserved
quantities of the dynamics. Left-right reflection symmetry is
implemented as $R:\Om\to \Om$,
$R\left(\substack{\eta\\\zeta}\right)=
\left(\substack{\eta\\-\zeta}\right)$. We
allow only the following elementary changes to occur at
neighboring sites $j,j+1$:
\[
\left(
\substack{\eta_{j}^{\phantom{a}}\\\zeta_{j}^{\phantom{a}}} ,
\substack{\eta_{j+1}^{\phantom{a}}\\\zeta_{j+1}^{\phantom{a}}}
\right)
\rightarrow
\left(
\substack{\eta_{j}^{\phantom{a}}\\\zeta_{j}^{\phantom{a}}\mp1} ,
\substack{\eta_{j+1}^{\phantom{a}}\\\zeta_{j+1}^{\phantom{a}}\pm1}
\right),
\qquad
\left(
\substack{\eta_{j}^{\phantom{a}}\\\zeta_{j}^{\phantom{a}}} ,
\substack{\eta_{j+1}^{\phantom{a}}\\\zeta_{j+1}^{\phantom{a}}}
\right)
\rightarrow
\left(
\substack{\eta_{j}^{\phantom{a}}\mp1\\\zeta_{j}^{\phantom{a}}} ,
\substack{\eta_{j+1}^{\phantom{a}}\pm1\\\zeta_{j+1}^{\phantom{a}}}
\right)
\]
with appropriate rates. Beside the conditions already imposed we also
assume that the one dimensional marginals of the steady state measures
factorize as follows:
\[
\pi(\omega)=
\pi\left(\substack{\eta\\\zeta}\right)
=
p(\eta)q(\zeta).
\]
The simplest case, with $\bar n=1$ and $\bar z = 1/2$, that is
with $\Omega=\{0,1\}\times\{-1/2,+1/2\}$, was introduced and fully
analyzed in \cite{tothvalko2} and  \cite{popkovschutz}. For a full
description (i.e. identification of the rates which satisfy the
imposed conditions, Eulerian hydrodynamic limit, etc. see those
papers.) It turns out that conditions
(\ref{cond:cons})-(\ref{cond:srevers}) impose some nontrivial
combinatorial constraints on the rates which are satisfied by a
finite parameter family of models. The number of free parameters
increases with $\bar n$ and $\bar z$. Since the concrete
expressions of the rates are not relevant for our further
presentation we omit the lengthy computations.

\end{enumerate}


\subsection{Expectations}
\label{subs:expectations}

Expectation, variance, covariance with respect to the measures
$\pin_{\tau,\theta}$ will be denoted by
$\expect_{\tau, \theta}(.)$, $\var_{\tau,\theta}(.)$,
$\cov_{\tau, \theta}(.)$.

We compute the expectations of the conserved quantities with
respect  to the  canonical measures, as functions of the
parameters
$\tau$ and
$\theta$:
\begin{eqnarray*}
&&
\rho(\tau,\theta)
:=
\expect_{\tau, \theta}(\eta)
=
\sum_{\omega \in \Om}
\eta(\omega) \pi_{\tau, \theta}(\omega)
=
G_{\tau}(\tau,\theta).
\\[2pt]
&&
u(\tau, \theta)
:=
\expect_{\tau, \theta}(\zeta)
=
\sum_{\omega \in \Om}
\zeta(\omega) \pi_{\tau, \theta}(\omega)
=
G_{\theta}(\tau, \theta),
\end{eqnarray*}
Elementary calculations show, that the matrix-valued function
\begin{eqnarray*}
\left(
\begin{array}{cc}
   \rho_{\tau} &  \rho_{\theta}
\\[4pt]
   u_{\tau} &  u_{\theta}
\end{array}
\right)
=
\left(
\begin{array}{cc}
   G_{\tau\tau} &  G_{\tau\theta}
\\[4pt]
   G_{\theta\tau}   &  G_{\theta\theta}
\end{array}
\right) =: G''(\tau, \theta)
\end{eqnarray*}
is equal to the covariance matrix $\cov_{\tau,\theta}(\eta,\zeta)$
and therefore it is strictly positive definite. It follows that
the function $(\tau,\theta)
\mapsto(\rho(\tau,\theta),u(\tau,\theta))$  is invertible. We
denote the inverse function by
$(\rho,u)\mapsto(\tau(\rho,u),\theta(\rho,u))$. Denote by
$(\rho,u)\mapsto S(\rho,u)$ the convex conjugate (Legendre
transform) of the strictly convex function $(\tau,\theta)\mapsto
G(\tau,\theta)$:
\begin{eqnarray}
\label{eq:thdentropy}
S(\rho,u)
:=
\sup_{\tau,\theta}
\big(
\rho\tau+ u\theta-G(\tau,\theta)
\big),
\end{eqnarray}
and
\beq
\label{eq:domain}
{\cal D}
&:=&
\{(\rho,u)\in\R_+\times\R:
S(\rho,u)<\infty\}
\\[5pt]
\notag
&=&
\text{co}\{(\eta,\zeta):\pi(\omega)>0\},
\eeq
where
$\text{co}$ stands for convex hull. The nondegeneracy condition
(\ref{eq:pinondeg})
implies that
$\partial{\cal D}\cap\{\rho=0\}=\{(0,u):|u|\le u_*\}$.  For
$(\rho,u)\in{\cal D}$ we have
\begin{eqnarray*}
\tau(\rho,u)=S_\rho(\rho,u),
\qquad
\theta(\rho,u)=S_u(\rho,u).
\end{eqnarray*}
In probabilistic terms:
$S(\rho,u)$ is the rate function of joint large deviations of
$(\sum_j\eta_j, \sum_j\zeta_j)$. In thermodynamic terms:
$S(\rho,u)$  corresponds to the equilibrium thermodynamic entropy. Let
\begin{eqnarray*}
\left(
\begin{array}{cc}
   \tau_{\rho} &  \tau_{u}
\\[4pt]
   \theta_{\rho} &     \theta_{u}
\end{array}
\right)
=
\left(
\begin{array}{cc}
   S_{\rho\rho} &  S_{\rho u}
\\[4pt]
   S_{u\rho}   &  S_{uu}
\end{array}
\right)
=:
S''(\rho,u).
\end{eqnarray*}
It is obvious that the matrices $G''(\tau,\theta)$ and
$S''(\rho,u)$ are strictly positive definite and are inverse of
each other:
\begin{eqnarray}
\label{eq:inverse}
G''(\tau,\theta)S''(\rho,u)
=
I
=
S''(\rho,u)G''(\tau,\theta),
\end{eqnarray}
where either
$(\tau,\theta)=(\tau(\rho,u),\theta(\rho,u))$  or
$(\rho,u)=(\rho(\tau,\theta),u(\tau,\theta))$.  With slight abuse of
notation we shall denote:
$\pi_{\tau(\rho,u),\theta(\rho,u)}=:\pi_{\rho,u}$,
$\pin_{\tau(\rho,u),\theta(\rho,u)}=:\pin_{\rho,u}$,
$\expect_{\tau(\rho,u),\theta(\rho,u)}=:\expect_{\rho,u}$,
etc.

As a general convention,  if
$\xi:\Om^m\to\R$  is a local function then its expectation  with
respect to the canonical measure
$\pi^{^{_{m}}}_{\rho,u}$ is denoted by
\begin{eqnarray*}
\Xi(\rho,u)
:=
\expect_{\rho,u} (\xi)
=
\sum_{\omega_1,\dots,\omega_m\in \Om^m}
\xi(\omega_1,\dots,\omega_m)
\pi_{\rho,u}(\omega_1)
\cdots
\pi_{\rho,u}(\omega_m).
\end{eqnarray*}

\subsection{Fluxes}
\label{subs:fluxes}

We introduce the fluxes of the conserved quantities. The
infinitesimal generators
$\Ln$ and
$\Kn$ act on the conserved quantities as follows:
\begin{eqnarray*}
\begin{array}{rll}
\Ln\eta_i=
&
- \psi(\omega_{i},\omega_{i+1}) +
\psi(\omega_{i-1},\omega_{i})
&
=:
-\psi_{i} + \psi_{i-1},
\\[7pt]
\Ln\zeta_i=
& - \phi(\omega_{i},\omega_{i+1}) +
\phi(\omega_{i-1},\omega_{i})
&
=:
-\phi_{i} + \phi_{i-1},
\\[7pt]
\Kn\eta_i=
&
- \psi^s(\omega_{i},\omega_{i+1}) +
\psi^s(\omega_{i-1},\omega_{i})
&
=:
-\psi^s_{i} + \psi^s_{i-1},
\\[7pt]
\Kn\zeta_i=
&
- \phi^s(\omega_{i},\omega_{i+1}) +
\phi^s(\omega_{i-1},\omega_{i})
&
=: -\phi^s_{i} + \phi^s_{i-1},
\end{array}
\end{eqnarray*}
where
\begin{eqnarray}
\label{eq:phipsidef}
&&
\begin{array}{rrl}
\psi(\omega_1,\omega_2)
&:=&
\dsum_{\omega'_1,\omega'_2\in \Om}
r(\omega_1,\omega_2;\omega'_1,\omega'_2)
\big(\eta(\omega'_2)-\eta(\omega_2)\big)
\\[15pt]
\phi(\omega_1,\omega_2)
&:=&
\dsum_{\omega'_1,\omega'_2\in \Om}
r(\omega_1,\omega_2;\omega'_1,\omega'_2)
\big(\zeta(\omega'_2)-\zeta(\omega_2)\big)
\end{array}
\\[5pt]
\label{eq:phispsisdef}
&&
\begin{array}{rrl}
\psi^s(\omega_1,\omega_2)
&:=&
\dsum_{\omega'_1,\omega'_2\in \Om}
s(\omega_1,\omega_2;\omega'_1,\omega'_2)
\big(\eta(\omega'_2)-\eta(\omega_2)\big)
\\[15pt]
\phi^s(\omega_1,\omega_2)
&:=&
\dsum_{\omega'_1,\omega'_2\in \Om}
s(\omega_1,\omega_2;\omega'_1,\omega'_2)
\big(\zeta(\omega'_2)-\zeta(\omega_2)\big)
\end{array}
\end{eqnarray}
Note that due to the left-right symmetry and conservations, i.e.
(\ref{eq:mirroretazeta}) and   conditions
(\ref{cond:cons}) and (\ref{cond:lrsymm}), the microscopic fluxes have
the following symmetries:
\beqs
\phi(\omega_1,\omega_2)=\phantom{-}\phi(R\omega_2,R\omega_1),
\\[5pt]
\psi(\omega_1,\omega_2)=-\psi(R\omega_2,R\omega_1).
\eeqs

In order to simplify some of our  further arguments (in
particular, see (\ref{eq:kns}) in subsection
\ref{subs:firsttermright1}) we impose one more microscopic
condition

\begin{enumerate}[{(A)}]

\setcounter{enumi}{\value{aux}}
\item
\label{cond:gradflux}
\emph{Gradient condition on symmetric fluxes:}
The microscopic fluxes of the symmetric part, defined in
(\ref{eq:phispsisdef})  satisfy the following gradient conditions
\beq
\label{eq:kappa}
&&
\psi^s(\omega_1,\omega_2)
=
\kap(\ome_1)-\kap(\ome_2)
=:
\kap_1-\kap_2
\\[5pt]
\notag
&&
\phi^s(\omega_1,\omega_2)
=
\chi(\ome_1)-\chi(\ome_2)
=:
\chi_1-\chi_2.
\eeq

\setcounter{aux}{\value{enumi}}
\end{enumerate}

\noindent
{\bf Remark:}
(1)
This is a technical assumption
(referring actually to the measure $\pi$)
which  simplifies considerably the
arguments of subsection \ref{subs:lentboundsproof}. The symmetric
part $\Kn$ has the role of enhancing convergence to local
equilibrium. Its effect is
\emph{not  seen} in the limit, so in principle we can choose it
conveniently.   Without this assumption we would be forced  to use all
the non-gradient  technology developed in
\cite{varadhan} (see also
\cite{kipnislandim}), which would make the paper even longer.
\\
(2) It is easy to see that $\eta(\omega_1)=\eta(\omega_2)=0$
implies $\psi^s(\omega_1,\omega_2)=0$ and thus (by choosing a
suitable additive constant) $\omega\mapsto\kappa(\omega)$ can be
chosen so that \beq \label{eq:kappavanish} \eta(\omega)=0
\,\,\,\Rightarrow\,\,\, \kappa(\omega)=0. \eeq

The
\emph{macroscopic fluxes} are:
\begin{eqnarray}
\label{eq:PhiPsidef}
\begin{array}{rrl}
\Psi(\rho,u)
&:=&
\expect_{\rho,u} (\psi)
\\[8pt]
&=&
\hskip-1mm
\dsum_{\substack{\omega_1,\omega_2,\\\omega'_1,\omega'_2}\in \Om}
r(\omega_1,\omega_2;\omega'_1,\omega'_2)
\big(\eta(\omega'_2)-\eta(\omega_2)\big)
\pi_{\rho,u}(\omega_1)
\pi_{\rho,u}(\omega_2),
\\[25pt]
\Phi(\rho,u)
&:=&
\expect_{\rho,u} (\phi)
\\[8pt]
&=&
\hskip-1mm
\dsum_{\substack{\omega_1,\omega_2,\\\omega'_1,\omega'_2}\in \Om}
r(\omega_1,\omega_2;\omega'_1,\omega'_2)
\big(\zeta(\omega'_2)-\zeta(\omega_2)\big)
\pi_{\rho,u}(\omega_1)
\pi_{\rho,u}(\omega_2).
\end{array}
\end{eqnarray}
These are smooth regular functions of the variables
$(\rho,u)\in{\cal D}$. Note that due to reversibility of
$\Kn$, for any value of
$\rho$ and
$u$
\begin{eqnarray*}
\expect_{\rho,u} (\psi^s)
=0=
\expect_{\rho,u} (\phi^s).
\end{eqnarray*}
(These identities hold true without assuming condition
(\ref{cond:gradflux}).)

For later use we mention here that according to
\cite{tothvalko2}, the  macroscopic fluxes
$\Psi(\rho,u)$ and
$\Phi(\rho,u)$ satisfy the following \emph{Onsager reciprocity
  relation}
\beq
\label{eq:onsager}
&&
\Psi_u(\rho,u) \var_{\rho,u}(\zeta)
-
\Phi_u(\rho,u) \cov_{\rho, u}(\eta,\zeta)
=
\\[5pt]
\notag
&&
\hskip5cm
\Phi_\rho(\rho,u) \var_{\rho,u}(\eta)
-
\Psi_\rho(\rho,u) \cov_{\rho, u}(\eta,\zeta).
\eeq

For the concrete examples presented at the end of subsection
\ref{subs:rates} the following domains $\dom$ and macroscopic rates
are gotten:

\begin{enumerate}[]

\item
\underline{\emph{$\{-1,0,+1\}$-model}}:
\beqs
&&
\dom=\{(\rho,u)\in\R_+\times\R\,:\,\rho+\abs{u}\le1\}
\\
&&
\Psi(\rho,u)=\rho u
\\
&&
\Phi(\rho,u)=\rho+u^2.
\eeqs

\item
\underline{\emph{Two lane models}} with $\bar n=1$:
\beqs
&&
\dom=\{(\rho,u)\in\R_+\times\R\,:\,\rho\le 1,\,\,\,\abs{u}\le\bar z\}
\\
&&
\Psi(\rho,u)=\rho(1-\rho)\psi(u)
\\
&&
\Phi(\rho,u)=\varphi_0(u)+\rho\varphi_1(u),
\eeqs
where $\psi(u)$ is odd, while $\varphi_{0}(u)$ and
$\varphi_{1}(u)$ are even functions of
$u$, determined by the jump rates of the model. In the simplest
particular case with $\bar z=1/2$
\beqs
&&
\Psi(\rho,u)=\rho(1-\rho)u
\\
&&
\Phi(\rho,u)=(\rho-\gamma)(1-u^2),
\eeqs
where $\gamma\in\R$ is the only model dependent parameter which
appears in the macroscopic fluxes. For details see \cite{tothvalko2}.
\end{enumerate}



\subsection{The hdl under Eulerian scaling}
\label{subs:hdl}

Given a system of interacting particles as defined in the previous
subsections, by  applying Yau's relative entropy method (see
\cite{yau1} or the monograph
\cite{kipnislandim}),  one shows that under Eulerian scaling the local
densities of the conserved quantities
$\rho(t,x), \, u(t,x)$ evolve according to the system of partial
differential equations:
\begin{eqnarray}
\label{eq:euler}
\left\{
\begin{array}{l}
  \pt \rho+\px \Psi(\rho, u)=0
\\[6pt]
  \pt u+\px \Phi(\rho,u)=0
\end{array}
\right.
\end{eqnarray}
where
$\Psi(\rho, u)$ and
$\Phi(\rho,u)$ are the macroscopic fluxes defined in
(\ref{eq:PhiPsidef}).

The precise statement of the hydrodynamical limit is as follows:
Consider a microscopic system which satisfies conditions
(\ref{cond:cons})-(\ref{cond:srevers}) of subsection
\ref{subs:rates}. Note that condition (\ref{cond:gradflux}) of
subsection \ref{subs:fluxes} is not assumed. Let $\Psi(\rho,u)$
and $\Phi(\rho,u)$ be the macroscopic fluxes computed for this
system and $\rho(t,x), u(t,x)$ $x\in\T$, $t\in[0,T]$ be
\emph{smooth} solution of the pde (\ref{eq:euler}). Let the
microscopic system of size $n$ be driven by the infinitesimal
generator
\begin{eqnarray*}
\Gn= n\Ln+ n^{1+\delta}\Kn,
\end{eqnarray*}
where
$\delta\in[0,1)$ is fixed. This means that the main, asymmetric  part
  of the   generator is speeded up by
$n$ and the additional symmetric part by
$n^{1+\delta}$. Let
$\mun_0$ be a probability distribution on
$\Omn$ which is the initial distribution of the microscopic system of
  size
$n$, and
\begin{eqnarray*}
\mun_t
:=
\mun_0 e^{t\Gn}
\end{eqnarray*}
the distribution of the system at (macroscopic) time
$t$. The
\emph{local equilibrium}  measure
$\nun_t$ (itself a probability measure on
$\Omn$) is defined by
\begin{eqnarray*}
\nun_t
:=
\prod_{j\in\Tn} \pi_{ \rho(t,\frac{j}{n}), u(t,\frac{j}{n})}.
\end{eqnarray*}
This measure \emph{mimics on a microscopic scale}  the macroscopic
evolution driven by the pde (\ref{eq:euler}).

We denote by $H(\mun_t|\pin)$, respectively, by $H(\mun_t|\nun_t)$
the relative entropy of the measure $\mun_t$ with respect to the
absolute reference measure $\pin$, respectively, with respect to
the local equilibrium measure $\nun_t$.

The precise statement of the Eulerian hydrodynamic limit is the
following

\begin{theorem*}
Assume conditions
(\ref{cond:cons})-(\ref{cond:srevers}) and let
$\delta\in[0,1)$ be fixed. If
\begin{eqnarray*}
H \left( \mun_0 \,\left|\, \nun_0
\right.\right) =o(n)
\end{eqnarray*}
then
\begin{eqnarray*}
H\left( \mun_t \,\left|\,\nun_t \right.\right) =o(n)
\end{eqnarray*}
uniformly for
$t\in[0,T]$.
\end{theorem*}

\noindent
{\bf Remark:}
Note that due to finiteness of the state space
$\Om$ the condition
\begin{eqnarray*}
H\left(\mun_0\,\left|\,\pin_{1,0}\right.\right)
=\Ordo(n)
\end{eqnarray*}
holds automatically.

The Theorem follows  from direct application of Yau's relative entropy
method. For the proof and its direct consequences see
\cite{yau1},
\cite{kipnislandim} or
\cite{tothvalko2}. For the main consequences of this Theorem see e.g.
Corollary 1 of \cite{tothvalko2}.



\section{Low density asymptotics and the main result:
hydrodynamic limit under intermediate scaling}
\label{section:intermediate}

\subsection{General properties and low density asymptotics of the
  macroscopic fluxes}
\label{subs:genprop}
\label{subs:lda}

The fluxes in the Euler equation (\ref{eq:euler}) are regular
smooth functions of in $(\rho,u)\in\overline{{\cal D}}$.

From the left-right symmetry of the microscopic models it follows that
\beq
\label{eq:parity}
\Phi(\rho,-u)=\Phi(\rho,u),
\qquad
\Psi(\rho,-u)=-\Psi(\rho,u).
\eeq
It is also obvious that for
$u\in[-u_*,u_*]$
\beq
\label{eq:psizero}
\Psi(0,u)=0.
\eeq

We make two assumptions about the low density asymptotics of the
macroscopic fluxes. Here is the first one:

\begin{enumerate} [(A)]
\setcounter{enumi}{\value{aux}}

\item
\label{cond:psi1'}
We assume that
$\Psi_{\rho u}(0,0)\not=0$.  Actually, by possibly redefining the time
scale and orientation of space, without loss of generality we assume
\beq
\label{eq:psi1'}
\Psi_{\rho u}(0,0)=1.
\eeq

\setcounter{aux}{\value{enumi}}
\end{enumerate}

From the Onsager relation
(\ref{eq:onsager}) and obvious parity  considerations  it also follows
that
\beq
\label{eq:phi1}
\Phi_{\rho}(0,0)
=
\Psi_{\rho u}(0,0)\var_{0,0}(\zeta)
=
1.
\eeq
Note, that here we rely on the choice
(\ref{eq:choiceofv0}) of the scaling factor
$v_0$ in
(\ref{eq:zetamap}).

We denote
\beq
\label{eq:gammadef}
\gamma:=\frac12\Phi_{uu}(0,0).
\eeq
Our results will hold for
$\gamma>1$ only.

From
(\ref{eq:parity}) and
(\ref{eq:psi1'}) it follows that
\beq
\label{eq:phi0'-psi1asym}
\Phi_{u}(0,u)-\Psi_{\rho}(0,u)
=
(2\gamma-1) u + \Ordo(|u|^3).
\eeq
The second condition imposed on the low density asymptotics of the
macroscopic fluxes is:

\begin{enumerate}[(A)]

\setcounter{enumi}{\value{aux}}

\item
\label{cond:phi0'-psi1}
For $u\in[-u_*,u_*]$, $u\not=0$ %
\beq \label{eq:phi0'-psi1} \Phi_{u}(0,u)-\Psi_{\rho}(0,u) \not=0,
\\
\Phi_{\rho} (0,u)\neq 0, \quad \Psi_{\rho u} (0,u)\neq
0\label{eq:phiro_psiru}
\eeq%

\setcounter{aux}{\value{enumi}}

\end{enumerate}

\noindent{\bf Remarks:}
(1)
(\ref{cond:psi1'})
is a very natural nondegeneracy condition: if
$\Psi_{\rho  u}(0,0)$  vanished then in the perturbation calculus to
be performed, higher order terms would be dominant and a different
scaling limit should be taken.
\\
(2) Due to (\ref{eq:parity}), (\ref{eq:psi1'}) and
(\ref{eq:phi0'-psi1asym}) conditions (\ref{eq:phi0'-psi1}),
(\ref{eq:phiro_psiru})  hold anyway in a neighborhood of $u=0$,
and this would suffice, but the forthcoming arguments, in
particular the proof of Lemmas \ref{lemma:charcurves} and
\ref{lemma:lentbounds} would be less transparent. We assume
condition (\ref{cond:phi0'-psi1})  for technical convenience only.
Condition (\ref{eq:phi0'-psi1}) amounts to forbidding other 
non-hyperbolic pointson 
$\overline{\partial{\cal D}\cap\{\rho=0\}}$, beside the point
$(\rho,u)=(0,0)$. Condition (\ref{eq:phiro_psiru}) reflects 
the natural monotonicity requirements (i) and (ii) formulated
about the microscopic models
at the beginning of Section \ref{section:models}.

We are interested in the behavior of the pde near the isolated
non-hyperbolic point $(\rho,u)=(0,0)$. The asymptotic expansion
for $\rho+u^2\ll 1$ of the macroscopic fluxes and their first
partial derivatives is
\begin{eqnarray}
\begin{array}{ll}
\Psi(\rho,u)
=
\rho u \big( 1 + \Ordo(\rho+u^2) \big),
&
\Phi(\rho,u)
=
(\rho  + \gamma u^2) \big( 1 + \Ordo(\rho+u^2) \big),
\\[8pt]
\Psi_\rho(\rho,u)
=
u \big( 1 + \Ordo(\rho+u^2) \big),
&
\Phi_\rho(\rho,u)
=
1 + \Ordo(\rho+u^2),
\\[8pt]
\Psi_u(\rho,u)
=
\rho  \big( 1 + \Ordo(\rho+u^2) \big),
&
\Phi_u(\rho,u)
=
2\gamma u \big( 1 + \Ordo(\rho+u^2) \big).
\end{array}\label{eq:macr_flx_asy}
\end{eqnarray}
We are looking for ``small solutions'' of the pde
(\ref{eq:euler}):  Let $\rho_0(x)$ and $u_0(x)$ be given profiles
and assume that $\rho^{\vareps}(t,x)$, $u^{\vareps}(t,x)$ is
solution of  the pde (\ref{eq:euler}) with initial condition
\begin{eqnarray*}%
\rho^{\vareps}(0,x)=\vareps^{2}\rho_0(x), \quad
u^{\vareps}(0,x)=\vareps \,u_0(x).
\end{eqnarray*}
Then, at least formally,
\begin{eqnarray*}
\vareps^{-2}\rho^{\vareps}(\vareps^{-1}t,x) \to \rho(t,x), \quad
\vareps^{-1}u^{\vareps}(\vareps^{-1}t,x) \to u(t,x),
\end{eqnarray*}
where $\rho(t,x)$, $u(t,x)$ is solution of the pde (\ref{eq:pde})
with initial condition
\begin{eqnarray*}
\rho(0,x)=\rho_0(x), \quad u(0,x)=u_0(x).
\end{eqnarray*}

\subsection{The main result}
\label{subs:mainres}

The asymptotic computations of subsection \ref{subs:lda} suggest
the scaling under which we  should derive the pde (\ref{eq:pde})
as hydrodynamic limit: fix a (small) positive $\beta$ and choose
the scaling
\begin{eqnarray*}
\begin{array}{lccc}
                           & \text{MICRO}      & \phantom{MMMM}
                                                & \text{MACRO}
\\[8pt]
\text{space}               &   nx              &&       x
\\[5pt]
\text{time}                &   n^{1+\beta}t    &&       t
\\[5pt]
\text{particle density}
      \phantom{MMMM}       &   n^{-2\beta}\rho &&       \rho
\\[5pt]
\text{`slope of the wall'} &   n^{-\beta}u     &&        u
\end{array}
\end{eqnarray*}
Ideally the result should be valid for
$0<\beta<1/2$  but we are able to prove much less than that.

Choose a model satisfying the conditions
(\ref{cond:cons})-(\ref{cond:gradflux}) of section
\ref{section:models} and conditions
(\ref{cond:psi1'}-\ref{cond:phi0'-psi1}) of subsection
\ref{subs:genprop}, and let
$\gamma$ be given by
(\ref{eq:gammadef}), corresponding to the microscopic system chosen.
Let the microscopic system of size
$n$ (defined on the discrete torus
$\Tn$) evolve (on macroscopic time scale) according to the
infinitesimal generator
\begin{eqnarray*}
\Gn=
n^{1+\beta}\Ln + n^{1+\beta+\delta}\Kn.
\end{eqnarray*}
with $\beta>0$ and some further conditions to be imposed on
$\beta$ and $\delta$ (see Theorem \ref{thm:main}).  Denote by
$\mun_t$ the true distribution of the microscopic system at
macroscopic time $t$:
\begin{eqnarray*}
\mun_t
:=
\mun_0 e^{t\Gn},
\end{eqnarray*}
where
$\mun_0$ is the initial distribution.

We use the translation invariant product measure
\begin{eqnarray*}
\pin
:=
\pin_{n^{-2\beta},0}
\end{eqnarray*}
as
\emph{absolute reference measure}.
Global entropy will be considered relative to this
measure, Radon-Nikodym derivatives of
$\mun_t$  and the local equilibrium measure
$\nun_t$ to be defined below, with respect to
$\pin$ will be used.

Given a smooth solution $\big(\rho(t,x),u(t,x)\big)$,
$(t,x)\in[0,T]\times\T$, of the pde (\ref{eq:pde}) define the
\emph{local equilibrium measure} $\nun_t$ on $\Omn$ as follows
\begin{eqnarray}
\label{eq:loceqm}
\nun_t
:=
\prod_{j\in\Tn}
\pi_{n^{-2\beta}  \rho(t,\frac{j}{n}), n^{-\beta} u(t,\frac{j}{n})}.
\end{eqnarray}
This time-dependent measure \emph{mimics on a microscopic level}
the macroscopic evolution governed by the pde (\ref{eq:pde}).

Our main result is the following

\begin{theorem}
\label{thm:main}
Assume that the microscopic system of interacting particles satisfies
conditions
(\ref{cond:cons})-(\ref{cond:gradflux}) of subsubsections
\ref{subs:rates}, \ref{subs:fluxes}  and the uniform log-Sobolev
condition
(\ref{cond:lsi}) of subsection
\ref{subs:apriori}. Additionally, assume that the macroscopic fluxes
satisfy conditions
(\ref{cond:psi1'}),
(\ref{cond:phi0'-psi1}) of subsection
\ref{subs:genprop} and
$\gamma>1$. Choose
$\beta\in(0,1/2)$ and
$\delta\in(1/2,1)$ so that
\begin{eqnarray}
\label{eq:bdcond}
2\delta-8\beta>1
\quad
\text{ and }
\quad
\delta+3\beta<1.
\end{eqnarray}

Let $\big(\rho(t,x),u(t,x)\big)$, $(t,x)\in[0,T]\times\T$, be
\emph{smooth} solution of the pde (\ref{eq:pde}),  such that
$\inf_{x\in\T}\rho(0,x)>0$ and let $\nun_t$, $t\in[0,T]$
be the corresponding local equilibrium measure defined in
(\ref{eq:loceqm}).

Under these conditions, if
%
%
%
\begin{eqnarray}
\label{eq:main0}
H\big(\mun_0 \,\big|\,\nun_0 \big) =o(n^{1-2\beta})
\end{eqnarray}
then
\begin{eqnarray}
\label{eq:main}
H\big(\mun_t \,\big|\,\nun_t \big)=o(n^{1-2\beta})
\end{eqnarray}
uniformly for
$t\in[0,T]$.
\end{theorem}

\noindent
{\bf Remarks:}
\\
(i)
From (\ref{eq:main0}) via the identity (\ref{eq:id0})
and the entropy inequality
it also follows that
\begin{eqnarray}
\label{eq:ent0}
H\big(\mun_0\,\big|\,\pin\big) = \Ordo(n^{1-2\beta}).
\end{eqnarray}
See the beginning of subsection \ref{subs:prepcomp}
\\
(i)
If $\gamma>3/4$, in smooth solutions
vacuum does not appear. That is $\inf_{x\in\T}\rho(0,x)>0$ implies
$\inf_{(t,x)\in[0,T]\times\T}\rho(t,x)>0$.
\\
(ii)
Although for the \emph{$\{-1,0,+1\}$-model} we have $\g=1$,
our proof can also be extended to cover this model. Actually, in
that case the proof is much simpler, since the Eulerian pde is
equal to the limit pde (\ref{eq:pde}) and thus the cutoff function
(see Section \ref{section:cutoff}) can be determined explicitly.

\begin{corollary}
\label{cor:main}
Assume the conditions of Theorem \ref{thm:main}.
Let $g:\T\to\R$ be a test function. Then
for any $t\in[0,T]$
\\
(i)
\beqs
n^{2\b-1}
\sum_{j\in\Tn} g(\frac{j}{n})\eta_j(t)
&\buildrel\prob\over\rightarrow&
\int_\T g(x) \rho(t,x)\,dx,
\\[5pt]
n^{\b-1}
\sum_{j\in\Tn} g(\frac{j}{n})\zeta_j(t)
&\buildrel\prob\over\rightarrow&
\int_\T g(x) u(t,x)\,dx.
\eeqs
(ii)
\beqs
H\big(\mun_0\,\big|\,\pin\big)
-
H\big(\mun_t\,\big|\,\pin\big)
=
o(n^{1-2\beta}).
\eeqs

\end{corollary}

See the proof of Corollary 1 in \cite{tothvalko2}.
\section{Notations and general preparatory  computations}
\label{section:proofprep}

This section completely standard in the context of the relative
entropy method. So we shall be sketchy.

\subsection{Notation}
\label{subs:notation}

We denote
\begin{eqnarray*}
\hn(t)
&:=&
n^{-(1-2\beta)}
H\big( \mun_t \,\big|\,\nun_t \big).
\\[5pt]
\sn(t)
&:=&
n^{-(1-2\beta)}
\left(
H\big( \mun_0 \,\big|\,\pin \big)
-
H\big( \mun_t \,\big|\,\pin \big)
\right).
\end{eqnarray*}
We know
\emph{a priori}  that
$t\mapsto\sn(t)$ is monotone increasing and due to
(\ref{eq:ent0})
\begin{eqnarray}
\label{eq:apriorientropybound}
\sn(t)
=
\Ordo(1),
\quad
\text{ uniformly for }
t\in[0,\infty).
\end{eqnarray}
In fact, from Theorem \ref{thm:main}  it follows (see Corollary
\ref{cor:main}) that as long as the solution $\rho(t,x), u(t,x)$
of the pde (\ref{eq:pde}) is smooth
\begin{eqnarray*}
\sn(t)
=
o(1),
\quad
\text{ uniformly for }
t\in[0,T].
\end{eqnarray*}

For
$(\rho,u)\in (0,\infty)\times(-\infty,\infty)$ denote
\begin{eqnarray*}
&&
\taun(\rho,u)
:=
\tau(n^{-2\beta}\rho, n^{-\beta} u)-\tau(n^{-2\beta}, 0)
\\[5pt]
&&
\thetan(\rho,u)
:=
n^\beta\theta(n^{-2\beta}\rho, n^{-\beta} u).
\end{eqnarray*}
Note that, for symmetry reasons $\theta(n^{-2\beta}, 0)=0$.
Mind that
$\tau$ is chemical potential rather than fugacity and for small
densities the fugacity
$\lambda:=e^\tau$ scales like
$\rho$, i.e. $\tau(n^{-2\beta}, 0)\sim -2\beta\log n$. If
$\rho>0$ and
$u\in\R$ are fixed then
$\taun(\rho,u)$ and
$\thetan(\rho,u)$ stay of order
$1$, as
$n\to\infty$.

Given the smooth solution
$\rho(t,x), u(t,x)$, with
$\rho(t,x)>0$ we shall use the notation
\begin{eqnarray*}
&&
\taun(t,x):=\taun(\rho(t,x), u(t,x)),
\\[5pt]
&&
\thetan(t,x):=\thetan(\rho(t,x), u(t,x)),
\\[5pt]
&&
v(t,x):=\log\rho(t,x).
\end{eqnarray*}
The following asymptotics hold uniformly in
$(t,x)\in[0,T]\times\T$:
\beq
\label{eq:asym}
\begin{array}{rl}
\taun(t,x)=
\phantom{\px}
v(t,x)+\Ordo(n^{-2\b}),
&
\phantom{\px}
\thetan(t,x)
=
\phantom{\px}
u(t,x)+\Ordo(n^{-2\b})
\\[5pt]
\px\taun(t,x)
=
\px v(t,x)+\Ordo(n^{-2\b}),
&
\px\thetan(t,x)
=
\px u(t,x)+\Ordo(n^{-2\b})
\\[5pt]
\pt\,\taun(t,x)
=
\pt \,v(t,x)+\Ordo(n^{-2\b}),
&
\pt\,\thetan(t,x)
=
\pt \,u(t,x)+\Ordo(n^{-2\b})
\end{array}
\eeq

The logarithm of the Radom-Nikodym derivative of the  time dependent
reference measure
$\nunt$ with respect to the absolute referencee measure
$\pin$ is denoted by
$f^n_t$:
\begin{eqnarray}
\nonumber
\fnt(\uome)
&:=&
\log \frac{d\nun_t}{d\pin}(\uo)
\\[5pt]
\label{eq:RadNikdef}
&=&
\sum_{j\in \Tn}
\Big\{
\taun(t,\frac{j}{n}) \eta_j
+
n^{-\beta}\thetan(t,\frac{j}{n}) \zeta_j
\\
\nonumber
&&
\hskip10mm
-
G(\taun(t,\frac{j}{n})+ \tau(n^{-2\b},0) ,
n^{-\beta}\thetan(t,\frac{j}{n}))
+
G(\tau(n^{-2\b},0) ,0)
\Big\}
\end{eqnarray}
%

\subsection{Preparatory computations}
\label{subs:prepcomp}

In order to obtain the main estimate
(\ref{eq:main}) our aim is to get a Gr\"{o}nwall type inequality:
we will prove that for every
$t\in[0,T]$
\begin{eqnarray}
\label{eq:gronwall}
\hn(t)- \hn(0)
=
\int_0^t \pt\hn(s)ds
\leq
C \int_0^t \hn(s)ds +
o(1),
\end{eqnarray}
where the error term is uniform in
$t\in[0,T]$. Because it is assumed that
$\hn(0)=o(1)$, the Theorem follows.

We start with  the identity
\begin{eqnarray}
\label{eq:id0}
H(\munt|\nunt)-H(\munt|\pin)
=
-
\int_{\Omn} \fnt \,d\munt.
\end{eqnarray}
From this identity, the explicit form of the Radon-Nikodym
derivative (\ref{eq:RadNikdef}), the asymptotics (\ref{eq:asym}),
via the entropy inequality and (\ref{eq:main0})
the a priori entropy bound (\ref{eq:ent0}) follows indeed, as
remarked after the formulation
of Theorem \ref{thm:main}.

Next we differentiate (\ref{eq:id0}) to obtain
\begin{eqnarray}
\notag
\pt\hn(t)
=
-
\int_{\Omn}
\left(
n^{3\b}
{\Ln\!\fnt}
+
n^{3\b+\d}
{\Kn\!\fnt}
+
n^{-1+2\b}
{\pt\fnt}
\right)
\,d\munt
-\pt \sn(t).%
\\
\label{eq:entprod}
\end{eqnarray}
Usually, an adjoint version  of (\ref{eq:entprod}) is being used
in form of an inequality.  In our case this form is needed.
We emphasize that the term $-\pt\sn(t)$ on the right hand side
will be of crucial importance.

We compute the three terms under the integral.
\begin{eqnarray}
\label{eq:elef}
n^{3\b}
\Ln\!\fnt(\uome)
&=&
\phantom{+}
\frac{1}{n}\sum_{j\in\Tn}
\px v(t,\frac{j}{n})
n^{3\b}\psi_j
+
\frac{1}{n}\sum_{j\in\Tn}
\px u(t,\frac{j}{n})
n^{2\b}\phi_j
\\[5pt]
\notag
&&
+
A^n_1(t,\uome)
+
A^n_2(t,\uome)
+
A^n_3(t,\uome)
+
A^n_4(t,\uome),
\end{eqnarray}
where
$A^n_i(t,\uome)$,
$i=1,\dots,4$ are error terms which will be  easy to estimate:
\begin{eqnarray*}
A^n_1(t,\uome)
&:=&
\frac{1}{n}\sum_{j\in\Tn}
\big(\px \taun(t,\frac{j}{n})-
\px v(t,\frac{j}{n})\big)n^{3\b}\psi_j,
\\[5pt]
A^n_2(t,\uome)
&:=&
\frac{1}{n}\sum_{j\in\Tn}
\big(\px \thetan(t,\frac{j}{n})-
\px u(t,\frac{j}{n})\big)n^{2\b}\phi_j,
\\[5pt]
A^n_3(t,\uome)
&:=&
\frac{1}{n}\sum_{j\in\Tn}
\big(\gradn\taun(t,\frac{j}{n})-
\px \taun(t,\frac{j}{n})\big)n^{3\b}\psi_j
\\[5pt]
A^n_4(t,\uome)
&:=&
\frac{1}{n}\sum_{j\in\Tn}
\big(\gradn\thetan(t,\frac{j}{n})-
\px\thetan(t,\frac{j}{n})\big)n^{2\b}\phi_j.
\end{eqnarray*}
Here and in the sequel $\gradn$ denotes the discrete gradient:
\[
\gradn f(x):=n\big(f(x+1/n)-f(x)\big).
\]
See subsection
\ref{subs:firsterrors}  for the estimate of the error terms
$A^n_j(t,\uome)$,
$j=1,\dots,12$.

Next,
\begin{eqnarray}
\notag
n^{3\b+\d}
\Kn\!\fnt(\uome)
&=&
n^{-1+3\b+\d}
\frac{1}{n}\sum_{j\in\Tn}
\big(
\left(\gradn\right)^2 \taun(t,\frac{j}{n})\kappa_j
+
\left(\gradn\right)^2  \thetan(t,\frac{j}{n})\chi_j
\big)
\\[5pt]
\label{eq:kaef}
&=:&
A^n_5(t,\uome)
\end{eqnarray}
is itself a numerical error term. Finally
\begin{eqnarray}
\label{eq:peteef}
n^{-1+2\b}
{\pt\fnt}(\uome)
&=&
\phantom{+}
\frac{1}{n}\sum_{j\in\Tn}
\pt v(t,\frac{j}{n})
\big(n^{2\b}\eta_j-\rho(t,\frac{j}{n})\big)
\\[5pt]
\notag
&&
+
\frac{1}{n}\sum_{j\in\Tn}
\pt u(t,\frac{j}{n})
\big(n^{\b}\zeta_j-u(t,\frac{j}{n})\big)
\\[5pt]
&&
\notag
+
A^n_6(t,\uome)
+
A^n_7(t,\uome),
\end{eqnarray}
where
\begin{eqnarray*}
A^n_6(t,\uome)
&:=&
\frac{1}{n}\sum_{j\in\Tn}
\big(\pt\tau^n(t,\frac{j}{n})-\pt v(t,\frac{j}{n})\big)
\big(n^{2\b}\eta_j-\rho(t,\frac{j}{n})\big),
\\[5pt]
A^n_7(t,\uome)
&:=&
\frac{1}{n}\sum_{j\in\Tn}
\big(\pt\theta^n(t,\frac{j}{n})-\pt u(t,\frac{j}{n})\big)
\big(n^{\b}\zeta_j-u(t,\frac{j}{n})\big).
\end{eqnarray*}
are again easy-to-estimate error terms.

\subsection{Blocks}
\label{subs:blocks}

We fix once and for all a weight function $a:\R\to\R$. It is
assumed that:
\\
(1) $a(x)>0$ for $x\in \left(-1,1\right)$ and $a(x)=0$ otherwise,
\\
(2)
it has total weight
$\int a(x)\,dx=1$,
\\
(3)
it is even:
$a(-x)=a(x)$, and
\\
(4)
it is twice continuously differentiable.

We choose a
\emph{mesoscopic} block size
$l=l(n)$ such that
\beq
\label{eq:blocksizebounds}
1\ll n^{(1+\delta+5\beta)/3}\ll l(n) \ll n^{\delta-\beta}\ll n.
\eeq
This can be done due to condition (\ref{eq:bdcond}) imposed on
$\b$ and $\d$. 

Given a local variable (depending on
$m$ consecutive spins)
\[
\xi_i
=
\xi_i(\uome)
=
\xi(\ome_i,\dots,\ome_{i+m-1}),
\]
its \emph{block average} \emph{at macroscopic space coordinate}
$x$ is defined as
\begin{eqnarray}
\label{eq:blocks}
\wih{\xi}(x)
=
\wih{\xi}(\uo,x)
:=
\frac{1}{l}\sum_ja\left(\frac{nx-j}{l}\right)\xi _j.
\end{eqnarray}
Since
$l=l(n)$, we do not denote explicitly dependence of the block average
on the mesoscopic block size
$l$.

Note that
$x\mapsto\wih{\xi}(x)$ is smooth
\[
\px\wih{\xi}(x)
=
\px\wih{\xi}(\uome,x)
=
\frac{n}{l}\,
\frac{1}{l}\sum_j a{'}\left(\frac{nx-j}{l}\right)\xi _j,
\]
and  it is straightforward that
\begin{eqnarray}
\label{eq:trivigradbound}
\sup_{\uome\in\Omn}\sup_{x\in\T}
\abs{\px\wih{\xi}(\uome,x)}
\le
C \,
\left(\max_{\ome_1,\dots\ome_{m}}\xi(\ome_1,\dots,\ome_m)\right)\,
\frac{n}{l}.
\end{eqnarray}
For a more sophisticated bound on
$\abs{\px\wih{\xi}(\uome,x)}$ see
(\ref{eq:plaingrad}).

We shall use the handy (but slightly abused) notation
\begin{eqnarray*}
\wih{\xi}(t,x):=\wih{\xi}(\Xn_t,x).
\end{eqnarray*}
This is the empirical block average process of the local observable
$\xi_i$.

For the scaled block average of the two conserved quantities we shall
also use the notation
\begin{eqnarray*}
\rn(t,x):=n^{2\b}\wih{\eta}(t,x),
\qquad
\un(t,x):=n^{\b} \wih{\zet}(t,x).
\end{eqnarray*}

Introducing   block averages the main terms on the right hand side of
(\ref{eq:elef}) and
(\ref{eq:peteef}) become:
\begin{eqnarray}
&&
\hskip-2cm
\label{eq:elefblocks}
\frac{1}{n}\sum_{j\in\Tn}
\px v(t,\frac{j}{n})
n^{3\b}\psi_{j}
+
\frac{1}{n}\sum_{j\in\Tn}
\px u(t,\frac{j}{n})
n^{2\b}\phi_{j}=
\\[5pt]
\notag
&&
\hskip2cm
\frac{1}{n}\sum_{j\in\Tn}
\px v(t,\frac{j}{n})
n^{3\b}\wih{\psi}(\frac{j}{n})
+
\frac{1}{n}\sum_{j\in\Tn}
\px u(t,\frac{j}{n})
n^{2\b}\wih{\phi}(\frac{j}{n})
\\[5pt]
\notag
&&
\hskip6cm
+
A^n_8(t,\uome)
+
A^n_9(t,\uome),
\end{eqnarray}
respectively
\begin{eqnarray}
\notag
&&
\hskip-3cm
\frac{1}{n}\sum_{j\in\Tn}
\pt v(t,\frac{j}{n})
\big(n^{2\b}\eta_j-\rho(t,\frac{j}{n})\big)
+
\frac{1}{n}\sum_{j\in\Tn}
\pt u(t,\frac{j}{n})
\big(n^{\b}\zeta_j-u(t,\frac{j}{n})\big)
=
\\[5pt]
\label{eq:peteefblocks}
&&
\hskip4cm
\phantom{+}
\frac{1}{n}\sum_{j\in\Tn}
\pt v(t,\frac{j}{n})
\big(n^{2\b}\wih{\eta}(\frac{j}{n})-\rho(t,\frac{j}{n})\big)
\\[5pt]
\notag
&&
\hskip4cm
+
\frac{1}{n}\sum_{j\in\Tn}
\pt u(t,\frac{j}{n})
\big(n^{\b}\wih{\zeta}(\frac{j}{n})-u(t,\frac{j}{n})\big)
\\[5pt]
\notag
&&
\hskip4cm
+
A^n_{10}(t,\uome)
+
A^n_{11}(t,\uome),
\end{eqnarray}
where the error terms are
\begin{eqnarray*}
&&
A^n_{8}(t,\uome)
:=
\frac{1}{n}
\sum_{j\in\Tn}
\Big(
\px v\big(t,\frac{j}{n}\big)-
\frac{1}{l}\sum_k a\big(\frac{j-k}{l}\big)
\px v\big(t,\frac{k}{n}\big)
\Big)
n^{3\b}\psi_j,
\\[5pt]
&&
A^n_{9}(t,\uome)
:=
\frac{1}{n}
\sum_{j\in\Tn}
\Big(
\px u\big(t,\frac{j}{n}\big)-
\frac{1}{l}\sum_k a\big(\frac{j-k}{l}\big)
\px u\big(t,\frac{k}{n}\big)
\Big)
n^{2\b}
\phi_j,
\\[5pt]
&&
A^n_{10}(t,\uome)
:=
\frac{1}{n}
\sum_{j\in\Tn}
\Big(
\pt v\big(t,\frac{j}{n}\big)-
\frac{1}{l}\sum_k a\big(\frac{j-k}{l}\big)
\pt v\big(t,\frac{k}{n}\big)
\Big)
n^{2\b}
\eta_j,
\\[5pt]
&&
A^n_{11}(t,\uome)
:=
\frac{1}{n}
\sum_{j\in\Tn}
\Big(
\pt u\big(t,\frac{j}{n}\big)-
\frac{1}{l}\sum_k a\big(\frac{j-k}{l}\big)
\pt u\big(t,\frac{k}{n}\big)
\Big)
n^{\b}
\zeta_j.
\end{eqnarray*}
These error terms will be estimated in subsection
\ref{subs:firsterrors}.

Since $[0,T]\times\T\ni(t,x)\mapsto (\rho(t,x),u(t,x))$, is
a \emph{smooth} solution of the pde (\ref{eq:pde}), we have
\begin{eqnarray*}
\pt v=-u\px v-\px u,
\qquad
\pt u=-\rho\px v -\gamma u \px u.
\end{eqnarray*}
Inserting these expressions into the main terms of
(\ref{eq:peteefblocks}) eventually we obtain for the integrand in
(\ref{eq:entprod})
\begin{eqnarray}
\label{eq:sumup0}
&&
\hskip-2cm
n^{3\b}
{\Ln\!\fnt}(\uome)
+
n^{3\b+\d}
{\Kn\!\fnt}(\uome)
+
n^{-1+2\b}
{\pt\fnt}(\uome)
=
\\[5pt]
\notag
&&
\phantom{M,}
\phantom{+}
\frac{1}{n}\sum_{j\in\Tn}
\px v(t,\frac{j}{n})
\Big\{
n^{3\b}\wih{\psi}(\frac{j}{n})
-
\rho(t,\frac{j}{n}) u(t,\frac{j}{n})
\\[5pt]
\notag
&&
\phantom{MMMMM}
-
u(t,\frac{j}{n})
\big(n^{2\b}\wih{\eta}(\frac{j}{n}) -\rho(t,\frac{j}{n}) \big)
-
\rho(t,\frac{j}{n})
\big(n^{\b}\wih{\zeta}(\frac{j}{n}) -u(t,\frac{j}{n}) \big)
\Big\}
\\[5pt]
\notag
&&
\phantom{M}
+
\frac{1}{n}\sum_{j\in\Tn}
\px u(t,\frac{j}{n})
\Big\{
n^{2\b}\wih{\phi}(\frac{j}{n})
-
\big(\rho(t,\frac{j}{n})+\gamma  u(t,\frac{j}{n})^2 \big)
\\[5pt]
\notag
&&
\phantom{MMMMM}
-
\big(n^{2\b}\wih{\eta}(\frac{j}{n}) -\rho(t,\frac{j}{n}) \big)
-
2\gamma u(t,\frac{j}{n})
\big(n^{\b}\wih{\zeta}(\frac{j}{n}) -u(t,\frac{j}{n}) \big)
\Big\}
\\[5pt]
\notag
&&
\phantom{MM}
+
\sum_{k=1}^{12} A^n_k(t,\uome),
\end{eqnarray}
where
\begin{eqnarray*}
A^n_{12}(t)
&:=&
\frac{1}{n}\sum_{j\in\Tn}
\big(
\big(\px v\big)
\rho u
+
\big(\px u\big)
\big(\rho+\gamma  u^2 \big)
\big)(t,\frac{j}{n})
\\[5pt]
&=&
\frac{1}{n}\sum_{j\in\Tn}
\px\big(\rho u + \frac{\gamma}{3} u^3\big)(t,\frac{j}{n})
\end{eqnarray*}
%

\subsection{The error terms $A^n_k$, $k=1,\dots,12$}
\label{subs:firsterrors}

\begin{lemma}
\label{lemma:firstentropybounds}
There exists a finite constant
$C$,  such that for any
$n\in\N$,
$t\in[0,T]$ and for any sequence of real numbers
$b_j$,
$j=1,\dots,n$ the following bounds hold:
\begin{eqnarray}
\label{eq:firstetab}
&&
\expect_{\mun_t}
\Big(
\frac{1}{n}\sum_{j\in\Tn}
b_j\eta_j
\Big)
\le
C
n^{-2\b}
\left(
\frac{1}{n}\sum_{j\in\Tn}
b_j
+
\sqrt{\frac{1}{n}\sum_{j\in\Tn}
b_j^2}
\,\right)
\\[5pt]
\label{eq:firstzetab}
&&
\expect_{\mun_t}
\Big(
\frac{1}{n}\sum_{j\in\Tn}
b_j\zeta_j
\Big)
\le
C
n^{-\b}
\sqrt{\frac{1}{n}\sum_{j\in\Tn}
b_j^2}
\\[5pt]
\label{eq:firstpsib}
&&
\expect_{\mun_t}
\Big(
\frac{1}{n}\sum_{j\in\Tn}
b_j\psi_j
\Big)
\le
C
n^{-2\b}
\left(
\frac{1}{n}\sum_{j\in\Tn}
b_j
+
\sqrt{\frac{1}{n}\sum_{j\in\Tn}
b_j^2}
\,\right)
\\[5pt]
\label{eq:firstphib}
&&
\expect_{\mun_t}
\Big(
\frac{1}{n}\sum_{j\in\Tn}
b_j\phi_j
\Big)
\le
C
\left(
\frac{1}{n}\sum_{j\in\Tn}
b_j
+
n^{-\b}
\sqrt{\frac{1}{n}\sum_{j\in\Tn}
b_j^2}
\,\right)
\end{eqnarray}
\end{lemma}

\begin{proof}
The proof relies on the entropy inequality
\begin{eqnarray}
\label{eq:firstentr}
&&
\hskip-8mm
\expect_{\mun_t}
\Big(
\frac{1}{n}\sum_{j\in\Tn}
b_j \big(\xi_j-\expect_{\pin}(\xi_j)\big)
\Big)
\le
\\[5pt]
\notag
&&
\hskip8mm
\frac{1}{\gamma n} H(\mun_t|\pin)
+%
\frac{1}{\gamma n}
\log \expect_{\pin}
\Big(
\exp
\Big\{
\gamma \sum_{j\in\Tn}
b_j\big(\xi_j-\expect_{\pin}(\xi_j)\big)
\Big\}
\Big),
\end{eqnarray}
where
$\xi_j$  stands for either of
$\eta_j$,
$\zeta_j$,
$\psi_j$ or
$\phi_j$. We note that all these variables are bounded and
\begin{eqnarray*}
\begin{array}{ll}
\big|\,\expect_{\pin}
\big(\eta_j\big)\,\big| \le C n^{-2\b},
\qquad
&
\var_{\pin}
\big(\eta_j\big) \le C n^{-2\b},
\\[10pt]
\big|\,\expect_{\pin}
\big(\zeta_j\big)\,\big| = 0,
\qquad
&
\var_{\pin}
\big(\zeta_j\big) \le C,
\\[10pt]
\big|\,\expect_{\pin}
\big(\psi_j\big)\,\big| = 0,
\qquad
&
\var_{\pin}
\big(\psi_j\big) \le C n^{-2\b},
\\[10pt]
\big|\,\expect_{\pin}
\big(\phi_j\big)\,\big| \le C,
\qquad
&
\var_{\pin}
\big(\phi_j\big) \le C.
\end{array}
\end{eqnarray*}
From these bounds and the entropy inequality
(\ref{eq:firstentr})  the statement of the lemma follows directly.
\end{proof}

Now we turn to the estimates on the error terms. We use the bounds
(\ref{eq:firstetab}),
(\ref{eq:firstzetab}),
(\ref{eq:firstpsib}) and
(\ref{eq:firstphib}) of Lemma
\ref{lemma:firstentropybounds}, the asymptotics
(\ref{eq:asym}) and uniform approximation of
$\px$  of smooth functions by their discrete derivative
$\grn$. Straightforward computations yield
\begin{eqnarray*}
\label{eq:a1bound}
&&
\expect_{\mun_t}
\big(
A^n_1(t)
\big)
\le C\big(n^{-1+\b}+ n^{-\b} \big)
=
o(1),
\\[5pt]
\label{eq:a2bound}
&&
\expect_{\mun_t}
\big(
A^n_2(t)
\big)
\le C \big(n^{-1+2\b} + n^{-\b} \big)
=
o(1),
\\[5pt]
\label{eq:a3bound}
&&
\expect_{\mun_t}
\big(
A^n_3(t)
\big)
\le Cn^{-1+\b}
=
o(1),
\\[5pt]
\label{eq:a4bound}
&&
\expect_{\mun_t}
\big(
A^n_4(t)
\big)
\le C\big(n^{-1+2\b} +n^{-1+\b} \big)
=
o(1),
\\[5pt]
\label{eq:a5bound}
&&
\expect_{\mun_t}
\big(
A^n_5(t)
\big)
\le C\big(n^{-1+\b+\d} +n^{-1+2\b+\d} \big)
=
o(1),
\\[5pt]
\label{eq:a6bound}
&&
\expect_{\mun_t}
\big(
A^n_6(t)
\big)
\le Cn^{-2\b}
=
o(1),
\\[5pt]
\label{eq:a7bound}
&&
\expect_{\mun_t}
\big(
A^n_7(t)
\big)
\le Cn^{-2\b}
=
o(1),
\\[5pt]
\label{eq:a8bound}
&&
\expect_{\mun_t}
\big(
A^n_8(t)
\big)
\le C\big(n^{-1+\b}+n^{\b}l^{-1}+ n^{-1+\b}l\big)
=
o(1),
\\[5pt]
\label{eq:a9bound}
&&
\expect_{\mun_t}
\big(
A^n_9(t)
\big)
\le C\big(n^{-1+2\b}+n^{\b}l^{-1}+ n^{-1+\b}l \big)
=
o(1),
\\[5pt]
\label{eq:a10bound}
&&
\expect_{\mun_t}
\big(
A^n_{10}(t)
\big)
\le C\big(n^{-1}+l^{-1}+  n^{-1}l\big)
=
o(1),
\\[5pt]
\label{eq:a11bound}
&&
\expect_{\mun_t}
\big(
A^n_{11}(t)
\big)
\le C\big(l^{-1}+  n^{-1}l\big)
=
o(1).
\end{eqnarray*}
Finally,
$A^n_{12}(t)$  is a simple numerical error term (no probability
involved):
\begin{eqnarray*}
A^n_{12}(t)
\le C n^{-1}
=
o(1).
\end{eqnarray*}
%


\subsection{Sumup}
\label{subs:sumup1}

Thus, integrating
(\ref{eq:entprod}), using
(\ref{eq:sumup0}) and the bounds of subsection
\ref{subs:firsterrors} we obtain
\begin{eqnarray}
\label{eq:sumup1}
\hn(t)
=
\int_0^t{\cal A}^n(s)\,ds
+
\int_0^t{\cal B}^n(s)\,ds
-
\sn(t)
+
o(1),
\end{eqnarray}
where
\begin{eqnarray}
\label{eq:curlya}
{\cal A}^n(t)
:=
\expect\Big(
\frac{1}{n}\sum_{j\in\Tn}
\Big\{
\big(\px v\big)
\big\{
n^{3\b}\wih{\psi}
-
\rho u
-
u
(\rn -\rho )
-
\rho
(\un -u )
\big\}
\Big\}
(t,\frac{j}{n})
\Big)
\end{eqnarray}
and
\begin{eqnarray}
\label{eq:curlyb}
{\cal B}^n(t)
:=
\expect\Big(
\frac{1}{n}\sum_{j\in\Tn}
\Big\{
\big(\px u\big)
\big\{
n^{2\b}\wih{\phi}
-
(\rho+\gamma  u^2 )
-
(\rn -\rho)
-
2\gamma u
(\un -u)
\big\}
\Big\}
(t,\frac{j}{n})
\Big)
\end{eqnarray}
The main difficulty is caused by
${\cal A}^{n}(t)$. The term
${\cal B}^{n}(t)$ is estimated exactly as it is done in
\cite{tothvalko1} for the one-component systems: since
$\Phi(\rho,u)=\rho+\gamma u^2$ is linear in
$\rho$ and quadratic in
$u$ no problem is caused by the low particle density. By repeating the
arguments of
\cite{tothvalko1} we obtain
\begin{eqnarray}
\label{eq:regi}
\int_0^t{\cal B}^{n}(s)ds
\le
C \int_0^t\hn(s)ds + o(1).
\end{eqnarray}

In the rest of the proof we concentrate on the essentially difficult
term
${\cal A}^{n}(t)$.

\section{Cutoff}
\label{section:cutoff}

We define the
\emph{rescaled macroscopic fluxes}
\begin{eqnarray}
\label{eq:scaledfluxes}
\Psin(\rho,u):=n^{3\beta}\Psi(n^{-2\beta}\rho,n^{-\beta}u),
\qquad
\Phin(\rho,u):=n^{2\beta}\Phi(n^{-2\beta}\rho,n^{-\beta}u).
\end{eqnarray}
defined on the scaled domain

\beq
\label{eq:scaleddomain}
{\cal D}^n
:=
\{(\rho,u): (n^{-2\b}\rho,n^{-\b}u)\in{\cal D}\}.
\eeq
The  first partial derivatives of the scaled fluxes are
\begin{eqnarray}
\label{eq:scaledfluxderiv}
\begin{array}{ll}
\Psin_\rho(\rho,u)=n^{\beta}\Psi_\rho(n^{-2\beta}\rho,n^{-\beta}u),
&
\Phin_\rho(\rho,u)=\Phi_\rho(n^{-2\beta}\rho,n^{-\beta}u),
\\[8pt]
\Psin_u(\rho,u)=n^{2\beta}\Psi_u(n^{-2\beta}\rho,n^{-\beta}u),
&
\Phin_u(\rho,u)=n^{\beta}\Phi_u(n^{-2\beta}\rho,n^{-\beta}u).
\end{array}
\end{eqnarray}
For any
$(\rho,u)\in\R_+\times\R$
\begin{eqnarray}
\label{eq:fluxlimits}
\begin{array}{ll}
\displaystyle
\lim_{n\to\infty}
\Psin(\rho,u)
=
\rho u,
\qquad
&
\displaystyle
\lim_{n\to\infty}
\Phin(\rho,u)
=
\rho+\gamma u^2,
\\[8pt]
\displaystyle
\lim_{n\to\infty}
\Psin_\rho(\rho,u)
=
u,
&
\displaystyle
\lim_{n\to\infty}
\Phin_\rho(\rho,u)
=
1,
\\[8pt]
\displaystyle
\lim_{n\to\infty}
\Psin_u(\rho,u)
=
\rho,
&
\displaystyle
\lim_{n\to\infty}
\Phin_u(\rho,u)
=2\gamma u.
\end{array}
\end{eqnarray}
The convergence is uniform in compact subsets of
$\R_+\times\R$

Note  that
\begin{eqnarray*}
&&
\Psin(\rn(t,x),\un(t,x))
=
n^{3\beta}\Psi(\wih{\eta}(t,x),\wih{\zet}(t,x)),
\\[5pt]
&&
\Phin(\rn(t,x),\un(t,x))
=
n^{2\beta}\Phi(\wih{\eta}(t,x),\wih{\zet}(t,x)).
\end{eqnarray*}
%

\subsection{The direct approach --- why it fails?}
\label{subs:naive}

The most natural thing is to write the summand in
${\cal A}^{n}(t)$ as
\begin{eqnarray}
\label{eq:directapproach}
&&
\hskip-25mm
n^{3\b}\wih{\psi}
-
\rho u
-
u
( \rn - \rho )
-
\rho
(\un - u )
=
\\ [5pt]
\notag
&&
n^{3\b}(\wih{\psi}-\Psi(\wih{\eta},\wih{\zet}) )
+
(\Psi^{n}(\rn,\un) - \rn\un )
+
(\rn-\rho )(\un -u )
\end{eqnarray}
By applying Varadhan's ``one block estimate'' and controlling the
error terms in the Taylor expansion of
$\Psi$, the first two terms on the right hand side can be dealt with.
However, the last term causes serious problems:
with proper normalization, it is distributed with respect to
the local equilibrium measure
$\nunt$, like a product of independent Poisson and  Gaussian random
variables, and thus it  does
\emph{not}  have a finite exponential  moment. Since the robust
estimates heavily rely on the entropy inequality where the finite
exponential moment is needed, we  have to choose another approach for
estimating
${\cal A}^n(t)$.

Instead of  writing plainly
({\ref{eq:directapproach}}), we introduce a cutoff.
We let
\[
M>\sup\{\rho(t,x)\lor|u(t,x)|:(t,x)\in[0,T]\times\T\}.
\]
The value of $M$ will be specified by the large deviation bounds given
in Proposition \ref{propo:ldbounds} (via Lemma
\ref{lemma:largedevi}).

Let
$I^n,J^n:\R_+\times\R\to\R$ be bounded functions so that
\beqs
&&
I^n+J^n=1,
\\[5pt]
&&
I^n(\rho,u)=1
\quad
\text{ for }
\quad
\rho\lor|u|\le M,
\\[5pt]
&&
I^n(\rho,u)=0
\quad
\text{ for `large' }
(\rho,u).
\eeqs
The last property will be specified later.

We  split the right hand side of
(\ref{eq:directapproach}) in a most natural way, according to this
cutoff:
\begin{eqnarray}
\label{eq:thecutoff}
&&
\hskip-25mm
n^{3\b}\wih{\psi}
-
\rho u
-
u
(\rn - \rho )
-
\rho
(\un - u )
=
\\ [5pt]
\notag
&&
\phantom{+}
n^{3\b}\wih{\psi} \J(\rn,\un)
-
(\rn u + \rho \un -\rho u) \J(\rn,\un)
\\ [5pt]
\notag
&&
+
n^{3\b}(\wih{\psi}-\Psi(\wih{\eta},\wih{\zet}) ) \I(\rn,\un)
\\ [5pt]
\notag
&&
+
(\Psi^{n}(\rn,\un) - \rn \un ) \I(\rn,\un)
+
(\rn-\rho )(\un -u ) \I(\rn,\un)
\end{eqnarray}
The second term on the right hand side is linear in the block
averages, so it does not cause any problem. The third term is
estimated by use of Varadhan's one block estimate. The fourth term
is Taylor approximation. Finally, the last term can be handled
with the entropy inequality  \emph{if the cutoff
$I^n(\rho,u)$ is strong enough} to tame the tail of the
Gaussian$\times$Poisson random variable.

The main difficulty is caused by the first term on the right hand
side. This term certainly can not be estimated with the robust method,
i.e. with entropy inequality: we would run into the same problem we
wanted to overcome by introducing the cutoff. The only way
this term may be  small is by some cancellation. It turns out that the
desired cancellations indeed occur (in form of a martingale appearing
in the space-time average) if and  only if
\beq
\label{eq:jn}
J^n(\rho,u) = S^n_\rho(\rho,u),
\eeq
where
$S^n(\rho,u)$ is a particular
\emph{Lax entropy of the scaled Euler equation}
\begin{eqnarray}
\label{eq:scaledeuler}
\left\{
\begin{array}{l}
\pt\rho+\px\Psin(\rho,u)=0
\\[7pt]
\pt u  +\px\Phin(\rho,u)=0,
\end{array}
\right.
\end{eqnarray}
with $\Psin(\rho,u)$ and $\Phin(\rho,u)$ defined in
(\ref{eq:scaledfluxes}). That is $\Sn$ is solution of the pde \beq
\label{eq:scaledlentpde} \Psin_{u}\Sn_{\rho\rho}
+\big(\Phin_u-\Psin_\rho\big)\Sn_{\rho u} -\Phin_\rho \Sn_{uu} =0.
\eeq

\subsection{The cutoff function}
\label{subs:cutoff}

In the present subsection we construct the cutoff function
(\ref{eq:jn}) and we state some estimates related to it. These
bounds will be of paramount importance in our further proof. They
are stated in the technical Lemmas \ref{lemma:charcurves},
\ref{lemma:lentbounds} and \ref{lemma:lentfluxbound}. The proof of
these lemmas is pure classical pde theory
and it is postponed to section
\ref{section:cutoffproofs}.

First, in subsubsection
\ref{subsubs:unscaled},
we formulate our construction and estimates in terms of Lax entropies
of the \emph{unscaled} Euler equation
(\ref{eq:euler}). Then in subsubsection
\ref{subsubs:scaled} we rescale these estimates in order to get the
necessary bounds on $\Sn$ and its derivatives.

\subsubsection{}
\label{subsubs:unscaled}

A Lax entropy/flux pair $S(\rho,u)$, $F(\rho,u)$ of the system
(\ref{eq:euler}) is solution of the system of pdes \beq
\label{eq:lentfluxpde} F_\rho=\Psi_\rho S_\rho+\Phi_\rho S_u,
\qquad F_u=\Psi_u S_\rho+\Phi_u S_u, \eeq defined on ${\cal D}$.
In particular the Lax entropy $S(\rho,u)$ solves the pde: \beq
\label{eq:lentpde} \Psi_{u} S_{\rho\rho}
+\big(\Phi_u-\Psi_\rho\big)S_{\rho u} -\Phi_\rho S_{uu} =0, \eeq

The linear pde (\ref{eq:lentpde}) is hyperbolic in ${\cal D}$.
One family of its characteristic  curves are solutions of the
following ODE, meant in the domain ${\cal D}$:
\begin{eqnarray}
\label{eq:charcurvesode}
\frac{d\rho}{du}=
\frac{
\sqrt{\big(\Phi_u-\Psi_\rho\big)^2 + 4\Phi_\rho\Psi_u}
-
\big(\Phi_u-\Psi_\rho\big)
}
{2\Phi_\rho},
\end{eqnarray}
The other family is obtained by reflecting $u$ to $-u$.

First we conclude that the line segment ${\cal D}\cap\{u=0\}$ is
\emph{not} characteristic for the hyperbolic pde
(\ref{eq:lentpde}).
That is: is intersects transversally the characteristic lines defined
by the differential equation (\ref{eq:charcurvesode}).
Indeed, from the Onsager relation
(\ref{eq:onsager}) and obvious parity considerations it follows,
that the right hand  side of (\ref{eq:charcurvesode}) restricted
to $\{u=0\}$ becomes
$\big(\var_{r,0}(\eta)/\var_{r,0}(\zeta)\big)^{1/2}$ and this
expression  is obviously finite for $r\in[0,\rho^*)$. It follows
that the Cauchy problem (\ref{eq:lentpde}), with the following
initial condition: \beq \label{eq:inicond} S(r,0) = s(r),\,\,
\quad S_u(r,0)=0,\,\, \qquad r\in[0,\rho^*) \eeq is \emph{well
posed}.

In our concrete problem the function
$s(r)$ will be chosen as follows:  we fix
$0<\ulr<\olr$, and define
\beq
\label{eq:smalls}
s(r)
=
\left\{
\begin{array}{ll}
0
&
\text{ if }  r\in[0, \ulr),
\\[8pt]
\displaystyle
\frac{r\log(r/\ulr)-(r-\ulr)}{\log(\olr/\ulr)}
&
\text{ if } r\in[\ulr,\olr),
\\[8pt]
\displaystyle
r-\frac{\olr-\ulr}{\log(\olr/\ulr)}
&
\text{ if } r\in[\olr,\infty).
\end{array}
\right.
\eeq
Note that
$s(r)$ and
$s'(r)$ are continuous.

We first analyze the global structure of the characteristic curves.
Due to the assumption
(\ref{cond:phi0'-psi1}) imposed on-, and regularity of the flux
functions
$\Phi$ and
$\Psi$, there exists some
$\rho_0>0$ such that the ODE
(\ref{eq:charcurvesode}) is regular in
$\{(\rho,u)\in {\cal  D}:\rho<\rho_0\text{ and }
(\rho,u)\not=(0,0)\}$.
We shall not be concerned about what happens outside this strip.
Denote by
$\sigma(u;r)$ the solution of the ODE
(\ref{eq:charcurvesode})
with initial condition
$\sigma(0;r)=r$.

\begin{lemma}
\label{lemma:charcurves} There exist constants $0<C_2<C_1<\infty$
and $r_0>0$ such that for any $r\in[0,r_0]$ \beq
\begin{array}{rl}
r+C_1\sqrt{r}u \le \sigma(u;r) \le r+C_2\sqrt{r}u & \quad \text{
if } u\le0,
\\[8pt]
r \le \sigma(u;r) \le r+C_1 \left(\sqrt{r} u\wedge
r^{\frac{4\g-3}{4\g-2}} u^{\frac1{2\g-1}} \right)& \quad \text{ if
} u\ge0.
\end{array}
\eeq The inequalities are valid as long as
$(\sigma(u;r),u)\in{\cal D}$.  The map $u\mapsto\sigma(u;r)$ is
regular and monotone increasing.
\end{lemma}

\noindent
See subsection
\ref{subs:charcurvesproof} for the proof of this lemma.

For
$r<r_0$ we partition the domain
${\cal D}$  in three parts as follows
\beqs
{\cal D}_1(r)
&:=&
\{(\rho,u)\in{\cal D}:\rho<\sigma(-|u|;r)\}
\\[5pt]
{\cal D}_2(r)
&:=&
\{(\rho,u)\in{\cal D}:\rho>\sigma(|u|;r)\}
\\[5pt]
{\cal D}_3(r)
&:=&
\{(\rho,u)\in{\cal D}:\sigma(-|u|;r)\le\rho\le\sigma(|u|;r)\}
\\[5pt]
&\phantom{:}=&
{\cal D}\setminus\big({\cal D}_1(r)\cup{\cal D}_2(r)\big).
\eeqs
See Figure 1 for a sketch of the domains ${\cal D}_1(r),{\cal
  D}_2(r),{\cal D}_3(r)$.

From Lemma
\ref{lemma:charcurves} it follows that
\beq
\notag
&&
\hskip-8mm
\{(\rho,u):0\le\rho< r-C_1\sqrt{r}|u|\}
\,
\subset
\,
{\cal D}_1(r)
\,
\subset
\,
\{(\rho,u):0\le\rho< r-C_2\sqrt{r}|u|\},
\\[-5pt]
\label{eq:D1D2bounds}
\\[-12pt]
\notag
&&
\hskip-8mm
\{(\rho,u):r+C_1\sqrt{r}|u|<\rho \}
\,
\subset
\,
{\cal D}_2(r)
\,
\subset
\,
\{(\rho,u):r<\rho\},
\eeq
and for
$0\le r\le r'\le r_0$
\beq
\label{eq:D12ordered}
{\cal D}_1(r)\,\subset\,{\cal D}_1(r'),
\qquad
{\cal D}_2(r')\,\subset\,{\cal D}_2(r).
\eeq

\begin{figure}[htp]
\centering
\includegraphics[width=(\textwidth-30mm)]{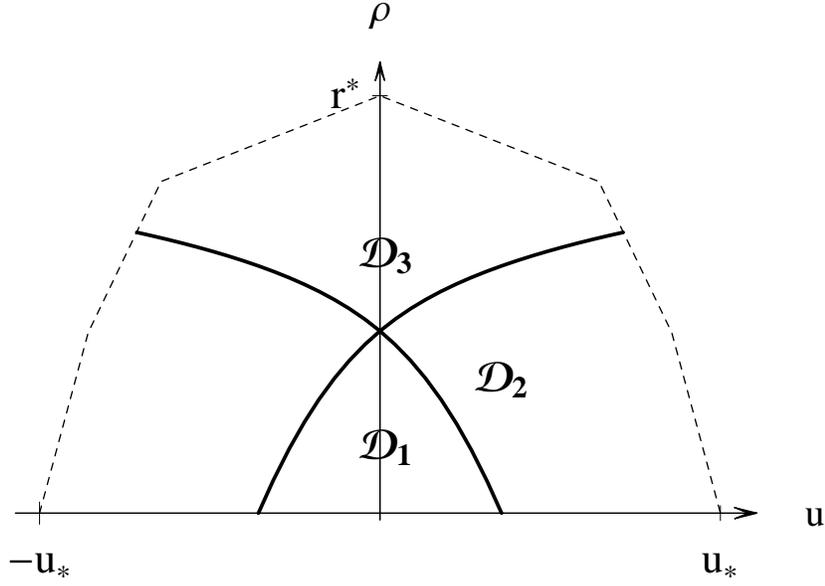}
\caption{$\mathcal{D}_1,\mathcal{D}_2,\mathcal{D}_3 $}
\label{fig:D123}
\end{figure}

From now $r_0$ is \emph{fixed for ever} and we denote
\[
\wt{\cal D}:={\cal D}_1(r_0).
\]
This domain is a \emph{rectangle} in characteristic coordinates with
diagonal
$\wt{\cal D}\cap\{u=0\}$, as opposed to
${\cal D}$ which may not be a full characteristic rectangle.
(Actually, choosing the characteristic coordinates in a natural
symmetric way, $z(\rho,u)=w(\rho,-u)$, the domain $\wt{\cal D}$ is a
square in characteristic coordinates.)
Note that
${\cal D}_3(r_0)\cap\{u\ge0\}$ and
${\cal D}_3(r_0)\cap\{u\le0\}$ are also characteristic rectangles.

Next we turn to the construction of a particular family of Lax
entropies
which will serve for obtaining the cutoff functions needed.
We fix $0<\ulr<\olr<r_0$.  and define $S:{\cal D}\to \R$ as follows:

\begin{enumerate}[(i)]

\item
\label{i}
In
$\wt{\cal D}$: $S(\rho,u)$ is solution of the Cauchy problem
(\ref{eq:lentpde})+(\ref{eq:inicond}) with
$s(r)$ given in
(\ref{eq:smalls}). Note that
\beq
\label{eq:lentinD12}
S(\rho,u)=
\left\{
\begin{array}{ll}
0
&
\text{ if } (\rho,u)\in {\cal D}_1(\ulr)\subset\wt{\cal D},
\\[8pt]
\displaystyle
\rho-\frac{\olr-\ulr}{\log(\olr/\ulr)}
&
\text{ if } (\rho,u)\in {\cal D}_2(\olr)\cap\wt{\cal D}.
\end{array}
\right.
\eeq

\item
\label{ii}
In
${\cal D}_2(\olr)$:
\beq
\label{eq:lentinD2}
S(\rho,u):=
\rho -\frac{\olr-\ulr}{\log(\olr/\ulr)},
\quad
\text{ if } (\rho,u)\in {\cal D}_2(\olr)
\eeq
Note that there is no contradiction: in
$\wt{\cal D}\cap{\cal D}_2(\ulr)$,
(\ref{i}) yields the same expression.

\item
In
${\cal D}_3(\olr)\setminus\wt{\cal D}$:
$S(\rho,u)$ is defined as solution of the Goursat problem
(\ref{eq:lentpde}) with boundary conditions on the characteristic
lines
$\partial\wt{\cal D}\cap{\cal D}_3(\olr)$, respectively,
$\partial{\cal D}_2(\olr)\setminus\wt{\cal D}$ provided by
(\ref{i}), respectively, (\ref{ii}).

\end{enumerate}

\noindent Note that $S(\rho,u)$ is solution of the pde
(\ref{eq:lentpde}), \emph{globally} in $\wt{\cal D}$.

We denote
\beqs
{\cal D}_3(\ulr,\olr)
&:=&
\phantom{\cup}
{\cal D}\setminus\big({\cal D}_1(\ulr)\cup{\cal D}_2(\olr)\big)
\\
&\phantom{:}=&
\phantom{\cup}
\big({\cal D}_1(\olr)\cap{\cal D}_2(\ulr)\big)
\cup
\big({\cal D}_1(\olr)\cap{\cal D}_3(\ulr)\big)
\\[5pt]
&&
\cup
\big({\cal D}_3(\olr)\cap{\cal D}_2(\ulr)\big)
\cup
\big({\cal D}_3(\olr)\cap{\cal D}_3(\ulr)\big).
\eeqs

\begin{figure}[htp]
\centering
\includegraphics[width=(\textwidth-30mm)]{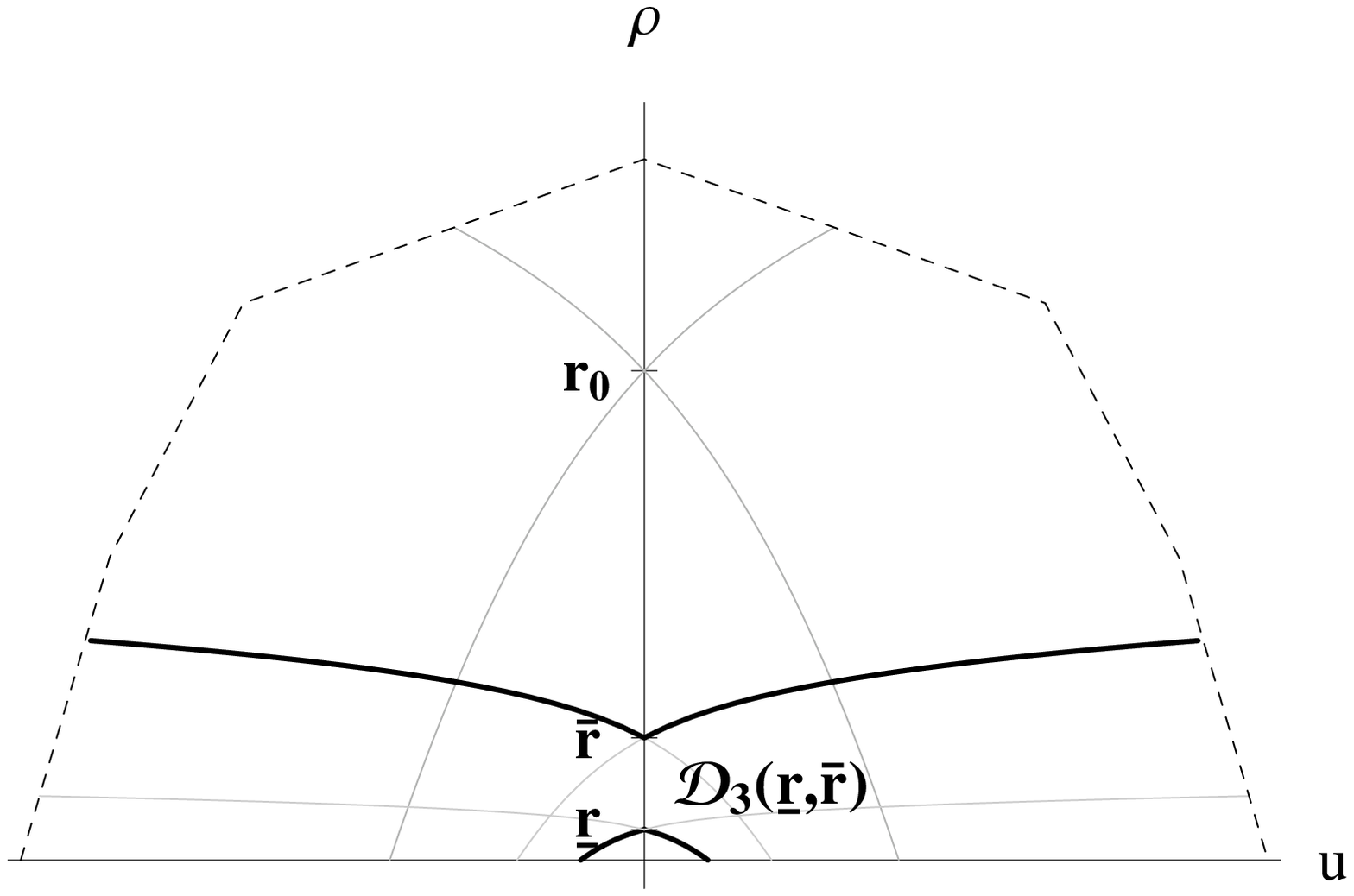}
\caption{$\mathcal{D}_3(\ulr,\olr)$} \label{fig:D3r}
\end{figure}

The following lemma provides the necessary bounds on the partial
derivatives of $S(\rho,u)$ (up to second order) in the domain
${\cal D}_3(\ulr,\olr)$.

\begin{lemma}
\label{lemma:lentbounds}
There exists a constant
$C<\infty$ such that for any
$0<\ulr<\olr<r_0$ and
$S(\rho,u)$ defined as above
the following global bounds hold:
\beq
\label{eq:s'rhobound}
\hskip-10mm
|S_\rho(\rho,u)-\ind_{{\cal D}_2(\olr)}(\rho,u)|
&\le&
C
\,
\ind_{{\cal D}_3(\ulr,\olr)}(\rho,u),
\\[5pt]
\label{eq:s'ubound}
|S_u(\rho,u)|
&\le&
\frac{C(\sqrt{\olr}-\sqrt{\ulr})}{\log(\olr/\ulr)}
\,
\ind_{{\cal D}_3(\ulr,\olr)}(\rho,u),
\\[5pt]
\label{eq:s''rhorhobound}
|S_{\rho\rho}(\rho,u)|
&\le&
\frac{C}{\log(\olr/\ulr)}
\,
\frac{1}{\ulr+\rho}
\,
\ind_{{\cal D}_3(\ulr,\olr)}(\rho,u),
\\[5pt]
\label{eq:s''rhoubound}
|S_{\rho u}(\rho,u)|
&\le&
\frac{C}{\log(\olr/\ulr)}
\,
\frac{1}{\sqrt{\ulr}+\sqrt{\rho}+\abs{u}}
\,
\ind_{{\cal D}_3(\ulr,\olr)}(\rho,u),
\\[5pt]
\label{eq:s''uubound}
|S_{uu}(\rho,u)|
&\le&
\frac{C}{\log(\olr/\ulr)}
\,
\ind_{{\cal D}_3(\ulr,\olr)}(\rho,u),
\eeq
\end{lemma}

\noindent
This lemma is proved  in subsection
\ref{subs:lentboundsproof}.

Beside the bounds on the partial derivatives of
$S(\rho,u)$ we shall also need a bound  on the function
$F(\rho,u)-\Psi(\rho,u) S_\rho(\rho,u)$. From
(\ref{eq:lentfluxpde}) and
(\ref{eq:lentinD12}) it follows that
\beqs
\label{eq:lentfluxinD12}
F(\rho,u)=
\left\{
\begin{array}{cl}
0
&
\text{ if } (\rho,u)\in {\cal D}_1(\ulr),
\\[8pt]
\displaystyle
\Psi(\rho,u)
&
\text{ if } (\rho,u)\in {\cal D}_2(\olr).
\end{array}
\right.
\eeqs

\begin{lemma}
\label{lemma:lentfluxbound}
With the assumptions and notations of Lemma
\ref{lemma:lentbounds}
\beq
\label{eq:fbound}
|F(\rho,u)-\Psi(\rho,u)S_\rho(\rho,u)|
\le
\frac{C \sqrt{\olr}}{\log(\olr/\ulr)}(\olr+u^2)
\ind_{{\cal D}_3(\ulr,\olr)}(\rho,u).
\eeq
\end{lemma}

\noindent
See subsection
\ref{subs:lentfluxboundproof}
for the proof of this lemma.

\medskip

\subsubsection{}
\label{subsubs:scaled}

The \emph{scaled} functions $\Sn(\rho,u)$, $\Fn(\rho,u)$ are
defined on the scaled domain ${\cal D}^n$ given in
(\ref{eq:scaleddomain}), as follows: fix $0<\ulr<\olr<\infty$ and
define the \emph{unscaled} Lax entropy/flux pair as in the
previous section but with \emph{downscaled} initial conditions
\beq \label{eq:downscaledinicond} S(r,0) = n^{-2\b}s(n^{2\b}r)
,\,\, \quad S_u(r,0)=0.\,\, \qquad r\in[0,\rho^*), \eeq with the
function $r\mapsto s(r)$ given in (\ref{eq:smalls}). Now,   define
the pair of scaled functions $\Sn,\Fn:{\cal D}^n\to\R$ as \beq
\Sn(\rho,u):=n^{2\b}S(n^{-2\b}\rho,n^{-\b}u), \qquad
\Fn(\rho,u):=n^{3\b}F(n^{-2\b}\rho,n^{-\b}u). \eeq It is
straightforward to check that $\Sn$, $\Fn$ form a Lax entropy/flux
pair of the pde (\ref{eq:scaledeuler}): %
\beq%
\label{eq:lentfluxpde_n}
\Fn_\rho=
\Psin_\rho \Sn_\rho+\Phin_\rho\Sn_u,
\qquad
\Fn_u=
\Psin_u \Sn_\rho+\Phin_u \Sn_u,
\eeq%
in particular $\Sn$ solves the pde (\ref{eq:scaledlentpde}).

We partition the scaled domain
\[
{\cal D}^n
=
{\cal D}^n_1(\ulr)
\cup
{\cal D}^n_2(\olr)
\cup
{\cal D}^n_3(\ulr,\olr)
\]
with the partition elements
\beqs
&&
{\cal D}^n_1(\ulr)
:=
\{(\rho,u)\in{\cal D}_n: (n^{-2\b}\rho,n^{-\b}u)\in{\cal
  D}_1(n^{-2\b}\ulr)\}
\\[5pt]
&&
{\cal D}^n_2(\olr)
:=
\{(\rho,u)\in{\cal D}_n: (n^{-2\b}\rho,n^{-\b}u)\in{\cal
  D}_2(n^{-2\b}\olr)\}
\\[5pt]
&&
{\cal D}^n_3(\ulr,\olr)
:=
\{(\rho,u)\in{\cal D}_n: (n^{-2\b}\rho,n^{-\b}u)\in{\cal
  D}_3(n^{-2\b}\ulr,n^{-2\b}\olr)\}.
\eeqs

In the following Proposition we summarize our main estimates
formulated in terms of the scaled objects. The  statement is a mere
corollary of the previous lemmas. It follows by simple scaling from
(\ref{eq:D1D2bounds}) and
(\ref{eq:s'rhobound})-(\ref{eq:fbound}).

\begin{proposition}
\label{prop:s's''bounds}
There exists a constant $C<\infty$, such that given any
$0<\ulr<\olr<\infty$ and the scaled Lax entropy/flux pair defined as
prescribed above, the following bounds hold uniformly in $n$:
\beq
\label{eq:scaledD1bound}
{\cal D}_1^n(\ulr)
&\supset&
\{(\rho,u):0\le\rho\le\ulr-C^{-1}\sqrt{\ulr}|u|\}
\\[5pt]
\label{eq:scaledD3bound}
{\cal D}_3^n(\ulr,\olr)
&\subset&
\{(\rho,u):\rho\le\olr+C\sqrt{\olr}|u|\}
\eeq
\beq
\label{eq:sn'rhobound}
\hskip-10mm
|\Sn_\rho(\rho,u)-\ind_{{\cal D}^n_2(\olr)}(\rho,u)|
&\le&
C
\,
\ind_{{\cal D}^n_3(\ulr,\olr)}(\rho,u),
\\[5pt]
\label{eq:sn'ubound}
|\Sn_u(\rho,u)|
&\le&
\frac{C(\sqrt{\olr}-\sqrt{\ulr})}{\log(\olr/\ulr)}
\,
\ind_{{\cal D}^n_3(\ulr,\olr)}(\rho,u),
\\[5pt]
\label{eq:sn''rhorhobound}
|\Sn_{\rho\rho}(\rho,u)|
&\le&
\frac{C}{\log(\olr/\ulr)}
\,
\frac{1}{\ulr+\rho}
\,
\ind_{{\cal D}^n_3(\ulr,\olr)}(\rho,u),
\\[5pt]
\label{eq:sn''rhoubound}
|\Sn_{\rho u}(\rho,u)|
&\le&
\frac{C}{\log(\olr/\ulr)}
\,
\frac{1}{\sqrt{\ulr}+\sqrt{\rho}+\abs{u}}
\,
\ind_{{\cal D}^n_3(\ulr,\olr)}(\rho,u),%
\\[5pt]
\label{eq:sn''uubound}
|\Sn_{uu}(\rho,u)|
&\le&
\frac{C}{\log(\olr/\ulr)}
\,
\ind_{{\cal D}^n_3(\ulr,\olr)}(\rho,u),
\eeq
\beq
\label{eq:fnbound}
\hskip-5mm
|\Fn(\rho,u)-\Psin(\rho,u)\Sn_\rho(\rho,u)|
&\le&
\frac{C \sqrt{\olr}}{\log(\olr/\ulr)} (\olr+u^2)
\ind_{{\cal D}^n_3(\ulr,\olr)}(\rho,u).
\eeq
\end{proposition}

Due to
(\ref{eq:scaledD1bound}) we can choose
$\ulr$ large enough to ensure that for all
$n$
\beq
\label{eq:ulrchoice}
\{(\rho,u):\rho\lor|u|\le M\}
\subset
{\cal D}^n_1(\ulr),
\eeq
where $M$ is specified by the large deviation bounds given in
Proposition \ref{propo:ldbounds} (via Lemma
\ref{lemma:largedevi}).

Further on, we choose $\olr$ so large that
\beq
\label{eq:olrchoice}
\frac{C}{\log(\olr/\ulr)}
<
\left(100 \sup_{(t,x)\in[0,T]\times\T}\abs{\log
  \rho(t,x)}\right)^{-1},
\eeq
and thus the bounds
(\ref{eq:sn''rhorhobound})-(\ref{eq:sn''uubound})
of Proposition \ref{prop:s's''bounds} are sufficient for our further
purposes.

\subsection{Outline of the further steps of proof}
\label{subs:outline2}

In section
\ref{section:tools}
we present the main probabilistic technical ingredients of the
forthcoming  proof. These are variants of  entropy inequalities and of
the celebrated one and two block estimates.

\noindent
In Section \ref{section:firstterm} we give an estimate
for the terms with 'large' values of $\rruu$, we prove that
\begin{eqnarray}
\label{eq:firstterm}
&&
\hskip-5mm
\abs{
\expect
\left(
\int_0^t\frac1n\sum_{j\in\Tn}
\Big\{
\big(\px v\big)
\big(n^{3\beta}\wih{\psi}\big)
\J(\rn,\un)
\Big\}
(s,\frac{j}{n})\,ds
\right)
}
\\[5pt]
\notag
&&
\hskip40mm
\le
\frac12\,\hn(t)
+
\frac12\,\sn(t)
+
C\, \int_0^t\hn(s)\,ds
+
o(1).
\end{eqnarray}

\noindent
In Section \ref{section:otherterms} we estimate the terms with
'small' values of $\rruu$, the section is divided into four
subsections.

\noindent
In subsection \ref{subs:secondterm} we prove
\begin{eqnarray}
\label{eq:secondterm}
&&
\hskip-8.9mm
\abs{
\expect
\left(
\frac1n\sum_{j\in\Tn}
\Big\{
\big(\px v\big)\,
\big(
\rn u+ \rho \un - \rho u
\big)\,
\J(\rn,\un)
\Big\}
(s,\frac{j}{n})
\right)
}
\\[5pt]
\notag
&&
\hskip75mm
\le
C\, \hn(s)
+o(1).
\end{eqnarray}

\noindent
In subsection \ref{subs:thirdterm} we prove the one block estimate
\begin{eqnarray}
\label{eq:thirdterm}
&&
\hskip-8mm
\abs{
\expect
\left(
\int_0^t\frac1n\sum_{j\in\Tn}
\Big\{
\big(\px v\big)\,
n^{3\beta}
\big(
\wih{\psi} - \Psi(\wih{\eta},\wih{\zeta})
\big)\,
\I(\rn,\un)
\Big\}
(s,\frac{j}{n})\,ds
\right)
}
\\
\notag
&&
\hskip100mm
=
o(1).
\end{eqnarray}

\noindent
In subsection \ref{subs:fourthterm} we control the Taylor
approximation
\begin{eqnarray}
\label{eq:fourthterm}
&&
\hskip-7.9mm
\abs{
\expect
\left(
\frac1n\sum_{j\in\Tn}
\Big\{
\big(\px v\big)\,
\big(\Psi^n(\rn,\un) - \rn\un\big)\,
\I(\rn,\un)
\Big\}
(s,\frac{j}{n})
\right)
}
\\[5pt]
\notag
&&
\hskip75mm
\le
C\, \hn(s)
+o(1).
\end{eqnarray}

\noindent
Finally, in subsection \ref{subs:fifthterm} we control the
fluctuations
\begin{eqnarray}
\label{eq:fifthterm}
&&
\hskip-6.5mm
\abs{
\expect
\left(
\frac1n\sum_{j\in\Tn}
\Big\{
\big(\px v\big)\,
\big(\rn-\rho\big)
\big(\un-u\big)\,
\I(\rn,\un)
\Big\}
(s,\frac{j}{n})
\right)
}
\\[5pt]
\notag
&&
\hskip80mm
\le
C\,
\hn(s)
+o(1).
\end{eqnarray}
%
Having all these done, from
(\ref{eq:curlya}),
(\ref{eq:thecutoff}) and the bounds
(\ref{eq:firstterm}),
(\ref{eq:secondterm}),
(\ref{eq:thirdterm}),
(\ref{eq:fourthterm}),
(\ref{eq:fifthterm}) it follows that
\begin{eqnarray}
\label{eq:uj}
\int_0^t{\cal A}^{n}(s)ds
\le
\frac12 \,\hn(t)
+
\frac12 \,\sn(t)
+
C \,\int_0^t\hn(s)ds
+
o(1).
\end{eqnarray}
Finally, from (\ref{eq:sumup1}), (\ref{eq:regi}), (\ref{eq:uj})
and noting that $\sn(t)\geq 0$ we get the desired Gr\"onwall
inequality (\ref{eq:gronwall}) and the Theorem follows. Note the
importance of the term $-\pt\sn(t)$ on the right hand side of
(\ref{eq:entprod}).

\section{Tools}
\label{section:tools}

\subsection{Fixed time estimates}
\label{subs:fixedtime}

In the estimates with fixed time
$s\in[0,T]$ we shall use the notation
\begin{eqnarray}
\label{eq:bigeldef}
L=L(n):=n^{-2\b}l.
\end{eqnarray}
Note that
$L\gg1$  as
$n\to\infty$.

The following general entropy estimate will be exploited all over:

\begin{lemma}
\label{lemma:fixedtime}
{\rm (Fixed time entropy inequality)}
\\
Let
$l\le n$,
${\cal V}:\Om^l\to \R$ and denote
${\cal V}_j(\uome) := {\cal V}(\omega_{j},\dots,\omega_{j+l-1})$.
Then for any
$\gamma>0$
\begin{eqnarray}
\notag
&&
\hskip-10mm
\expect\Big(
\frac{1}{n}  \sum_{j\in\Tn}
{\cal V}_j (\Xn_s)
\Big)
\le
\frac{1}{\gamma} \hn(s)
+
\frac{1}{\gamma \,L}
\,
\frac{1}{n}\sum_{j\in\Tn}
\log\expect_{\nu^n_s}
\left(
\exp\big\{\gamma\,L\,{\cal V}_j\big\}
\right).
\\[-5pt]
&&
\label{eq:fixedtime}
\end{eqnarray}
\end{lemma}

\noindent
This lemma is standard tool in the context of relative entropy
method. For its proof we refer the reader to the original paper
\cite{yau1} or the monograph \cite{kipnislandim}.

\begin{proposition}
\label{propo:ldbounds}
{\rm (Fixed time large deviation bounds)}
\\
(i)
For any
$\varepsilon>0$ there exists
$M<\infty$ such that for any
$s\in[0,T]$
\begin{eqnarray}
\label{eq:ldboundrho+u}
\expect\Big(
\frac1n \sum_{j\in\Tn}
\big\{
\big(1+\rn + \abs{\un}\big)
\ind_{ \{ \rn \lor \abs{\un} > M \} }
\big\}
(s,\frac{j}{n})
\Big)
\le
\vareps\,\hn(s) + o(1).
\end{eqnarray}
(ii)
There exist
$C<\infty$ and
$M<\infty$ such that for any
$s\in[0,T]$
\begin{eqnarray}
\label{eq:ldboundusq}
\expect\Big(
\frac1n\sum_{j\in\Tn}
\big\{
\abs{\un}^2
\ind_{\{\rn\lor\abs{\un}>M\}}
\big\}(s,\frac{j}{n})
\Big)
\le
C\,\hn(s) + o(1).
\end{eqnarray}
\end{proposition}

The proof of Proposition
\ref{propo:ldbounds} is postponed to subsection
\ref{subs:ldboundsproof}. It relies on the entropy
inequality
(\ref{eq:fixedtime}) of Lemma
\ref{lemma:fixedtime}, the stochastic dominations formulated in Lemma
\ref{lemma:stochdom} (see subsection
\ref{subs:ldboundsproof}) and standard large deviation bounds.

We shall refer to
(\ref{eq:ldboundrho+u}) and
(\ref{eq:ldboundusq}) as
\emph{large deviation bounds}.

\begin{proposition}
\label{propo:kurschakbounds}
{\rm (Fixed time fluctuation bounds)}
\\
For any
$M<\infty$ there exists a
$C<\infty$ such that the following bounds hold:
\begin{eqnarray}
\label{eq:kurschakboundu}
&&
\expect\Big(
\frac1n \sum_{j\in\Tn}
\abs{\un-u}^2
(s,\frac{j}{n})
\Big)
\le
C\,\hn(s) + o(1),
\\[5pt]
\label{eq:kurschakboundrho}
&&
\expect\Big(
\frac1n \sum_{j\in\Tn}
\big\{
\abs{\rn-\rho}^2
\ind_{ \{ \rn \le  M \} }
\big\}
(s,\frac{j}{n})
\Big)
\le
C\,\hn(s) + o(1).
\end{eqnarray}
\end{proposition}

The proof of Proposition
\ref{propo:kurschakbounds} is postponed to subsection
\ref{subs:kurschakboundsproof}. It relies on the
entropy inequality
(\ref{eq:fixedtime}) of Lemma
\ref{lemma:fixedtime}, and  Gaussian fluctuation  estimates.

We shall refer to
(\ref{eq:kurschakboundu}) and
(\ref{eq:kurschakboundrho}) as
\emph{fluctuation bounds}.

\subsection{Convergence to local equilibrium and a priori bounds}
\label{subs:apriori}

The hydrodynamic limit relies on macroscopically fast convergence to
(local) equilibrium in blocks of mesoscopic size
$l$. Fix the block size
$l$ and
$(N,Z)\in\N\times(w_0/2)\Z$ with the restriction
$N\in[0, l\max\eta]$,
$Z\in[l\min\zeta, l\max\zeta]$ and denote
\begin{eqnarray*}
\Omega^l_{N,Z}
&:=&
\big\{
\uome\in\Omega^l\,:\,
\sum_{j=1}^l\eta_j=N, \sum_{j=1}^l\zeta_j=Z
\big\},
\\[5pt]
\pil_{N,Z}(\underline\omega)
&:=&
\pil_{\lambda,\theta}(\underline\omega\,|\, \sum_{j=1}^l\eta_j=N.
\sum_{j=1}^l\zeta_j=Z),
\end{eqnarray*}
Expectation with respect to the measure
$\pi^l_{N,Z}$ is denoted by
$\expect^l_{N,Z}\big(\cdot\big)$. For
$f:\Omega^l_{N,Z}\to\R$ let
\begin{eqnarray*}
&&
\hskip-3mm
K^l_{N,Z}f(\uome)
:=
\sum_{j=1}^{l-1}
\sum_{\ome',\ome''}
s(\ome_j,\ome_{j+1};\ome',\ome'')
\big(f(\Theta_{j,j+1}^{\ome',\ome''}\uome)-f(\uome)\big),
\\[5pt]
&&
\hskip-3mm
D^l_{N,Z}(f)
:=
\frac12
\sum_{j=1}^{l-1}
\expect^l_{N,Z}
\left(
\sum_{\ome',\ome''}
s(\ome_j,\ome_{j+1};\ome',\ome'')
\big(f(\Theta_{j,j+1}^{\ome',\ome''}\uome)-f(\uome)\big)^2
\right).
\end{eqnarray*}
In plain words:
$\Omega^l_{N,Z}$ is the hyperplane of configurations
$\uome\in\Omega^l$ with fixed values of the conserved quantities,
$\pi^l_{N,Z}$ is the
\emph{microcanonical distribution} on this hyperplane,
$K^l_{N,Z}$ is the symmetric infinitesimal generator restricted to the
hyperplane
$\Omega^l_{N,Z}$, and finally
$D^l_{N,Z}$ is the Dirichlet form associated to
$K^l_{N,Z}$. Note, that
$K^l_{N,Z}$ is defined with \emph{free boundary conditions}.

The convergence to local equilibrium is
\emph{quantitatively controlled} by the following uniform  logarithmic
Sobolev estimate, assumed to hold:

\begin{enumerate}[(A)]
\setcounter{enumi}{\value{aux}}

\item
\label{cond:lsi}
\emph{Logarithmic Sobolev inequality:}
There exists a finite constant
$\aleph$ such that for any
$l\in\N$,
$(N,Z)\in\N\times(w_0/2)\Z$ with  the restriction
$N\in[0, l\max\eta]$,
$Z\in[l\min\zeta, l\max\zeta]$, and any
$h:\Omega^l_{N,Z}\to\R_+$ with
$\expect^l_{N,Z}(h)=1$ the following bound holds:
\begin{eqnarray}
\label{eq:lsi}
\expect^l_{N,Z}\big(h\log h\big)
\le
\aleph \, l^2
D^l_{N,Z} \left(\sqrt{h}\right).
\end{eqnarray}

\setcounter{aux}{\value{enumi}}
\end{enumerate}

\medskip
\noindent {\bf Remark:}
The uniform logarithmic Sobolev inequality (\ref{eq:lsi}) is expected
to hold for a very wide range of locally finite interacting particle
systems, though we do not know about a fully general proof.
In \cite{yau2}   the logarithmic Sobolev
inequality  is proved  for  symmetric $K$-exclusion
processes.  This implies that (\ref{eq:lsi}) holds for the
two lane models defined in subsection \ref{section:models}. In
\cite{fritztoth} Yau's method of proving logarithmic Sobolev
inequality is applied and the logarithmic Sobolev inequality is
stated for random stirring models with arbitrary number of colors.
In particular, (\ref{eq:lsi}) follows for the $\{-1,0,+1\}$-model
defined in subsection \ref{section:models}.

\smallskip

The following  large deviation bound  goes back to Varadhan
\cite{varadhan}.  See also the monographs
\cite{kipnislandim}  and
\cite{fritz1}.

\begin{lemma}
\label{lemma:varadhan}
{\rm (Time-averaged entropy inequality, local equilibrium)}
\\
Let
$l\le n$,
${\cal V}:\Om^l\to \R_+$ and denote
${\cal V}_j(\uome) := {\cal V}(\omega_{j},\dots,\omega_{j+l-1})$.
Then for any
$\gamma>0$
\begin{eqnarray}
\label{eq:varadhan}
&&
\hskip-7mm
\expect\left(
\frac{1}{n}  \sum_{j\in\Tn} \int_0^t
{\cal V}_j (\Xn_s)\, ds
\right)
\le
\\[5pt]
\notag
&&
\hskip10mm
\frac{\aleph\,l^3}{2\,\gamma \,n^{1+3\b+\d}}\,
\Big(
\sn(t)
+
\frac{2\,n^{1+3\b+\d}\,t}{\aleph l^3}
\max_{ N, Z }
\log
\expect^{l}_{N,Z}
\big(\exp\left\{\gamma {\cal V} \right\}\big)
\Big).
\end{eqnarray}
\end{lemma}

\noindent
{\bf Remarks:}
(1)
Since
\[
\frac{n^{1+3\b+\d}}{l^3}
=o(1),
\]
in order to apply efficiently Lemma
\ref{lemma:varadhan} one has to chose
$\gamma=\gamma(n)$ so that
\[
\expect^{l}_{N,Z}
\big(\exp\left\{\gamma {\cal V} \right\}\big)
=\Ordo(1),
\]
\emph{uniformly} in the  block size
$l=l(n)\in\N$, and in
$N\in[0,l\max\eta]$ and
$Z\in[l\min\zeta,l\max\zeta]$.
\\
(2)
Assuming only uniform bound of size
$\sim \left(\aleph l^{2}\right)^{-1}$ on the spectral gap of
$K^l_{N,Z}$ (rather  than the stronger logarithmic Sobolev inequality
(\ref{eq:lsi})) and using Rayleigh-Schr\"odinger perturbation (see
Appendix 3 of
\cite{kipnislandim}) we would get
\begin{eqnarray*}
&&
\hskip-5mm
\expect\left(
\frac{1}{n}  \sum_{j\in\Tn} \int_0^t
{\cal V}_j (\Xn_s)\, ds
\right)
\le
\\[5pt]
&&
\hskip5mm
\frac{\aleph\,l^3\,\Vert {\cal V}\Vert_\infty }{2\,n^{1+3\b+\d}}
\sn(t)
+
t \, \Vert {\cal V}\Vert_\infty
\Big(
\frac{\max_{N,Z} \expect^l_{N,Z}\big({\cal V}\big)}
     {\Vert{\cal V}\Vert_\infty}
+
\frac{\max_{N,Z} \var^l_{N,Z}\big({\cal V}\big)}
     {\Vert{\cal V}\Vert_\infty^2}
\Big),
\end{eqnarray*}
which would not be sufficient for our needs.
\\
(3)
The proof of the bound
(\ref{eq:varadhan}) explicitly relies on the logarithmic Sobolev
inequality
(\ref{eq:lsi}). It appears  in
\cite{yau3} and it is reproduced in several places, see e.g.
\cite{fritz1}, \cite{fritz2}. We do  not repeat it here.

\smallskip

The main probabilistic ingredients of our proof are summarized in
Proposition
\ref{propo:aprioribounds} which is  consequence of Lemma
\ref{lemma:varadhan}. These are variants of the celebrated
\emph{one block estimate}, respectively,
\emph{two blocks estimate} of Varadhan and  co-authors.

\begin{proposition}
\label{propo:aprioribounds}
{\rm (Time-averaged block replacement and gradient bounds)}
\\
Given a local variable $\xi:\Ome^m\to\R$ there exists a constant $C$
such that the following bounds hold:
\\
(i)
\begin{equation}
\label{eq:plainobrepl}
\expect
\left(
\int_0^t\int_{\T}
\big|
\big\{
\wih{\xi}
-
\Xi(\wih{\eta},\wih{\zeta})
\big\}
(s,x)
\big|^2
\,dx\,ds
\right)
\le
C\,
\frac{l^2}{n^{1+3\b+\d}}
\,
\big(
\sn(t) + o(1)
\big).
\end{equation}
(ii)
\begin{equation}
\label{eq:plaingrad}
\expect
\left(
\int_0^t\int_{\T}
\big|
\px \wih{\xi}(s,x)
\big|^2
\,dx\,ds
\right)
\le
C
n^{1-3\b-\d}
\big(
\sn(t) + o(1)
\big).
\end{equation}
(iii)
Further on, if $\xi:\Ome\to\R$ (that is: it depends on a single spin)
and $\xi(\ome)=0$ whenever $\eta(\ome)=0$ then the following stronger
version of the gradient bound holds:
\begin{equation}
\label{eq:turbograd}
\expect
\left(
\int_0^t\int_{\T}
\frac{\big|
\px \wih{\xi}(s,x)
\big|^2}
{\wih{\eta}(s,x)}
\,dx\,ds
\right)
\le
C
n^{1-3\b-\d}
\big(
\sn(t) + o(1)
\big).
\end{equation}
\end{proposition}

The proof of Proposition
\ref{propo:aprioribounds} is postponed to subsection
\ref{subs:aprioriboundsproof}. It relies on the
large  deviation bound
(\ref{eq:varadhan}) and some elementary probability  estimates stated
in Lemma
\ref{lemma:epl} (see subsection
\ref{subs:aprioriboundsproof}).

We shall refer to
(\ref{eq:plainobrepl}), respectively,
(\ref{eq:plaingrad}) and (\ref{eq:turbograd}) as the
\emph{block replacement  bounds}, respectively, the
\emph{gradient  bounds}.

We shall apply (\ref{eq:plainobrepl}) to $\xi=\phi$ and
$\xi=\psi$. From  (\ref{eq:plaingrad}) it follows that
\begin{eqnarray}
\label{eq:plaingradu}
\expect
\left(
\int_0^t\int_{\T}
\big|
\px \wih{u}(s,x)
\big|^2
\,dx\,ds
\right)
\le
C
n^{1-\b-\d}
\big(
\sn(t) + o(1)
\big),
\\[5pt]
\label{eq:plaingradrho}
\expect
\left(
\int_0^t\int_{\T}
\big|
\px \wih{\rho}(s,x)
\big|^2
\,dx\,ds
\right)
\le
C
n^{1+\b-\d}
\big(
\sn(t) + o(1)
\big).
\end{eqnarray}
Using (\ref{eq:turbograd}) the last bound is improved to
\beq
\label{eq:turbogradrho}
\expect
\left(
\int_0^t\int_{\T}
\frac{
\big|
\px \wih{\rho}(s,x)
\big|^2}{\wih{\rho}(s,x)}
\,dx\,ds
\right)
\le
C
n^{1-\b-\d}
\big(
\sn(t) + o(1)
\big).
\eeq
The bound (\ref{eq:turbograd}) will also be applied to $\xi=\kappa$
(see (\ref{eq:kappa}) and (\ref{eq:kappavanish})) to get
\beq
\label{eq:turbogradkappa}
\expect
\left(
\int_0^t\int_{\T}
\frac{
\big|%
n^{2\b} \px \wih{\kappa}(s,x)
\big|^2}{\wih{\rho}(s,x)}%
\,dx\,ds
\right)
\le
C
n^{1-\b-\d}
\big(
\sn(t) + o(1)
\big).
\eeq

\section{Control of the large values of $(\rho,u)$: 
proof of (\ref{eq:firstterm})}
\label{section:firstterm}

\subsection{Preparations}
\label{subs:firsttermprep}

In the present section we prove
(\ref{eq:firstterm}). First we replace
$\frac1n\sum_{j\in\Tn}\cdots$ by
$\int_\T\cdots\,dx$. Note that given a smooth function
$F:\T\to\R$
\begin{eqnarray}
\label{eq:summinint}
\abs{
\frac1n\sum_{j\in\Tn} F(\frac{j}{n})
-
\int_\T F(x) \,dx
}
\le
\frac{1}{n} \left(\int_\T \abs{\px F(x)}^2\,dx\right)^{1/2}.
\end{eqnarray}
Hence it follows that
\begin{eqnarray}
\label{eq:sumtoint}
&&
\hskip-20mm
\expect
\left(
\int_0^t\frac1n\sum_{j\in\Tn}
\Big\{
\big(\px v\big)\,
\big(n^{3\beta}\wih{\psi}\big)\,
\J(\rn,\un)
\Big\}
(s,\frac{j}{n})\,ds
\right)
=
\\[5pt]
\notag
&&
\hskip1cm
\expect
\left(
\int_0^t\int_\T
\Big\{
\big(\px v\big)\,
\big(n^{3\beta}\wih{\psi}\big)\,
\J(\rn,\un)
\Big\}
(s,x)\,dx\,ds
\right)
+
A^n_{13},
\end{eqnarray}
where
$A^n_{13}$ is again a simple numerical error term:
\begin{eqnarray}
\notag
\abs{A^n_{13}}
&\le&
C n^{3\beta-1}
\Bigg\{
1+
\sqrt{
\expect
\left(
\int_0^s\int_\T
\abs{\px \wih{\psi}(s,x)}^2\,dx\,ds
\right)}
+
\sqrt{
\expect
\left(
\int_0^t\int_\T
\abs{\px \rn (s,x)}^2\,dx\,ds
\right)}
\\[5pt]
\notag
&&
\hskip18mm
+
\sqrt{
\expect
\left(
\int_0^t\int_\T
\abs{\px \un (s,x)}^2\,dx\,ds
\right)}
\,\Bigg\}
\\[5pt]
\label{eq:an13bound}
&=&
\Ordo(n^{5\beta} l^{-1})
=o(1).
\end{eqnarray}
In the last step we use the most straightforward gradient bound
(\ref{eq:trivigradbound}). (Using the gradient bound
(\ref{eq:plaingrad}) we  could obtain the much better upper bound
\begin{eqnarray*}
\abs{A^n_{13}}
=
\Ordo(n^{(-1-\delta+7\beta)/2})
=
o(n^{5\beta} l^{-1}),
\end{eqnarray*}
but we do not need this sharper estimate at this stage.)

So, we have to prove that the first term on the right hand side of
(\ref{eq:sumtoint}) is negligible. Recall that $\J=\S_\rho$. We
start with the application of the  martingale identity:
\begin{eqnarray}
\notag
&&
\hskip-8mm
\expect
\left(
\int_\T
\left\{
\big\{
v
\S(\rn,\un)
\big\}
(t,x)
-
\big\{
v
\S(\rn,\un)
\big\}
(0,x)
-
\int_0^t
\big\{
(\pt v)
\S(\rn,\un)
\big\}(
s,x)
\,ds
\right\}
\,dx
\right)
=
\\[5pt]
&&
\hskip30mm
\notag
\phantom{+\,\,}
\expect
\left(
\int_0^t
 \int_\T
v(s,x)
\left(
n^{1+\b}\Ln \S(\rn(x),\un(x))
\right)(\Xn_s)
\,dx\,ds
\right)
\\[-5pt]
\label{eq:1term_elso}
&&
\\[-5pt]
\notag
&&
\hskip30mm
\phantom{}
+
\expect
\left(
\int_0^t
 \int_\T
v(s,x)
\left(
n^{1+\b+\d}\Kn \S(\rn(x),\un(x))
\right)(\Xn_s)
\,dx\,ds
\right)
\end{eqnarray}
%
\subsection{The left hand side of (\ref{eq:1term_elso})}
\label{subs:firsttermleft}

From
(\ref{eq:sn'rhobound}),
(\ref{eq:sn'ubound}), and
(\ref{eq:ulrchoice}),
we conclude that
\[
\abs{\S(\rho,u)}\le
C\, \big(\rho+\abs{u}\big) \, \ind_{\{\rho\lor\abs{u}>M\}}.
\]
Hence, using  the large deviation bound
(\ref{eq:ldboundrho+u}) it follows that, by choosing
$M$ sufficiently large  we obtain
\begin{eqnarray*}
\expect\Big(
\frac1n \sum_{j\in\Tn}
\abs{\Sn(\rn,\un)(s,\frac{j}{n})}
\Big)
\le
\vareps\,\hn(s) + o(1).
\end{eqnarray*}
Hence, applying again
(\ref{eq:summinint})  we get
\begin{eqnarray}
\label{eq:lhsof1term_elso}
\abs{\,\text{l.h.s. of
(\ref{eq:1term_elso})}\,}
\le
\frac12\,\hn(t)
+
C\,\int_0^t\hn(s)\,ds
+
o(1).
\end{eqnarray}
{\bf Remark:}
Note that this is the point where
$M$ and thus the lower edge of the cutoff is fixed.
Also note the importance of the factor $1/2$ in fromt of $\hn(t)$ on
the right hand side.

\subsection{The right hand side of (\ref{eq:1term_elso}):
first computations}
\label{subs:firsttermright1}

First we compute how the infinitesimal generators
$n^{1+\b}\Ln$ and
$n^{1+\b+\d}\Kn$ act on the function
$\uome\mapsto \S(\rn(x),\un(x))$:
\begin{eqnarray}
\label{eq:lns}
&&
\hskip-8mm
n^{1+\b}\Ln
\S(\rn(x),\un(x))
=
\\[5pt]
\notag
&&
\hskip11.5mm
\Big\{
\S_{\rho}(\rn,\un)
\big(n^{3\b}\px\wih{\psi}\big)
+
\S_{u}(\rn,\un)
\big(n^{2\b}\px\wih{\phi}\big)
\Big\}(x)
+
A^n_{14}(x),
\\[10pt]
\label{eq:kns}
&&
\hskip-8mm
n^{1+\b+\d}\Kn
\S(\rn(x),\un(x))
=
\\[5pt]
\notag
&&
\hskip0mm
n^{-1+\b+\d}
\Big\{
\S_{\rho}(\rn,\un)
\big(n^{2\b}\px^2\wih{\kappa}\big)
+
\S_{u}(\rn,\un)
\big(n^\b\px^2\wih{\chi}\big)
\Big\}(x)
+
A^n_{15}(x),
\end{eqnarray}
where
$A^n_{14}(x)$ and
$A^n_{15}(x)$ are the following
\emph{numerical error terms}:
\begin{eqnarray*}
&&
A^n_{14}(x)
=
A^n_{14}(\uome,x)
:=
n^{1+\b}
\sum_{j\in\Tn}
\sum_{\ome', \ome'\in \Omega}
r(\ome_{j},\ome_{j+1};\ome',\ome')
\times
\\[5pt]
&&
\hskip20mm
\Big\{
\Sn
\big(
\rn(x)
+
\frac{n^{2\b}}{l}
\big(a\big(\frac{nx-j}{l}\big)-a\big(\frac{nx-j-1}{l}\big)\big)
(\eta'-\eta_{j}),
\\[5pt]
&&
\hskip28.5mm
\un(x)
+
\frac{n^{\b}}{l}
\big(a\big(\frac{nx-j}{l}\big)-a\big(\frac{nx-j-1}{l}\big)\big)
(\zeta'-\zeta_{j})
\big)
\\[5pt]
&&
\hskip23mm
-
\Sn\big(\rn(x),\un(x)\big)
\\[5pt]
&&
\hskip23mm
-
\Sn_{\rho}\big(\rn(x),\un(x)\big)
\frac{n^{2\b}}{l^2}
a^{\prime}\big(\frac{nx-j}{l}\big)(\eta'-\eta_{j})
\\[5pt]
&&
\hskip23mm
-
\Sn_{u}\big(\rn(x),\un(x)\big)
\frac{n^{\b}}{l^2}
a^{\prime}\big(\frac{nx-j}{l}\big)(\zet'-\zet_{j})
\Big\}.
\end{eqnarray*}
\begin{eqnarray*}
&&
A^n_{15}(x)
=
A^n_{15}(\uome,x)
:=
n^{1+\b+\d}
\sum_{j\in\Tn}
\sum_{\ome', \ome'\in \Om}
s(\ome_{j},\ome_{j+1};\ome',\ome')
\times
\\[5pt]
&&
\hskip20mm
\Big\{
\Sn
\big(
\rn(x)
+
\frac{n^{2\b}}{l}
\big(a\big(\frac{nx-j}{l}\big)-a\big(\frac{nx-j-1}{l}\big)\big)
(\eta'-\eta_{j}),
\\[5pt]
&&
\hskip28.5mm
\un(x)
+
\frac{n^{\b}}{l}
\big(a\big(\frac{nx-j}{l}\big)-a\big(\frac{nx-j-1}{l}\big)\big)
(\zeta'-\zeta_{j})
\big)
\\[5pt]
&&
\hskip23mm
-
\Sn\big(\rn(x),\un(x)\big)
\\[5pt]
&&
\hskip23mm
-
\Sn_{\rho}\big(\rn(x),\un(x)\big)
\frac{n^{2\b}}{l^2}
a^{\prime}\big(\frac{nx-j}{l}\big)(\eta'-\eta_{j})
\\[5pt]
&&
\hskip23mm
-
\Sn_{u}\big(\rn(x),\un(x)\big)
\frac{n^{\b}}{l^2}
a^{\prime}\big(\frac{nx-j}{l}\big)(\zeta'-\zet_{j})
\Big\}
\\[5pt]
&&
\hskip15mm
+
\Sn_{\rho}\big(\rn(x),\un(x)\big)
\times
\\[5pt]
&&
\hskip2cm
\frac{n^{1+3\b+\d}}{l^2}
\sum_{j\in\Tn}
\Big\{
a^{\prime}\big(\frac{nx-j}{l}\big)(\kap_{j+1}-\kap_{j})
+
\frac{1}{l}
a^{\prime\prime}\big(\frac{nx-j}{l}\big)\kap_{j}
\Big\}
\\[5pt]
&&
\hskip15mm
+
\Sn_{u}\big(\rn(x),\un(x)\big)
\times
\\[5pt]
&&
\hskip2cm
\frac{n^{1+2\b+\d}}{l^2}
\sum_{j\in\Tn}
\Big\{
a^{\prime}\big(\frac{nx-j}{l}\big)(\chi_{j+1}-\chi_{j})
+
\frac{1}{l}
a^{\prime\prime}\big(\frac{nx-j}{l}\big)\chi_{j}
\Big\}.
\end{eqnarray*}
These error terms are easily estimated: using the fact that the second
partial derivatives of
$\S$ are uniformly bounded and
$\zeta$ and
$\eta$ are bounded,  by simple Taylor expansion after tedious but
otherwise straightforward computations we find:
\begin{eqnarray}
\label{eq:a14}
&&
\sup_{\uome\in\Omn}\sup_{x\in\T}
\abs{A^n_{14}(\uome,x)}
\le
C n^{1+3\b}l^{-2}
=
o(1),
\\[5pt]
\label{eq:a15}
&&
\sup_{\uome\in\Omn}\sup_{x\in\T}
\abs{A^n_{15}(\uome,x)}
\le
C n^{1+5\b+\d}l^{-3}
=
o(1).
\end{eqnarray}
No probabilistic arguments are involved in these bounds. The global
(averaged and integrated)  error introduced by these terms will be of
the same order.

Next we do some further transformations on the main terms coming
from the right hand sides of (\ref{eq:lns}) and (\ref{eq:kns}).
Performing  integrations by part, introducing the macroscopic
fluxes and using (\ref{eq:lentfluxpde_n}) we obtain:
\begin{eqnarray}
\notag
\label{eq:lnsipart}
&&
\hskip-18mm
-
\int_\T
v(x)
\Big\{
\S_{\rho}(\rn,\un)
\px\big(n^{3\b}\wih{\psi}\big)
+
\S_{u}(\rn,\un)
\px\big(n^{2\b}\wih{\phi}\big)
\Big\}(x)
\,dx
=
\\[5pt]
\notag
&&
\hskip15mm
\phantom{+\,}
\int_\T
\px v(x)
\Big\{
\big(n^{3\b}\wih{\psi}\big)
\S_{\rho}(\rn,\un)
\Big\}(x)
\,dx
\\[5pt]
\notag
&&
\hskip15mm
+
\int_\T
\px v(x)
\Big\{
\F(\rn,\un)
-
\Psi^{n}(\rn,\un) \S_{\rho}(\rn,\un)
\Big\}(x)
\,dx
\\[5pt]%
\notag
&&%
\hskip15mm
+%
\int_\T
\px v(x)
\Big\{
n^{2\beta}
\S_{u}(\rn,\un)
\big(\wih{\phi}-\Phi(\wih{\eta},\wih{\zeta})\big)
\Big\}(x)
\,dx
\\[5pt]%
\notag
&&%
\hskip15mm
+
\int_\T
 v(x)
\Big\{
\hskip2mm
n^{3\beta}
\S_{\rho\rho}(\rn,\un)
\big(\px\rn\big)
\big(\wih{\psi}-\Psi(\wih{\eta},\wih{\zeta})\big)
\\
\notag
&&
\hskip32mm
+
n^{3\beta}
\S_{\rho u}(\rn,\un)
\big(\px\un\big)
\big(\wih{\psi}-\Psi(\wih{\eta},\wih{\zeta})\big)
\\[5pt]
\notag
&&
\hskip32mm
+
n^{2\beta}
\S_{u\rho}(\rn,\un)
\big(\px\rn\big)
\big(\wih{\phi}-\Phi(\wih{\eta},\wih{\zeta})\big)
\\[5pt]
\notag
&&
\hskip32mm
+
n^{2\beta}
\S_{uu}(\rn,\un)
\big(\px\un\big)
\big(\wih{\phi}-\Phi(\wih{\eta},\wih{\zeta})\big)
\Big\}(x)
\,dx
\end{eqnarray}
Note that, since
$\J=\S_\rho$, the first term on the right hand side is exactly the
expression in the main term on the right hand side of
(\ref{eq:sumtoint}). Estimating the other terms on the right hand
side of
(\ref{eq:lnsipart}) is the object of the next subsection.

Now we turn to the main term on the right hand side of
(\ref{eq:kns}). Here, straightforward integration by parts yields
\begin{eqnarray}
\label{eq:knsipart}
&&
\hskip-18mm
-
\int_\T
v(x)
\Big\{
\S_{\rho}(\rn,\un)
\big(n^{2\b}\px^2\wih{\kappa}\big)
+
\S_{u}(\rn,\un)
\big(n^{\b}\px^2\wih{\chi}\big)
\Big\}(x)
\,dx
=
\\[5pt]
\notag
&&
\phantom{+\,}
\int_\T
\px v(x)
\Big\{
\S_{\rho}(\rn,\un)
\big(n^{2\b}\px\wih{\kappa}\big)
+
\S_{u}(\rn,\un)
\big(n^{\b}\px\wih{\chi}\big)
\Big\}(x)
\,dx
\\[5pt]
\notag
&&
+
\int_\T
v(x)
\Big\{
\phantom{+\,\,\,}
\S_{\rho\rho}(\rn,\un)
\big(\px\rn\big) \big(n^{2\b}\px\wih{\kappa}\big)
\\[5pt]
\notag
&&
\phantom{
+
\int_\T
v(x)
\Big\{
}
+
\S_{\rho u}(\rn,\un)
\left(
\big(\px\un\big) \big(n^{2\b}\px\wih{\kappa}\big)
+
\big(\px\rn\big) \big(n^{\b}\px\wih{\chi}\big)
\right)
\\[5pt]
\notag
&&
\phantom{
+
\int_\T
v(x)
\Big\{
}
+
\S_{uu}(\rn,\un)
\big(\px\un\big) \big(n^{\b}\px\wih{\chi}\big)
\Big\}(x)
\,dx
\end{eqnarray}
We will estimate the terms emerging from
the right hand side in the next subsection.

\subsection{The right hand side of (\ref{eq:1term_elso}):
bounds}
\label{subs:firsttermbounds}

\subsubsection{}
\label{subsubs:f-psis_rho}
We note that
\[
\abs{
\Fn(\rn,\un)-\Psi^n(\rn,\un)\Sn_\rho(\rn,\un)}
\le
C\big(1+\abs{\un}^2\big)
\ind_{\{\rn\lor\abs{\un}>M\}},
\]
see (\ref{eq:fnbound}) and
(\ref{eq:ulrchoice}).
Hence, applying the large deviation bounds (\ref{eq:ldboundrho+u})
and (\ref{eq:ldboundusq}) we obtain
\begin{eqnarray}
\label{eq:firstbound}
&&
\hskip-7mm
\expect\left(
\int_\T
\abs{
\Big\{
\F(\rn,\un)
-
\Psi^{n}(\rn,\un) \S_{\rho}(\rn,\un)
\Big\}(s,x)
}
\,dx
\right)
\\[5pt]
\notag
&&
\hskip75mm
\le
C\, \hn(s)
+
o(1).
\end{eqnarray}
%

\subsubsection{}
\label{subsubs:blockrepl}
We use
\[
\abs{\S_{u}(\rn,\un)}
\le C,
\]
see
(\ref{eq:sn'ubound}) and the first block replacement bound
(\ref{eq:plainobrepl}) to  obtain:
\begin{eqnarray}
\label{eq:secondbound}
&&
\hskip-7mm
\expect\left(
\int_0^t\int_\T
\abs{
\Big\{
n^{2\beta}
\S_{u}(\rn,\un)
\big(\wih{\phi}-\Phi(\wih{\eta},\wih{\zeta})\big)
\Big\}(s,x)
}
\,dx\,ds
\right)
\\[5pt]
\notag
&&
\hskip65mm
\le
C \,l \, n^{(-1-\d+\b)/2}
=o(1).
\end{eqnarray}
%

\subsubsection{}
\label{subsubs:gradandblock}

Next we use
\begin{eqnarray}
\notag
&&
\hskip-8mm
\abs{\S_{\rho\rho}(\rn,\un)}
\le
\frac{C}{\log(\olr/\ulr)}
\frac{1}{\ulr+\rn},
\qquad
\abs{\S_{\rho u}(\rn,\un)} \le
\frac{C}{\log(\olr/\ulr)}
\frac{1}{\sqrt{\ulr}+\sqrt{\rn}}
\\
&&
\hskip27mm
\abs{\S_{u u}(\rn,\un)} \le
\frac{C}{\log(\olr/\ulr)},
\label{eq:secondderbounds}
\end{eqnarray}
see (\ref{eq:sn''rhorhobound}), (\ref{eq:sn''rhoubound}),
respectively, (\ref{eq:sn''uubound}), and note that here \emph{we
do not exploit} the fact that the constant factors on the right
hand side are actually small. These, together with the   block
replacement bounds (\ref{eq:plainobrepl}), the gradient  bounds
(\ref{eq:plaingradu}), (\ref{eq:turbogradrho}) and the bound
(\ref{eq:apriorientropybound}) on the relative entropy $s^n(t)$
yield  the following four estimates:
\begin{eqnarray}
\notag
&&
\hskip-7mm
\expect\left(
\int_0^t\int_\T
\abs{
\Big\{
n^{3\beta}
\S_{\rho\rho}(\rn,\un)
\big(\px\rn\big)
\big(\wih{\psi}-\Psi(\wih{\eta},\wih{\zeta})\big)
\Big\}(s,x)
}
\,dx\,ds
\right)
\le
C\,l\,n^{\b-\d}
=o(1),
\\[5pt]
\notag
&&
\hskip-7mm
\expect\left(
\int_0^t\int_\T
\abs{
\Big\{
n^{3\beta}
\S_{\rho u}(\rn,\un)
\big(\px\un\big)
\big(\wih{\psi}-\Psi(\wih{\eta},\wih{\zeta})\big)
\Big\}(s,x)
}
\,dx\,ds
\right)
\le
C\,l\,n^{\b-\d}
=o(1),
\\[5pt]
\label{eq:thirdbound}
&&
\hskip-7mm
\expect\left(
\int_0^t\int_\T
\abs{
\Big\{
n^{2\beta}
\S_{u\rho}(\rn,\un)
\big(\px\rn\big)
\big(\wih{\phi}-\Phi(\wih{\eta},\wih{\zeta})\big)
\Big\}(s,x)
}
\,dx\,ds
\right)
\le C\,l\,n^{-\d} =o(1),
\\[5pt]
\notag
&&
\hskip-7mm
\expect\left(
\int_0^t\int_\T
\abs{
\Big\{
n^{2\beta}
\S_{uu}(\rn,\un)
\big(\px\un\big)
\big(\wih{\phi}-\Phi(\wih{\eta},\wih{\zeta})\big)
\Big\}(s,x)
}
\,dx\,ds
\right)
\le
C\,l\,n^{-\d}
=o(1).
\end{eqnarray}
%

\subsubsection{}
\label{subsubs:symmetric1}
Using
\[
\abs{\S_\rho(\rn,\un)}
\le C,
\qquad
\abs{\S_u(\rn,\un)}
\le C,
\]
see (\ref{eq:sn'rhobound}) and (\ref{eq:sn'ubound}), and the
gradient bounds (\ref{eq:plaingradu}), (\ref{eq:plaingradrho}) we
obtain the following two bounds:
\begin{eqnarray}
\notag
&&
\hskip-7mm
n^{-1+\b+\d}
\expect\left(
\int_0^t\int_\T
\abs{
\Big\{
\S_{\rho}(\rn,\un)
\big(\px\rn\big)
\Big\}(s,x)
}
\,dx\,ds
\right)
\\[5pt]
\notag
&&
\hskip69mm
\le
C\,n^{(-1+\d+3\b)/2}
=o(1),
\\[-5pt]
&&
\label{eq:fourthbound}
\\[-5pt]
\notag
&&
\hskip-7mm
n^{-1+\b+\d}
\expect\left(
\int_0^t\int_\T
\abs{
\Big\{
\S_{u}(\rn,\un)
\big(\px\un\big)
\Big\}(s,x)
}
\,dx\,ds
\right)
\\[5pt]
\notag
&&
\hskip69mm
\le
C\,n^{(-1+\d+\b)/2}
=o(1).
\end{eqnarray}
%

\subsubsection{}
\label{subsubs:symmetric2}

The following bounds are of paramount importance and they 
are sharp.
We use (\ref{eq:secondderbounds}) again and note
 that here we exploit it in its
\emph{full  power}: the constant factor on 
the right hand side is small.
These and the gradient bounds
(\ref{eq:plaingradu}) and
(\ref{eq:turbogradrho}) yield the following three bounds:
\begin{eqnarray}
\notag
&&
\hskip-7mm
n^{-1+\b+\d}
\expect\left(
\int_0^t\int_\T
\abs{
\Big\{
\S_{\rho\rho}(\rn,\un)
\big(\px\rn\big) \big(n^{2\b}\px\wih{\kappa}\big)
\Big\}(s,x)
}
\,dx\,ds
\right)
\\[3pt]
\notag
&&
\hskip87mm
\le
c \, \sn(t)
+
o(1),
\\[2pt]
\notag
&&
\hskip-7mm
n^{-1+\b+\d}
\expect\left(
\int_0^t\int_\T
\abs{
\Big\{
\S_{\rho u}(\rn,\un)
\big(\px\un\big)\big(n^{2\b}\px\wih{\kappa}\big)
\Big\}(s,x)
}
\,dx\,ds
\right)
\\[3pt]
\notag
&&
\hskip87mm
\le
c \,\sn(t)
+
o(1),
\\[2pt]
\label{eq:fifthbound}
&&
\hskip-7mm
n^{-1+\b+\d}
\expect\left(
\int_0^t\int_\T
\abs{
\Big\{
\S_{\rho u}(\rn,\un)
\big(\px\rn\big)\big(n^{\b}\px\wih{\chi}\big)
\Big\}(s,x)
}
\,dx\,ds
\right)
\\[3pt]
\notag
&&
\hskip87mm
\le
c \,\sn(t)
+
o(1),
\\[2pt]
\notag
&&
\hskip-7mm
n^{-1+\b+\d}
\expect\left(
\int_0^t\int_\T
\abs{
\Big\{
\S_{u u}(\rn,\un)
\big(\px\un\big)\big(n^{\b}\px\wih{\chi}\big)
\Big\}(s,x)
}
\,dx\,ds
\right)
\\[3pt]
\notag
&&
\hskip87mm
\le
c\, \sn(t)
+
o(1).
\end{eqnarray}
The ratio $\olr/\ulr$ is chosen so large that
\beq
\label{olrchoice2}
c \sup_{(t,x)\in[0,T]\times\T}\abs{v(t,x)}<\frac12.
\eeq
\subsection{Sumup}
\label{subs:firstermright2}

The identities
(\ref{eq:lns}),
(\ref{eq:kns}),
(\ref{eq:lnsipart}),
(\ref{eq:knsipart}) and the bounds
(\ref{eq:a14}),
(\ref{eq:a15}),
(\ref{eq:firstbound}),
(\ref{eq:secondbound}),
(\ref{eq:thirdbound}),
(\ref{eq:fourthbound}),
(\ref{eq:fifthbound}) yield
\begin{eqnarray}
&&\hskip-18mm
\notag
\abs{
\expect
\left(
\int_0^t\int_\T
\Big\{
\big(\px v\big)\,
\big(n^{3\beta}\wih{\psi}\big)\,
\S_\rho(\rn,\un)
\Big\}
(s,x)\,dx\,ds
\right)
-
\Big(
\text{r.h.s. of (\ref{eq:1term_elso})}
\Big)
}
\\[5pt]
&&\hskip60mm
\label{eq:rhsof1term_elso}
\le
\frac12\,\sn(t)
+
c\,\int_0^t\hn(s)\,ds
+
o(1).
\end{eqnarray}
Finally, from
(\ref{eq:sumtoint}),
(\ref{eq:an13bound}),
(\ref{eq:1term_elso}),
(\ref{eq:lhsof1term_elso}) and
(\ref{eq:rhsof1term_elso}) we obtain
(\ref{eq:firstterm}).

\section{Control of the small values of $(\rho,u)$: 
proof of the bounds (\ref{eq:secondterm}) to
  (\ref{eq:fifthterm}) }
\label{section:otherterms}

\subsection{Proof of (\ref{eq:secondterm})}
\label{subs:secondterm}

We exploit the straightforward inequality
\begin{eqnarray*}
\abs{\Jn(\rn,\un)}
=
\abs{\Sn_\rho(\rn,\un)}
\le
C\,
\ind_{\{\rn\lor\abs{\un}>M\}},
\end{eqnarray*}
see
(\ref{eq:sn'rhobound}) and
(\ref{eq:ulrchoice}),
and boundedness of the functions
$\rho(t,x)$,
$u(t,x)$,
$\px v(t,x)$. Thus, applying the large deviation bound
(\ref{eq:ldboundrho+u}) we readily obtain
(\ref{eq:secondterm}).

\subsection{Proof of (\ref{eq:thirdterm})}
\label{subs:thirdterm}

This is very similar to what has been done in various parts of
subsection
\ref{subs:firsttermbounds}. We use  the block replacement
bound
(\ref{eq:plainobrepl}) and the bound
\begin{eqnarray}
\label{eq:ibound}
\quad
\abs{\In(\rn,\un)}
=
\abs{1-\Sn_\rho(\rn,\un)}
\le
C
\quad
\end{eqnarray}
which follows from (\ref{eq:sn'rhobound}).
%
%
We readily obtain
\begin{eqnarray*}
&&
\hskip-8mm
\expect
\left(
\int_0^t\int_\T
\abs{
\Big\{
n^{3\beta}
\big(
\wih{\psi} - \Psi(\wih{\eta},\wih{\zeta})
\big)\,
\big\vert\I(\rn,\un)\big\vert
\Big\}
(s,x)}\,ds\,dx
\right)
\\
&&
\hskip60mm
\le
C\,l\,n^{(-1-\d+3\b)/2}
=o(1),
\end{eqnarray*}
which proves
(\ref{eq:thirdterm}).


\subsection{Proof of (\ref{eq:fourthterm})}
\label{subs:fourthterm}

We write
\begin{eqnarray}
\label{eq:idecomp}
\In(\rn,\un)=
\ind_{\{\rn\lor\abs{\un}\le M\}}
+
\ind_{\{\rn\lor\abs{\un}> M\}}
\In(\rn,\un),
\end{eqnarray}
and note that, by Taylor expansion of the function
$(\rho,u)\mapsto\Psi(\rho,u)$
\begin{eqnarray*}
\abs{\Psi^n(\rn,\un) - \rn \un}
\ind_{\{\rn\lor\abs{\un}\le M\}}
\le
C\,n^{-2\b}.
\end{eqnarray*}
On the other hand
\begin{eqnarray*}
\abs{\Psi^n(\rn,\un)}
\le
C\,\rn\,\abs{\un}
\end{eqnarray*}
and
\begin{eqnarray}
\label{eq:ibound'}
\rn\,  \abs{\In(\rn,\un)}
\le C (1+\abs{\un}),
\end{eqnarray}
see
(\ref{eq:D1D2bounds}) and (\ref{eq:sn'rhobound}). Thus
\begin{eqnarray*}
\abs{\Psi^n(\rn,\un) - \rn \un}
\abs{\In(\rn,\un)}
\le
C\,
\Big(
n^{-2\b}
+
\big(\abs{\un}+\abs{\un}^2\big)
\ind_{\{\rn\lor\abs{\un}> M\}}
\Big).
\end{eqnarray*}
From this,  using the large deviation bounds
(\ref{eq:ldboundrho+u}) and
(\ref{eq:ldboundusq}) we obtain
(\ref{eq:fourthterm}).

\subsection{Proof of (\ref{eq:fifthterm})}
\label{subs:fifthterm}

We use again
(\ref{eq:idecomp})  and
(\ref{eq:ibound'})  and get
\begin{eqnarray*}
\abs{
\big(\rn-\rho\big)
\big(\un-u\big)\,
\I(\rn,\un)}
\le
\phantom{+\,}
\abs{
\big(\rn-\rho\big)
\big(\un-u\big)}
\ind_{\{\rn\lor\abs{\un}\le M\}}
\\[5pt]
+
C\left(1+\abs{\un}+ \abs{\un}^2\right)
\ind_{\{\rn\lor\abs{\un}> M\}}
\end{eqnarray*}
Now the fluctuation bounds
(\ref{eq:kurschakboundu}),
(\ref{eq:kurschakboundrho}), and the large deviation bounds
(\ref{eq:ldboundrho+u}),
(\ref{eq:ldboundusq}) together yield
(\ref{eq:fifthterm}).

\section{Construction of the cutoff function: proofs}
\label{section:cutoffproofs}

\subsection{Proof of Lemma \ref{lemma:charcurves}}
\label{subs:charcurvesproof}

\begin{proof}
We sketch the proof for $u\ge0$ and leave the very similar $u\le0$
case for the reader. Let $\rho_1>0$, $u_1>0$  be so chosen that
for $(\rho,u)\in[0,\rho_1]\times[0,u_1]$  the following bounds
hold with a fixed $c>0$:
\beqs && \abs{\Phi_u(\rho,u)-\Psi_\rho(\rho,u)-(2\g-1)u} \le c u
(u^2+\rho),
\\[5pt]
&& \abs{ \Phi_\rho(\rho,u)-1}\le c (u^2+\rho),
\\[5pt]
&& \abs{\Psi_u(\rho,u)-\rho}\le c \rho (u^2+\rho), \eeqs and \beqs
\Phi_u(\rho,u)-\Psi_\rho(\rho,u)\not=0 \quad \text{ for }
(\rho,u)\not=(0,0). \eeqs This can be done due to the
$(\rho,u)\to(0,0)$ asymptotics of the macroscopic fluxes $\Phi$
and $\Psi$.

It follows that as long as
$(\sigma(u;r),u)\in[0,\rho_1]\times[0,u_1]$
\beqs%
\frac{d \rho}{d u} \le \frac{2\rho
(1+c'(\rho+u^2))}{\sqrt{(2\g-1)^2
u^2+4\rho}+(2\g-1)u}. \eeqs%
This implies
\[
\sigma(u;r)\le r+C \left(\sqrt{r} u \wedge r^{\frac{4\g-3}{4\g-2}}
u^{\frac1{2\g-1}}\right)
\]
with a positive C, as long as
$(\sigma(u;r),u)\in[0,\rho_1]\times[0,u_1]$.

From our assumptions it also follows that for $\rho\le\rho_1$ and
$u>u_1$
\[
\frac{d \rho}{d u} \le b\,\rho,
\]
where
\[
b:= \sup_{\substack{\rho<\rho_1\\u> u_1}} \frac{2\Psi_u(\rho,u)}
{\rho\,\big(\Phi_u(\rho,u)-\Psi_\rho(\rho,u)\big)} < \infty.
\]
Hence it follows that for $u\ge u_1$
\[
\sigma(u;r) \le \sigma(u_1;r) \exp\{b(u-u_1)\}
\]
as long as $\sigma(u;r)\le\rho_1$.

Putting these two arguments together the upper bound
\beqs
\sigma(u;r)
\le
r+C_1 \left(\sqrt{r} u
\wedge
r^{\frac{4\g-3}{4\g-2}} u^{\frac1{2\g-1}}\right),
\,\,\,
u\ge0,
\,\,\,
r<r_0,
\eeqs
follows with
\[
r_0
:=
\sup\left\{r:
\left( r + C \left(\sqrt{r} u_1
       \wedge
       r^{\frac{4\g-3}{4\g-2}} u_1^{\frac1{2\g-1}}\right)
\right)
\exp\{b(u^*-u_1)\}\le\rho_1\right\},
\]
and
\[
C_1
=
\frac {r_0^{\frac1{4\g-2}}\big(\exp\{b(u^*-u_1)\}-1\big)}
{u_1} + C\exp\{b(u^*-u_1)\}.
\]
\end{proof}

\subsection{Proof of Lemma \ref{lemma:lentbounds}}
\label{subs:lentboundsproof}

Note first that given the bounds (\ref{eq:macr_flx_asy}) and
condition (\ref{cond:phi0'-psi1}) of subsection
\ref{subs:genprop}, (\ref{eq:s''uubound}) follows directly from
(\ref{eq:s''rhoubound}), (\ref{eq:s''rhorhobound}) and
(\ref{eq:lentpde}). So, we shall concentrate on
(\ref{eq:s'rhobound})-(\ref{eq:s''rhoubound}) only.

By differentiating in the pde (\ref{eq:lentpde}) and applying
straightforward transformations we obtain the following
differential equations for $S_{\rho}$, $S_{u}$, $S_{\rho\rho}$,
and $S_{\rho u}$, respectively:
\beq
\label{eq:s'rhopde}
&& \hskip-25mm
\Psi_{u} {(S_{\rho})}_{\rho\rho} +
\big(\Phi_{u}-\Psi_{\rho}\big) {(S_{\rho})}_{\rho u} - \Phi_{\rho}
{(S_{\rho})}_{u u}+
\\
\notag
&& \hskip5mm
+ \Phi_{\rho}
\left(\frac{\Psi_{u}}{\Phi_{\rho}}\right)_{\rho}
{(S_{\rho})}_{\rho} + \Phi_{\rho}
\left(\frac{\Phi_{u}-\Psi_{\rho}}{\Phi_{\rho}}\right)_{\rho}
{(S_{\rho})}_{u} =0,
\eeq
\beq
\label{eq:s'upde}
&&
\hskip-25mm
\Psi_{u} {(S_{u})}_{\rho\rho} +
\big(\Phi_{u}-\Psi_{\rho}\big) {(S_{u})}_{\rho u} - \Phi_{\rho}
{(S_{u})}_{u u}+
\\
\notag
&&
\hskip5mm
+ \Psi_{u}
\left(\frac{\Phi_{u}-\Psi_{\rho}}{\Psi_{u}}\right)_{u}
{(S_{u})}_{\rho} - \Psi_{u}
\left(\frac{\Phi_{\rho}}{\Psi_{u}}\right)_{u} {(S_{u})}_{u} =0,
\eeq
\beq
&&
\hskip-25mm
\Psi_{u} {(S_{\rho\rho})}_{\rho\rho} +
\big(\Phi_{u}-\Psi_{\rho}\big) {(S_{\rho\rho})}_{\rho u} -
\Phi_{\rho} {(S_{\rho\rho})}_{u u}+
\label{eq:s''rhorhopde}
\\
\notag
&&
\hskip5mm
+ 2 \Phi_{\rho}
\left(\frac{\Psi_{u}}{\Phi_{\rho}}\right)_{\rho}
{(S_{\rho\rho})}_{\rho} + 2 \Phi_{\rho}
\left(\frac{\Phi_{u}-\Psi_{\rho}}{\Phi_{\rho}}\right)_{\rho}
{(S_{\rho\rho})}_{u} - \Phi_{\rho}
\left(\frac{\Psi_{u}}{\Phi_{\rho}}\right)_{\rho\rho} S_{\rho\rho}
\\
\notag
&&
\hskip75mm = - \Phi_{\rho}
\left(\frac{\Phi_{u}-\Psi_{\rho}}{\Phi_{\rho}}\right)_{\rho\rho}
S_{\rho u},
\eeq
\beq
\label{eq:s''rhoupde}
&&
\hskip-25mm
\Psi_{u} {(S_{\rho u})}_{\rho\rho} +
\big(\Phi_{u}-\Psi_{\rho}\big) {(S_{\rho u})}_{\rho u} -
\Phi_{\rho}%
{(S_{\rho u})}_{u u}+
\\
\notag
&&
\hskip5mm
+ \left\{ \Psi_{u}
\left(\frac{\Phi_{u}-\Psi_{\rho}}{\Psi_{u}}\right)_{u} +
\Phi_{\rho} \left(\frac{\Psi_{u}}{\Phi_{\rho}}\right)_{\rho}
\right\} {(S_{\rho u})}_{\rho}+
\\
\notag
&&
\hskip5mm
+ \left\{ \Phi_{\rho}
\left(\frac{\Phi_{u}-\Psi_{\rho}}{\Phi_{\rho}}\right)_{\rho} -
\Psi_{u} \left(\frac{\Phi_{\rho}}{\Psi_{u}}\right)_{u} \right\}
{(S_{\rho u})}_{u}+
\\
\notag
&&
\hskip5mm
+ \left\{ \Psi_{u}
\left(\frac{\Phi_{\rho}}{\Psi_{u}}\right)_{u}
\left(\frac{\Phi_{u}-\Psi_{\rho}}{\Phi_{\rho}}\right)_{\rho} +
\Phi_{\rho}
\left(\frac{\Phi_{u}-\Psi_{\rho}}{\Phi_{\rho}}\right)_{\rho u}
\right\} S_{\rho u}
\\
\notag
&&
\hskip40mm
=
- \left\{ \Psi_{u}
\left(\frac{\Phi_{\rho}}{\Psi_{u}}\right)_{u}
\left(\frac{\Psi_{u}}{\Phi_{\rho}}\right)_{\rho} + \Phi_{\rho}
\left(\frac{\Psi_{u}}{\Phi_{\rho}}\right)_{\rho u} \right\}
S_{\rho \rho}.
\eeq %
Because of the conditions outlined in subsection
\ref{subs:genprop} all the coefficients are smooth functions, if
$\rho$ is small enough.

The respective initial conditions are
\beq
 \label{eq:s'rhoic}&&
\hskip-10mm
S_{\rho}(\rho,0)= \frac{\log(r/\ulr)}{\log(\olr/\ulr)}
\ind_{\{\rho\in[\ulr,\olr]\}} +
\ind_{\{\rho\in[\olr,r_0]\}},%
\quad {(S_{\rho})}_{u}(\rho,0)= 0,
\\[5pt]
\label{eq:s'uic} &&
\hskip-10mm
 S_{u}(\rho,0)= 0, \quad {(S_{u})}_u(\rho,0)=
\frac{1}{\log(\olr/\ulr)} \,
\frac{\Psi_u(\rho,0)}{\rho \, \Phi_{\rho}(\rho,0)}
\ind_{\{\rho\in[\ulr,\olr]\}},
\\[5pt]
\label{eq:s''rhorhoic} &&
\hskip-10mm
 S_{\rho\rho}(\rho,0)=
\frac{\rho^{-1}}{\log(\olr/\ulr)} \ind_{\{\rho\in[\ulr,\olr]\}},
\quad {(S_{\rho\rho})}_u(\rho,0)= 0,
\\[5pt]
 &&\hskip-10mm
 S_{\rho u}(\rho,0)= 0, \notag\\%
&&\hskip-10mm
{(S_{\rho u})}_{u}(\rho,0)=
\frac{1}{\log(\olr/\ulr)}%
\left(
\frac{\Psi_u(\rho,0)}{\rho \, \Phi_{\rho}(\rho,0)}
\right)_{\rho}
+\frac{1}{\log(\olr/\ulr)}\,
\frac{\Psi_{u}(\rho,0)}{\rho\,\Phi_\rho(\rho,0)}
\left\{\delta(\rho-\ulr)-\delta(\rho-\olr)\right\}.
\label{eq:s''rhouic}
\eeq%
Observe, that because of the asymptotics (\ref{eq:macr_flx_asy})
we have
\beq%
\frac{\Psi_u(\rho,0)}{\rho\,\Phi_{\rho}(\rho,0)}=1+\Ordo(\rho),\quad
{\left(
\frac{\Psi_u(\rho,0)}{\rho\,\Phi_{\rho}(\rho,0)}
\right)}_\rho(\rho,0)=\Ordo(1) 
\label{eq:boundary_asy}.
\eeq%
In order to understand the pdes
(\ref{eq:s'rhopde})-(\ref{eq:s''rhoupde}) first we analyze in
general the pde \beq \label{eq:genpde} \Psi_{u} f_{\rho\rho} +
\big(\Phi_{u}-\Psi_{\rho}\big) f_{\rho u} - \Phi_{\rho} f_{u u} +
A f_{\rho} + B f_{u} + G f = 0, \eeq the functions
$\{A\,,\,B\,,\,G\}=\{A\,,\,B\,,\,G\}(\rho,u)$ being given on the
left hand side of the pdes
(\ref{eq:s'rhopde})-(\ref{eq:s''rhoupde}). It is easy to check
that $A$ and $G$ are even, $B$ is odd with respect to $u$ and also
that $\kappa:=A(0,0)$ is $1,2\g-1,2,2\g$, respectively, in the
four cases.

We solve in $\wt{\cal D}$ the Cauchy problem (\ref{eq:genpde})
with the initial condition \beq \label{eq:geninicond}
f(\rho,0)=s(\rho), \qquad f_u(\rho,0)=t(\rho), \qquad
\rho\in[0,r_0]. \eeq The functions $s(\rho)$ and  $t(\rho)$ will
be identified with the various expressions in
(\ref{eq:s'rhoic})-(\ref{eq:s''rhouic}). Then, in ${\cal
D}_3(\olr)\setminus\wt{\cal D}$ we solve the  Goursat problem
(\ref{eq:genpde}) with boundary conditions \beq
\label{eq:genboundcond} f(\rho,u)= \left\{
\begin{array}{cl}
f(r_1,0)%
& \text{ on } \quad \partial{\cal D}_3(\olr)\setminus\wt{\cal D},
\\[8pt]
\begin{array}{l}
\text{given by the solution of}
\\
\text{the previous Cauchy problem}
\end{array}
& \text{ on } \quad \partial\wt{\cal D}\cap {\cal D}_3(\olr).
\end{array}
\right. \eeq
Mind that $f(r_1,0)=1$ in the case of (\ref{eq:s'rhopde}) and
$f(r_1,0)=0$ in the cases (\ref{eq:s'upde})-(\ref{eq:s''rhoupde}).

The pde (\ref{eq:genpde}) is hyperbolic in the domains considered.
Its Jacobian matrix is
\begin{eqnarray}
\label{eq:jacobian} D = D(\rho,u) := \left(
\begin{array}{cc}
\ \Psi_\rho & \Psi_u
\\[8pt]
\Phi_\rho & \Phi_u
\end{array}
\right).
\end{eqnarray}
The eigenvalues of $D(\rho,u)$ are
\begin{eqnarray}
\label{eq:eigenvalues} 
\left.
\begin{array}{c}
\lambda\\\mu
\end{array}
\right\} = \pm \frac12 \left\{ \sqrt{\big(\Phi_u-\Psi_\rho\big)^2
+ 4\Phi_\rho\Psi_u} \mp \big(\Phi_u-\Psi_\rho\big) \right\}
+\Phi_u
\end{eqnarray}
Mind that from the Onsager relation (\ref{eq:onsager}) it follows
that for any $(\rho,u)\in{\cal D}$ $\big(\Phi_u-\Psi_\rho\big)^2 +
4\Phi_\rho\Psi_u\ge0$.

The characteristic coordinates (or Riemann invariants)
$w=w(\rho,u)$, $z=z(\rho,u)$ of the pde (\ref{eq:genpde})  are
determined, up to a functional relation \beq \label{eq:wzgauge}
\wt w(\rho,u) = g(w(\rho,u)), \qquad \wt z(\rho,u) = h(z(\rho,u)),
\eeq by the eigenvalue equations
\begin{eqnarray}\label{eq:riemanninv_eq}
(w_\rho,w_u) D= \lambda \, (w_\rho,w_u) , \qquad (z_\rho,z_u) D=
\mu \, (z_\rho,z_u).
\end{eqnarray}
Due to the gauge invariance (\ref{eq:wzgauge}) we can choose the
characteristic coordinates $w$ and $z$ so that \beqs w(0,u)=u
\ind_{\{u>0\}}, \qquad z(0,u)=-u \ind_{\{u<0\}}. \eeqs This choice
determines uniquely the characteristic coordinates. Observe, that
as a corollary of Lemma
\ref{lemma:charcurves} we have %
\beq%
\label{eq:riem_asy1} c_1 \left(\sqrt{\rho}+u\right)&<&w\rruu< c_2
\left(\sqrt{\rho}+u\right),\\\label{eq:riem_asy2} \rho&<&c\,
z\rruu^{\frac{4\g-3}{2\g-1}}\, w\rruu^{\frac1{2\g-1}}.
\eeq%
We denote
\beqs%
 & w^\circ:=w(r_0,0)=z(r_0,0)=:z^\circ &
\\[5pt]
& \ulw:=w(\ulr,0)=z(\ulr,0)=:\ulz &
\\[5pt]
& \olw:=w(\olr,0)=z(\olr,0)=:\olz &
\eeqs%
In characteristic coordinates %
\beqs && \wt{\cal D} = [0,w^\circ]\times[0,z^\circ],
\\[5pt]
&& \wt{\cal D}\cap\{u=0\} =
[0,w^\circ]\times[0,z^\circ]\cap\{w=z\},
\\[5pt]
&& {\cal D}_3(\olr)\setminus\wt{\cal D} =
[w^\circ,u_*]\times[0,\olz]\cup[0,\olw]\times[z^\circ,u_*],
\\[5pt]
&&
\partial {\cal D}_3(\olr)\setminus\wt{\cal D}=
\{(w,\olz):w\in[w^\circ,u_*]\} \cup \{(\olw,z):z\in[z^\circ,u_*]\}
\\[5pt]
&&
\partial \wt {\cal D}\cap {\cal D}_3(\olr)=
\{(w^\circ,z):z\in[0,z^\circ]\} \cup
\{(w,z^\circ):w\in[0,w^\circ]\}. \eeqs The pde (\ref{eq:genpde})
written in characteristic coordinates reads
\beq %
\label{eq:genchpde} f_{wz}+\alpha f_w +\beta f_z +\nu f=0,
\eeq%
where \beq &&\notag \alpha=\alpha(w,z):= \frac{\lambda_z - A u_z +
B \rho_z}{\lambda-\mu}
\\[5pt]
&& \beta=\beta(w,z):= -\frac{\mu_w  - A u_w +  B
\rho_w}{\lambda-\mu}\label{eq:coeff_ch_def}
\\[5pt]
&&\notag \nu=\nu(w,z):= \frac{G (\rho_w u_z -\rho_z u_w)
}{\lambda-\mu} \eeq (Now all the functions on the right are
understood as functions of $(w,z)$.) In characteristic coordinates
the initial conditions of the Cauchy problem in $\wt{\cal D}$ are
\beq \label{eq:genchinicond} \notag && f(v,v) = s(\rho(v,v)) =:
\wt s(v),
\\[-5pt]
\\[-5pt]
\notag && \big(f_w-f_z)(v,v)
= 2u_w(v,v)\, t(\rho(v,v)) =: \wt t(v). \eeq%
The Cauchy problem (\ref{eq:genchpde})+(\ref{eq:genchinicond}) in
the domain $\wt{\cal D}\cap\{z\le w\}$ is solved by \beq
\label{eq:riemannsolver} f(w_0,z_0) &=& \phantom{+} \frac12
\varphi(z_0,z_0) \wt s(z_0) + \frac12 \varphi(w_0,w_0) \wt s(w_0)
\\[5pt]
\notag && + \frac12 \int_{z_0}^{w_0} \varphi(v,v) \wt t(v)  dv +
\frac12 \int_{z_0}^{w_0} \big(\varphi_w-\varphi_z\big)(v,v) \wt
s(v) dv
\\[5pt]
\notag && + \int_{z_0}^{w_0}
\varphi(v,v)(\beta(v,v)-\alpha(v,v))\wt s(v)dv \eeq where the
\emph{Riemann function} $\varphi$ is a solution of the adjoint
Goursat problem: \beq \label{eq:gen_hyp_pde_adj}
\varphi_{wz}-(\alp \varphi)_w-(\bet \varphi)_z+\nu \varphi=0, \eeq
with boundary conditions: \beq \label{eq:gen_adj_bound} \left\{
\begin{array}{l}
\varphi(w_0,t)=\exp \int_{z_0}^t \alp(w_0,v)dv, \quad
t\in[z_0,w_0],
\\[5pt]
\varphi(s,z_0)=\exp \int_{w_0}^s \bet(v,z_0)dv, \quad\,
s\in[z_0,w_0].
\end{array}
\right. \eeq %
Actually, the Riemann function  depends also on $(w_0,z_0)$:
$\varphi(w,z)= \varphi(w_0,z_0;w,z)$. In order to avoid heavy
typography  we omit explicit notation of this dependence. Note
that in our cases (because of the left-right reflection
 symmetry of the  respective pde), we will have
 $\bet(v,v)=\alpha(v,v)$, thus  on the right hand side of
 (\ref{eq:riemannsolver}) the last term cancels.

\noindent Also, if we consider the non-homogeneous pde
\[
h_{wz}+\a h_w+\b h_z+\nu h=g
\]
with the same initial conditions then it is solved by %
\beq h(w_0,z_0)=f(w_0,z_0)+{\int\int}_{\bigtriangleup_{w_0,z_0}}
g(s,t) \varphi(s,t)ds dt, \label{eq:genchpde_inh_solve}
\eeq%
where $\bigtriangleup_{w_0,z_0}$ is the triangle with vertices
$(z_0,z_0), (w_0,w_0), (w_0,z_0)$.

For details see any advanced textbook   on partial differential
equations, e.g. \cite{john}, \cite{garabedian}, \cite{evans},
\cite{smoller}.

In order to estimate $f$ we will give a uniform estimate on the
Riemann function $\varphi$.

\begin{proposition}[Bounds on the Riemann
    function]\label{prop:bound_rsolver}
Let $\varphi$ be the Riemann function associated to the equation
(\ref{eq:genpde}) (where the coefficients $A,B,G$ are given on the
left hand side of (\ref{eq:s'rhopde}), (\ref{eq:s'upde}),
(\ref{eq:s''rhorhopde}) or (\ref{eq:s''rhoupde})) and $w_0>z_0>0$.
Then the following bounds hold uniformly  for $(w_0,z_0)\in
\wt{\cal D}$ with $z_0<t<s<w_0$
and $(s,t)\in \wt{\cal D}$:%
\beq\label{eq:entr_rsolv_bounds_wz}%
\abs{\varphi(s,t)}&<&c
\left(\frac{s}{w_0}\right)^{\frac{\kappa-1}{2
\g-1}}\\%
\notag \abs{(\pp_w\varphi-\pp_z\varphi)(s,s)}&<& c \frac1{w_0}
\left(\frac{s}{w_0}\right)^{\frac{\kappa-1}{2 \g-1}-1}.
\eeq%
\end{proposition}

Using Proposition \ref{prop:bound_rsolver} with
(\ref{eq:riemannsolver}), the initial conditions
(\ref{eq:s'rhoic})-(\ref{eq:s''rhouic}) and with Lemma
\ref{lemma:charcurves} we can estimate $f$ in $\wt{\cal D}$ which
gives (\ref{eq:s'rhobound}) and (\ref{eq:s'ubound}) in this
domain.

The equations (\ref{eq:s''rhorhopde}) and (\ref{eq:s''rhoupde})
are not closed for $S_{\rho\rho}$ and $S_{\rho u}$, respectively,
with the previous method one can only prove the required estimates
for the solution of the respective homogeneous pdes. However, with
(\ref{eq:genchpde_inh_solve}) it is easy to show that these
estimates can be extended for $S_{\rho\rho}$ and $S_{\rho u}$,
too.

As an example, we show how to get the bound
(\ref{eq:s''rhorhobound}) for the \emph{homogeneous} solution $f$
of (\ref{eq:s''rhorhopde}) with initial conditions
(\ref{eq:s''rhorhoic}).

For the initial conditions we have the following bounds (using
(\ref{eq:riem_asy1}),(\ref{eq:riem_asy2})): %
\beqs \abs{\wt s(v)}&<& \frac{c}{\log(\olr/\ulr) v^2}
\ind_{\{v\in[\ulz,\olz]\}},\\
\abs{\wt t(v)}&=&0.\eeqs From Proposition \ref{prop:bound_rsolver}
we have that for any $z_0<w_0$
\[
\abs{\varphi(v,v)}< c \left(\frac{v}{w_0}\right)^{\frac1{2\g-1}},
\quad \abs{\partial_w \varphi(v,v)-\partial_z \varphi(v,v)}<
\frac{c}{w_0} \left(\frac{v}{w_0}\right)^{\frac1{2\g-1}-1}.
\]
Together with (\ref{eq:riemannsolver}) we get
\beq%
\label{eq:f_bound} \abs{f(w_0,z_0)}< \frac{c}{\log(\olr/\ulr)}
w_0^{-\frac1{2\g-1}} \ulz^{2-\frac1{2\g-1}}
\ind_{\{z_0\in[\ulz,\olz], \ulz\leq w_0\}},
\eeq%
for any $z_0<w_0$. This means
\[
\abs{f(w_0,z_0)}< \frac{c}{\log(\olr/\ulr)} \ulz^{2}
\ind_{\{z_0\in[\ulz,\olz], \ulz\leq w_0\}},
\]
which (using (\ref{eq:riem_asy1})) translates to the $\rruu$
coordinates the following way:
\[
\abs{f\rruu}< \frac{c}{\log(\olr/\ulr)} \frac1{\ulr} \ind_{{\cal
D}_3(\ulr,\olr)}(\rho,u).
\]
(We only get this for $u\ge 0$, but from the symmetry of the pde
this is also true for $u<0$.) Also from (\ref{eq:f_bound}) we get
\[
\abs{f(w_0,z_0)}< \frac{c}{\log(\olr/\ulr)} w_0^{-\frac1{2\g-1}}
z_0^{2-\frac1{2\g-1}} \ind_{\{z_0\in[\ulz,\olz], \ulz\leq w_0\}},
\]
which (using (\ref{eq:riem_asy2})) gives
\[
\abs{f\rruu}< \frac{c}{\log(\olr/\ulr)} \frac1{\rho} \ind_{{\cal
D}_3(\ulr,\olr)}(\rho,u).
\]
Putting together the two bounds on $f$ we get the required
estimate of (\ref{eq:s''rhorhobound}).

Now back to the proof of Proposition \ref{prop:bound_rsolver}. The
following lemma will be a basic tool for our estimates:
\begin{lemma}[Goursat estimate]\label{lem:goursat}
Suppose the functions $A(w,z), B(w,z),C(w,z)$ are defined on
$[x_1,x_2]\times[y_1,y_2]$ where $0\leq y_1<y_2\leq x_1<x_2 \leq
\infty$ and have the following
properties:%
\beq%
\sup_{z\in[y_1,y_2]}\int_{x_1}^{x_2} \abs{B(s,z)}ds &<&\frac16,\notag\\
\sup_{w\in[x_1,x_2]}\int_{y_1}^{y_2} \abs{A(w,t)}dt &<&\frac16
,\label{eq:coeff_int_bounds}\\
\int_{x_1}^{x_2}\int_{y_1}^{y_2}\abs{C(s,t)} ds dt
&<&\frac16.\notag
\eeq%
Let $U$ be a solution of%
\beq U_{wz}-(A U)_w-(B U)_z+ C U=0\label{eq:pde_gen}
\eeq%
on the rectangle $[x_1,x_2]\times[y_1,y_2]$ and let
\[
\sup_{y_1\leq z\leq y_2}\abs{U(x_1,z)}+\sup_{x_1\leq w \leq
x_2}\abs{U(w,y_2)}=M.
\]
Then
\[
\sup_{(w,z)\in [x_1,x_2]\times[y_1,y_2]} \abs{U(w,z)}\leq 5 M.
\]
\end{lemma}

\begin{proof}
Denote $U(w,y_2)=f(w)$, $U(x_1,z)=g(z)$. Denote %
\beqs%
\wt{f}(w)&=&f(w)-\int_{x_1}^{w} B(s,y_2)f(s)ds,\\
\wt{g}(z)&=&g(z)-\int_{y_2}^{z} A(x_1,t)g(t)dt.
\eeqs%
Then $U$ satisfies the following integral equation
(for every $(w,z)\in [x_1,x_2]\times[y_1,y_2]$): %
\beq\notag%
U(w,z)&=&\wt{f}(w)+\wt{g}(z)-U(x_1,y_2)
+\int_{y_2}^{w} A(w,t)
U(w,t)dt+\int_{x_1}^{z} B(s,z) U(s,z) ds\notag
\\&&-\int_{x_1}^{w}\int_{y_2}^{z} C(s,t) U(s,t)
dsdt.\label{eq:gen_integrodiff}
\eeq%
Taking absolute values after some trivial estimates we get:
\[
\abs{U(w,z)}\leq \frac{13}{6}M+\frac12  \sup_{(s,t)\in
[x_1,w]\times[z,y_2]} \abs{U(s,t)},
\]
from that the needed bound follows immediately.
\end{proof}

\begin{remark*} Observe, that if we have some additional
estimates on $A_z,B_w$ and $C$ then differentiating
(\ref{eq:gen_integrodiff}) with respect to $w$ or $z$ and using
the Grönwall inequality, one could also get bounds on $\abs{U_w}$
and $\abs{U_z}$ in $[x_1,x_2]\times[y_1,y_2]$. \end{remark*}

\begin{proof}[Proof of Proposition \ref{prop:bound_rsolver}]
$\phantom{ }$\\
\noindent(i) We first prove the proposition in the case when
$\Psi(\rho,u)=\rho u$ and $\Phi(\rho,u)=\rho+\g u^2$. In that case
the coefficient functions in (\ref{eq:genpde}) take the following
form: $A=\kappa$, $B=C=0$. Also, the Riemann invariants can be
computed explicitly:%
\beq
\label{def:rieminv_inf}
\begin{array}{rcl}
w(\rho,u)
&:=&
\left(
\frac{\sqrt{(2 \gamma-1) u^2+4\rho}+(2\gamma -1)u}
     {4\g-2}
\right)^{\frac{2\g-1}{4\g-3}}
\left(
\sqrt{(2 \gamma-1) u^2+4 \rho}-(2\gamma -2)u\right)
          ^{\frac{2\gamma-2}{4\gamma-3}}
\\
z(\rho,u)
&:=&
\left(
\frac{\sqrt{(2 \gamma-1) u^2+4 \rho}-(2\gamma-1)u}
     {4\g-2}
\right)^{\frac{2\g-1}{4\g-3}}
\left(
\sqrt{(2 \gamma-1)u^2+4 \rho}+(2\gamma -2)u\right)
              ^{\frac{2\gamma-2}{4\gamma-3}}.
\end{array}
\eeq
One can easily check that
the equations (\ref{eq:riemanninv_eq}) hold.
\noindent We define%
\beq\label{def:rsolv_transf}%
\b_0(w):=\b(w,0)\quad \textup{and} \quad U(w,z):= \varphi(w,z)
\exp\left(-\int_{w_0}^{w} \b_0(s) ds\right) \eeq%
(for $0\leq z <w$). From the definitions one can calculate that
\[
\b_0(w)=\frac{\kappa-1}{2\g-1} \frac{u_w
(w,0)}{u(w,0)}=\frac{\kappa-1}{2\g-1} w^{-1}.
\]
Thus, it is enough to prove, that%
\beq\label{eq:entr_rsolv_bounds_aim}%
\abs{U(s,t)}<C \textup{,}\quad \abs{\pp_w U(s,s)}<C
\frac1{s}\textup{,}\quad \abs{\pp_z U(s,s)}<C
\eeq%
for $z_0\leq t \leq s \leq w_0$ uniformly, with a constant
depending only on $\kappa$. To show
(\ref{eq:entr_rsolv_bounds_aim}) we will apply Lemma
\ref{lem:goursat}.

From (\ref{eq:gen_hyp_pde_adj}),
(\ref{eq:gen_adj_bound}) and (\ref{def:rsolv_transf}) we get that%
\beq%
U_{wz}-(\alp U)_w-\left((\b-\b_0) U\right)_z-\alpha \b_0 U=0,
\eeq%
and%
\beqs%
U(w_0,z)&=&\exp\left(\int_{z_0}^z \alpha(w_0,t) dt\right),\\
U(w,z_0)&=&\exp\left(\int_{w_0}^w (\b(s,z_0)-\b(s,0))ds\right). %
\eeqs%
Using the explicit formulas for the Riemann-invariants one can get
estimates for the integrals needed for Lemma \ref{lem:goursat}.
Suppose $[x_1,x_2]\times[y_1,y_2]\subseteq [v,w_0]\times[v,z_0]$
for some
$z_0<v<w_0$. Then it can be shown that for $z_0 \leq z\leq w \leq w_0$ %
\beq\notag%
\abs{\int_{y_1}^{y_2} \alpha(w,t) dt}&<& c \frac{x_2-x_1}{w}, \\
\abs{\int_{x_1}^{x_2} \b(s,z)-\b(s,0) ds}&<& c z
\left(\frac1{x_1}-\frac1{x_2}\right), \label{eq:coeff_scaling_3}\\
\abs{\int_{x_1}^{x_2} \int_{y_1}^{y_2} (\alpha \b_0)(s,t)
ds\,dt}&<& c\left(\frac1{x_1}-\frac1{x_2}\right) (y_2-y_1),\notag
\eeq%
where $c$ only depends on $\kappa$. From the first and second
inequality it follows that $\abs{U(w_0,z)}$ and $\abs{U(w,z_0)}$
can be bounded by a constant depending only on $\kappa$. Now fix
$s,t$ with $z_0<t<s<w_0$. Using (\ref{eq:coeff_scaling_3}) one can
partition $[v,w_0]\times[z_0,v]$ into smaller rectangles in a way
that the number of rectangles only depends on the value of
$\kappa$ and on each small rectangle the conditions of Lemma
\ref{lem:goursat} hold (with $A=\alpha$, $B=\b-\b_0$, $C=-\alpha
\b_0$). Applying successively Lemma \ref{lem:goursat} for the
small rectangles (starting with the vertex $(w_0,z_0)$) one gets
$\abs{U(s,t)}<c(\kappa)$. Using the remark after Lemma
\ref{lem:goursat} one can also get the results of
(\ref{eq:entr_rsolv_bounds_aim}) for the partial derivatives of
$U$.

\noindent(ii) In the general case we do not know the explicit
forms of the coefficients in (\ref{eq:genpde}) only their
asymptotics:
\beqs%
A(\rho,u)&=&\kappa (1+\Ordo(\rho+u^2))\\ B(\rho,u)&=&c_1 u
(1+\Ordo(\rho+u^2))\\ C(\rho,u)&=&c_2 (1+\Ordo(\rho+u^2)). \eeqs
We also do not have explicit formulas for the Riemann-invariants,
but because $\Psi\rruu=\rho u(1+\Ordo(\rho+u^2))$ and
$\Phi\rruu=(\rho+\g u^2)(1+\Ordo(\rho+u^2))$ if $\rho, \abs{u}\ll
1$ the level-lines will approximate the respective level-lines of
the system examined in (i). We will follow the steps of the proof
for the specific case. We define $\b_0$ and $U$ as in
(\ref{def:rsolv_transf}). From the asymptotics we have
$\b_0(w)=\frac{\kappa-1}{2\g-1} w^{-1}+\Ordo(1)$ which means it is
again enough to prove (\ref{eq:entr_rsolv_bounds_aim}). We have
the following equation for $U$:
\beq%
U_{wz}-(\alp U)_w-\left((\b-\b_0) U\right)_z+(\nu-\alpha \b_0)
U=0,
\eeq%
with%
\beqs%
U(w_0,z)&=&\exp\left(\int_{z_0}^z \alpha(w_0,t) dt\right),\\
U(w,z_0)&=&\exp\left(\int_{w_0}^w (\b(s,z_0)-\b(s,0))ds\right). %
\eeqs%
If we can prove similar bounds for the integrals of coefficients
as in (\ref{eq:coeff_scaling_3}) then using Lemma
\ref{lem:goursat} the required estimates follow. From
(\ref{eq:coeff_ch_def}) we have that%
\beqs %
\int_{x_1}^{x_2} \int_{y_1}^{y_2} \nu(s,t) ds\,dt=
\int\int_{\Delta} \frac{G}{\lambda-\mu} d\rho \,du \eeqs%
where $\Delta$ is the domain corresponding to the
$[x_1,x_2]\times[y_1,y_2]$ rectangle on the $\rruu$ plane. If
$r_0$ is small enough, then for $\rruu \in \wt{\cal D}$
\[\abs{\frac{C}{\lambda-\mu}}<c \frac{1}{\sqrt{u^2+\rho}}\] and since the
right-hand side is integrable this gives uniform bounds on the
previous integral. To get bounds on the integrals%
\beqs \int_{y_1}^{y_2} \alpha(w,t) dt,\,\,\int_{x_1}^{x_2}
\b(s,z)-\b(s,0) ds,\,\, \int_{x_1}^{x_2} \int_{y_1}^{y_2} (\alpha
\b_0)(s,t) ds\,dt,
\eeqs%
we observe that because of the asymptotics described earlier if
$r_0\rightarrow 0$ these integrals will be uniformly close to the
respective integrals of the (i) case. Thus, again if $r_0$ is
small enough (but fixed!) then the arguments of (i) may be
repeated.
\end{proof}
We have proved that the bounds
(\ref{eq:s'rhobound})-(\ref{eq:s''rhoubound}) hold in $\wt{\cal
D}$ if $r_0$ is small enough, now we have to extend this to the
domain ${\cal D}_3(\olr)\setminus\wt{\cal D}$ for the solution of
the respective Goursat-problems with boundary conditions
(\ref{eq:genboundcond}). The solution of these Goursat-problems
can be expressed by integral equations similar to
(\ref{eq:gen_integrodiff}) (see \cite{john}, \cite{garabedian}).
Thus if we have some estimates for the solution on the boundary
(which we have from (\ref{eq:genboundcond}) and the previous
estimates of the Cauchy problem) then these may be extended (up to
a constant multiplier) if we have uniform bounds on the integrals
of the respective coefficients.

If $\olr$ is small enough then for $\rruu \in {\cal
D}_3(\olr)\setminus\wt{\cal D}$ we have
\[\abs{\Phi_{\rho}\rruu-\Psi_u\rruu}>c>0\] which also implies that in
this domain $\lambda-\mu>c'>0$. From this it follows, that we can
fix a small enough $\olr_0>0$ such that in the domain ${\cal
D}_3(\olr_0)\setminus\wt{\cal D}$ all the respective coefficients
are well-defined smooth functions. Moreover, since this domain is
compact, they are all bounded with a fixed constant which means
that we have uniform bounds on the respective integrals. Using
similar arguments as in the estimate of the solution of the
Cauchy-problem one can get the bounds
(\ref{eq:s'rhobound})-(\ref{eq:s''rhoubound}) also in this domain
which completes the proof of Lemma \ref{lemma:lentbounds}.

\subsection{Proof of Lemma \ref{lemma:lentfluxbound}}
\label{subs:lentfluxboundproof}
\begin{proof}
Note that
\beqs
&&
\big(F-\Psi S_\rho\big)_\rho=
\Phi_\rho S_u-\Psi S_{\rho\rho},
\\[5pt]
&&
\big(F-\Psi S_\rho\big)_u=
\Phi_u S_u-\Psi S_{\rho u}.
\eeqs
Also, there exists a constant
$C<\infty$ such that for any
$(\rho,u)\in{\cal D}$
\beqs
|\Psi(\rho,u)|\le C\rho |u|,
\quad
|\Phi_\rho(\rho,u)|\le C,
\quad
|\Phi_u(\rho,u)|\le C|u|.
\eeqs
From these and the bounds
(\ref{eq:s'ubound}),
(\ref{eq:s''rhorhobound}),
(\ref{eq:s''rhoubound}) of Lemma
\ref{lemma:lentbounds} it follows that
\beqs
&&
\big|\big(F-\Psi S_\rho\big)_\rho\big|(\rho,u)
\le
\frac{C}{\log(\olr/\ulr)}\big(\sqrt{\olr}+|u|\big)
\,
\ind_{{\cal D}_3(\ulr,\olr)}(\rho,u),
\\[5pt]
&&
\big|\big(F-\Psi S_\rho\big)_u\big|(\rho,u)
\le
\frac{C}{\log(\olr/\ulr)}\big(\rho+\sqrt{\olr}|u|\big)
\,
\ind_{{\cal D}_3(\ulr,\olr)}(\rho,u).
\eeqs
Integrating these and  using
(\ref{eq:D1D2bounds}), the bound
(\ref{eq:fbound}) follows.
\end{proof}

\section{Proof of the ``Tools''}
\label{section:ingred}

\subsection{Proof of the large deviation bounds (Proposition
  \ref{propo:ldbounds})}
\label{subs:ldboundsproof}

Recall the definition (\ref{eq:bigeldef}) of $L$. The following
lemma follows from simple coupling arguments.

\begin{lemma}
\label{lemma:stochdom}
{\rm (Stochastic dominations)}
\\
There exists a constant
$C$ depending only on
$\max_{(s,x)\in[0,T]\times\T} \rho(s,x)$ and
$\max_{(s,x)\in[0,T]\times\T} \abs{u(s,x)}$ such that   for any fixed
$(s,x)\in[0,T]\times\T$ the following stochastic dominations hold:
\begin{eqnarray}
\label{eq:stochboundrho}
&&
\prob_{\nun_s}\Big(
\rn(x)>z
\Big)
\le
\prob\Big(
{\mathrm{POI}(L)}>(z/C)L
\Big),
\\[5pt]
\label{eq:stochboundu}
&&
\prob_{\nun_s}\Big(
\abs{\un(x)}>z
\Big)
\le
\prob\Big(
\abs{\mathrm{GAU}}>
\big((z/C)-1\big){\sqrt{L}}
\Big),
\end{eqnarray}
where
$\mathrm{POI}(L)$ is a Poissonian random variable with expectation
$L$, and
$\mathrm{GAU}$ is a standard Gaussian random variable.
\end{lemma}



\begin{lemma}
\label{lemma:largedevi}
{\rm (Large deviation bounds)}
\\
(i)
For any
$\gamma<\infty$ there exists
$M<\infty$, such that for any
$n$,
$j\in\Tn$ and
$s\in[0,T]$
\begin{eqnarray}
\notag
&&
\log\expect_{\nun_s}
\Big(
\exp\big\{
\gamma L \rn(\frac{j}{n}) \ind_{\{\rn(\frac{j}{n})>M\}}
\big\}
\Big)
\le1,
\\[5pt]
\notag
&&
\log\expect_{\nun_s}
\Big(
\exp\big\{
\gamma L \rn(\frac{j}{n}) \ind_{\{|\un(\frac{j}{n})|>M\}}
\big\}
\Big)
\le1,
\\[-5pt]
\label{eq:largedevi1}
\\[-5pt]
\notag
&&
\log\expect_{\nun_s}
\Big(
\exp\big\{
\gamma L \big|\un(\frac{j}{n})\big| \ind_{\{\rn(\frac{j}{n})>M\}}
\big\}
\Big)
\le1,
\\[5pt]
\notag
&&
\log\expect_{\nun_s}
\Big(
\exp\big\{
\gamma L \big|\un(\frac{j}{n})\big| \ind_{\{|\un(\frac{j}{n})|>M\}}
\big\}
\Big)
\le1.
\end{eqnarray}
(ii)
For any
$\gamma\in(0,1/(8C^2))$ there exists
$M<\infty$, such that  for any
$n$,
$j\in\Tn$ and
$s\in[0,T]$
\begin{eqnarray}
\notag
&&
\log\expect_{\nun_s}
\Big(
\exp\big\{
\gamma L \big|\un(\frac{j}{n})\big|^2 \ind_{\{\rn(\frac{j}{n})>M\}}
\big\}
\Big)
\le1,
\\[-5pt]
\label{eq:largedevi2}
\\[-5pt]
\notag
&&
\log\expect_{\nun_s}
\Big(
\exp\big\{
\gamma L \big|\un(\frac{j}{n})\big|^2 \ind_{\{|\un(\frac{j}{n})|>M\}}
\big\}
\Big)
\le1.
\end{eqnarray}
\end{lemma}

\begin{proof}
(i)
We prove the first bound of
(\ref{eq:largedevi1}),  the other ones are done very similarly.

Let
$Z_L$ be a
$POI(L)$-distributed random variable. Using the stochastic domination
(\ref{eq:stochboundrho})  we obtain
\begin{eqnarray*}
&&
\log\expect_{\nun_s}
\Big(
\exp\big\{
\gamma L \rn(\frac{j}{n}) \ind_{\{\rn(\frac{j}{n})>M\}}
\big\}
\Big)
\\[5pt]
&&
\hskip15mm
\le
\log
\Big(
1+
\expect_{\nun_s}
\Big(
\exp\big\{
\gamma L \rn(\frac{j}{n})\big\}
\ind_{\{\rn(\frac{j}{n})>M\}}
\Big)
\Big)
\\[5pt]
&&
\hskip15mm
\le
\sqrt{
\expect_{\nun_s}
\Big(
\exp\big\{
2\gamma L \rn(\frac{j}{n})\big\}
\Big)
}
\sqrt{
\prob_{\nun_s}
\Big(
\rn(\frac{j}{n})>M
\Big)
}
\\[5pt]
&&
\hskip15mm
\le
\sqrt{
\expect
\Big(
\exp\big\{
\gamma C Z_L
\big\}
\Big)
}
\sqrt{
\prob
\Big(
(C/M)Z_L>L
\Big)
}
\\[5pt]
&&
\hskip15mm
\le
\exp
\big\{
\frac{L}{2}
\big(
(e^{(\alpha C)/M}-1)
+
(e^{\gamma C}-1)
-\alpha
\big)
\big\},
\end{eqnarray*}
where $\alpha$ is arbitrary positive number. In the last step
Markov's inequality is being used. Now, choosing
$\alpha>\exp(2\gamma C)$ and $M>(C\alpha)/(\ln2)$ we obtain
(\ref{eq:largedevi1}).

(ii)
Again, we prove the first bound in
(\ref{eq:largedevi2}).  The other one is done in an identical way.

Let again
$Z_L$ be a
$POI(L)$-distributed and
$X$ be a standard Gaussian random variable. Using the stochastic
dominations
(\ref{eq:stochboundrho})  and
(\ref{eq:stochboundu}) we obtain
\begin{eqnarray*}
&&
\log\expect_{\nun_s}
\Big(
\exp\big\{
\gamma L \big|\un(\frac{j}{n})\big|^2 \ind_{\{\rn(\frac{j}{n})>M\}}
\big\}
\Big)
\\[5pt]
&&
\hskip15mm
\le
\log
\Big(
1+
\expect_{\nun_s}
\Big(
\exp\big\{
\gamma L \big|\un(\frac{j}{n})\big|^2\big\}
\ind_{\{\rn(\frac{j}{n})>M\}}
\Big)
\Big)
\\[5pt]
&&
\hskip15mm
\le
\sqrt{
\expect_{\nun_s}
\Big(
\exp\big\{
2\gamma L \big|\un(\frac{j}{n})\big|^2\big\}
\Big)
}
\sqrt{
\prob_{\nun_s}
\Big(
\rn(\frac{j}{n})>M
\Big)
}
\\[5pt]
&&
\hskip15mm
\le
\sqrt{
\expect
\Big(
\exp\big\{
4\gamma C^2 \big(X^2+L\big)
\big\}
\Big)
}
\sqrt{
\prob
\Big(
Z_L>(M/C)L
\Big)
}
\\[5pt]
&&
\hskip15mm
\le
\big(1-8\gamma C^2\big)^{-1/4}
\exp
\big\{
\frac{L}{2}
\big(
4\gamma C^2
+
(e^{(\alpha C)/M}-1)
-\alpha
\big)
\big\},
\end{eqnarray*}
where $\alpha$  is arbitrary positive number. Given
$\gamma<1/(8C^2)$, we  choose $\alpha$ sufficiently large and
$M>(C\alpha)/(\ln2)$  to  obtain (\ref{eq:largedevi2}).
\end{proof}

Now we turn to the proof of Proposition
\ref{propo:ldbounds}:

\begin{proof}
The bounds
(\ref{eq:ldboundrho+u}), respectively,
(\ref{eq:ldboundusq}) follow directly from the entropy inequality
(\ref{eq:fixedtime}) of Lemma
\ref{lemma:fixedtime} and the bounds
(\ref{eq:largedevi1}), respectively,
(\ref{eq:largedevi2})  of Lemma
\ref{lemma:largedevi}. Recall that
$L\gg1$, as
$n\to\infty$.
\end{proof}

\subsection{Proof of the fluctuation bounds (Proposition
  \ref{propo:kurschakbounds})}
\label{subs:kurschakboundsproof}

Within this proof we need the notation
\begin{eqnarray*}
\unt(s,x)
&:=&
\frac{n^\b}{l}
\sum_{k}
a\big(\frac{nx-k}{l}\big)
\left(\zeta_k - n^{-\b}u\big(s,\frac{k}{n}\big)\right)
=
\un(x) - \expect_{\nun_s}\big(\un(x)\big),
\\[5pt]
\rnt(s,x)
&:=&
\frac{n^{2\b}}{l}
\sum_{k}
a\big(\frac{nx-k}{l}\big)
\left(\eta_k - n^{-2\b}\rho\big(s,\frac{k}{n}\big)\right)
=
\rn(x) - \expect_{\nun_s}\big(\rn(x)\big).
\end{eqnarray*}
Since
\begin{eqnarray*}
&&
\hskip-6mm
\Big|\,
\big|\,
\un\big(s,\frac{j}{n}\big)-u\big(s,\frac{j}{n}\big)
\,\big|
-
\big|\,
\unt\big(s,\frac{j}{n}\big)
\,\big|
\,\Big|
\\[5pt]
&&
\hskip0.5mm
\le
\big|\,
1-\frac{1}{l}\sum_{k}a\big(\frac{j-k}{l}\big)
\,\big|
\,
\big|\,
u\big(s,\frac{j}{n}\big)
\,\big|
+
\frac{1}{l}\sum_{k}a\big(\frac{j-k}{l}\big)
\,
\big|\,
u\big(s,\frac{j}{n}\big)-u\big(s,\frac{k}{n}\big)
\,\big|
\\[5pt]
&&
\hskip0.5mm
\le
C\left(\frac{1}{l}+\frac{l}{n}\right)
=
o(1),
\end{eqnarray*}
and, similarly
\begin{eqnarray*}
&&
\hskip-6mm
\Big|\,
\big|\,
\rn\big(s,\frac{j}{n}\big)-\rho\big(s,\frac{j}{n}\big)
\,\big|
-
\big|\,
\rnt\big(s,\frac{j}{n}\big)
\,\big|
\,\Big|
\\[5pt]
&&
\hskip0.5mm
\le
\big|\,
1-\frac{1}{l}\sum_{k}a\big(\frac{j-k}{l}\big)
\,\big|
\,
\big|\,
\rho\big(s,\frac{j}{n}\big)
\,\big|
+
\frac{1}{l}\sum_{k}a\big(\frac{j-k}{l}\big)
\,
\big|\,
\rho\big(s,\frac{j}{n}\big)-\rho\big(s,\frac{k}{n}\big)
\,\big|
\\[5pt]
&&
\hskip0.5mm
\le
C\left(\frac{1}{l}+\frac{l}{n}\right)
=
o(1),
\end{eqnarray*}
we have to prove
\begin{eqnarray}
\label{eq:kurschakboundut}
\expect\Big(
\frac1n \sum_{j\in\Tn}
\big|\,
\unt (s,\frac{j}{n})
\,\big|^2
\Big)
\le
C\,\hn(s) + o(1),
\end{eqnarray}
respectively,
\begin{eqnarray}
\label{eq:kurschakboundrhot}
&&
\expect\Big(
\frac1n \sum_{j\in\Tn}
\big|\,
\rnt(s,\frac{j}{n})
\,\big|^2
\,
\ind_{ \{ |\,\rnt(s,\frac{j}{n})\,| \,\le  M \} }
\Big)
\le
C\,\hn(s) + o(1).
\end{eqnarray}

\begin{lemma}
\label{lemma:kurschak}
(i)
There exists
$\gamma>0$  (sufficiently small) such that for all
$n$,
$j\in\Tn$  and
$s\in[0,T]$
\begin{eqnarray}
\label{eq:kurschaku}
\log\expect_{\nun_s}
\Big(
\exp\big\{
\gamma L
\,
\big|\,
\unt\big(s,\frac{j}{n}\big)
\,\big|^2
\big\}
\Big)
\le 1
\end{eqnarray}
(ii)
For any
$M<\infty$  there exists
$\gamma>0$  (sufficiently small) such that for all
$n$,
$j\in\Tn$ and
$s\in[0,T]$
\begin{eqnarray}
\label{eq:kurschakrho}
\log\expect_{\nun_s}
\Big(
\exp\big\{
\gamma L
\,
\big|\,
\rnt\big(s,\frac{j}{n}\big)
\,\big|^2
\ind_{\{\,|\,\rnt(s,\frac{j}{n})\,|\,\le M\}}
\big\}
\Big)
\le 1.
\end{eqnarray}
\end{lemma}

\begin{proof}
(i)
Let
$X$ be a standard Gaussian random variable, which is independent of
all other random variables appearing in this paper,  and denote by
$\langle\dots\rangle$  expectation with respect to
$X$.
\begin{eqnarray}
\label{eq:kur}
&&
\log\expect_{\nun_s}
\Big(
\exp\big\{
\gamma L
\,
\big|\,
\unt\big(s,\frac{j}{n}\big)
\,\big|^2
\big\}
\Big)
\\[5pt]
\notag
&&
\hskip15mm
=
\log\expect_{\nun_s}
\Big(
\exp\big\{
\frac{\gamma}{l}
\,
\big|\,
\sum_k
a\big(\frac{j-k}{l}\big)
\big(\zeta_k-\expect_{\nun_s}(\zeta_k)\big)
\,\big|^2
\big\}
\Big)
\\[5pt]
\notag
&&
\hskip15mm
=
\log\Big\langle
\expect_{\nun_s}
\Big(
\exp\big\{
X
\sqrt{\frac{2\gamma}{l}}
\,
\sum_k
a\big(\frac{j-k}{l}\big)
\big(\zeta_k-\expect_{\nun_s}(\zeta_k)\big)
\big\}
\Big)
\Big\rangle.
\end{eqnarray}
Now, note that  the random variables
$\zeta_k-\expect_{\nun_s}(\zeta_k)$,
$k\in\Tn$, are uniformly bounded and under the distribution
$\prob_{\nun_s}$ they are independent and have zero mean. Hence there
exists a finite constant
$C$  such that for any collection of real numbers
$\lambda_k$,
$k\in\Tn$
\begin{eqnarray*}
\expect_{\nun_s}
\Big(
\exp\big\{
\sum_k
\lambda_k
\big(\zeta_k-\expect_{\nun_s}(\zeta_k)\big)
\big\}
\Big)
\le
\exp\big\{C\sum_k \lambda_k^2\big\}.
\end{eqnarray*}
Further on,
there exists a finite constant $C$ such that for any  $l$
\begin{eqnarray}
\label{eq:sumsq}
\frac{1}{l}\sum_k \big|\,a\big(\frac{k}{l}\big)\,\big|^2\le C.
\end{eqnarray}
From these it follows that for some finite constant
$C$,
\begin{eqnarray*}
\text{r.h.s. of (\ref{eq:kur})}\,\,
\le
\log\left\langle
\exp\big\{
C\,\gamma\,X^2
\big\}
\right\rangle.
\end{eqnarray*}
Choosing
$\gamma$  sufficiently small in this last inequality we obtain
(\ref{eq:kurschaku}).

(ii)
Note first that, given
$M<\infty$  fixed,  there exists a zero mean bounded random variable
$Y$ such that for any
$r\in\R$
\begin{eqnarray*}
r^2\ind_{\{\,|\,r\,|\, \le \,M\,\}}
\le
\log
\expect\Big(\exp\big\{\,r\,Y\,\big\}\Big).
\end{eqnarray*}
Let
$Y_1,Y_2,\dots$  be i.i.d. copies of
$Y$ which are also independent of all other random variables appearing
in this paper, and denote by
$\langle\dots\rangle$ expectation with respect to these. Then we have
\begin{eqnarray}
\label{eq:kurs}
&&
\log\expect_{\nun_s}
\Big(
\exp\big\{
\gamma L
\,
\big|\,
\rnt\big(s,\frac{j}{n}\big)
\,\big|^2
\ind_{\{\,|\,\rnt(s,\frac{j}{n})\,|\,\le M\}}
\big\}
\Big)
\\[5pt]
\notag
&&
\hskip15mm
\le
\log
\Big\langle
\expect_{\nun_s}
\Big(
\exp\big\{
\frac{\sum_{p=1}^{\lceil \gamma L \rceil}Y_p}{L}
\sum_k
a\big(\frac{j-k}{l}\big)
\big(\eta_k-\expect_{\nun_s}(\eta_k)\big)
\big\}
\Big)
\Big\rangle.
\end{eqnarray}
Next note that for any
$\overline\lambda<\infty$ there exists a constant
$C<\infty$ such that for any
$n\in\N$,  any
$s\in[0,T]$  and any collection of real numbers
$\lambda_k\in[-\overline\lambda,\overline\lambda]$,
$k\in\Tn$
\begin{eqnarray*}
\expect_{\nun_s}
\Big(
\exp\big\{
\sum_k
\lambda_k\,
\big(\eta_k-\expect_{\nun_s}(\eta_k)\big)
\big\}
\Big)
\le
\exp\big\{
C
n^{-2\beta}
\sum_k
\lambda_k^2
\big\}.
\end{eqnarray*}
Hence, using again
(\ref{eq:sumsq}),
\begin{eqnarray*}
&&
\text{r.h.s.  of (\ref{eq:kurs})}\,\,
\le
\log
\Big\langle
\exp\big\{C\gamma\,
\left(
\big(Y_1+\dots_+Y_{\lceil \gamma L \rceil}\big)/\sqrt{\lceil \gamma L \rceil}
\right)^2
\big\}
\Big\rangle.
\end{eqnarray*}
Now, since the i.i.d. random variables
$Y_1,Y_2,\dots$  are bounded and have zero mean, choosing
$\gamma$  sufficiently small this last  expression can be made
arbitrarily small, uniformly in
$L$. Hence
(\ref{eq:kurschakrho}).
\end{proof}

Now back to the proof of Proposition
\ref{propo:kurschakbounds}.

\begin{proof}
From
(\ref{eq:fixedtime}) and
(\ref{eq:kurschaku}), respectively, from
(\ref{eq:fixedtime}) and
(\ref{eq:kurschakrho}) we deduce
(\ref{eq:kurschakboundut}), respectively,
(\ref{eq:kurschakboundrhot}). Finally, these two bounds and the
arguments at the beginning of the present subsection imply
(\ref{eq:kurschakboundu}), respectively,
(\ref{eq:kurschakboundrho}).
\end{proof}

\subsection{Proof of the block replacement and gradient bounds
  (Proposition   \ref{propo:aprioribounds})}
\label{subs:aprioriboundsproof}

\subsubsection{An elementary probability lemma}
\label{subsubs:epl}

Let
$(\Omega,\pi)$ be a finite probability space and
$\omega_i$,
$i\in\Z$ i.i.d.
$\Omega$-valued random variables with distribution
$\pi$. Further on let
\begin{eqnarray*}
\begin{array}{ll}
\vzeta:\Omega\to\R^d,
\quad
&
\vzeta_i:=\vzeta(\omega_i),
\\[5pt]
\xi:\Omega^{m}\to\R,
\quad
&
\xi_i:=\xi(\omega_{i}\dots,\omega_{i+m-1}).
\end{array}
\end{eqnarray*}
For
$\vx\in\text{co}(\text{Ran}(\vzeta))$ denote
\[
\Xi(\vx):=
\frac
{\expect_\pi\big(\xi_1\exp\{\sum_{i=1}^{m}\vlam\cdot\vzeta_i\}\}\big)}
{\expect_\pi\big(\exp\{\vlam\cdot\vzeta_1\}\}\big)^{m}},
\]
where
$\text{co}(\cdot)$  stands for `convex hull' and
$\vlam\in\R^d$ is chosen so that
\[
\frac
{\expect_\pi\big(\vzeta_1\exp\{\vlam\cdot\vzeta_1\}\}\big)}
{\expect_\pi\big(\exp\{\vlam\cdot\vzeta_1\}\}\big)}
=
\vx.
\]
For
$l\in\N$ we denote
\emph{plain} block averages by
\[
\overline{\vzeta}_l:=\frac1l\sum_{j=1}^l\vzeta_j.
\]
Finally, let
$b:[0,1]\to\R$ be a fixed smooth function and denote
\[
M(b):=\int_0^1 b(s)\,ds,
\qquad
V(b):=M(b^2)-M(b)^2.
\]
We also define  the block averages
\emph{weighted by $b$}
\[
\langle b\,,\,\vzeta\rangle_l
:=
\frac1l\sum_{j=0}^l
b(j/l)\vzeta_j,
\quad
\langle b\,,\,\xi\rangle_l
:=
\frac1l\sum_{j=0}^l
b(j/l)\xi_j,
\]

The following lemma  relies on elementary probability arguments:

\begin{lemma}
\label{lemma:epl1}
{\rm (Microcanonical exponential moments of block averages)}
\\
There exists a constant
$C<\infty$, depending only on
$m$, on the joint distribution of
$(\xi_i,\vzeta_i)$ and on the function
$b$, such that the following bounds hold uniformly in
$l\in\N$ and
$\vx\in({\rm{Ran}}(\vzeta)+\dots+{\rm{Ran}}(\vzeta))/l$:
\\
(i)
If
$M(b)=0$, then
\begin{eqnarray}
\label{eq:epl11}
\expect
\Big(
\exp
\big\{
\gamma \sqrt{l}
\langle b\,,\,\xi\rangle_l
\big\}
\,\Big|\,
\overline{\vzeta}_l=\vx
\Big)
\le
\exp\{C(\gamma^2+\gamma/\sqrt l)\}.
\end{eqnarray}
(ii)
If
$M(b)=1$ then
\begin{eqnarray}
\label{eq:epl12}
\expect
\Big(
\exp
\big\{
\gamma\sqrt{l}
\big(
{\langle b\,,\,\xi\rangle}_l
-
\Xi(\langle b\,,\,\vzeta\rangle_l)
\big)
\big\}
\,\Big|\,
\overline{\vzeta}_l=\vx
\Big)
\le
\exp\{C(\gamma^2+\gamma/\sqrt l)\}.
\end{eqnarray}
\end{lemma}

\begin{proof}
We prove the lemma with
$m=1$,  that is with
$\left(\xi_i\right)_{i=1}^l$  independent rather than
$m$-dependent. The
$m$-dependent case follows by applying Jensen's inequality in a rather
straightforward way.
\\
(i)
In order to simplify the argument we make the assumption that the
function
$s\mapsto b(s)$  is odd:
\begin{eqnarray}
\label{eq:bodd}
b(1-s)=-b(s).
\end{eqnarray}
The same argument works if the function
$s\mapsto b(s)$  can be rearranged (by  permutation of finitely many
subintervals of
$[0,1]$)  into a piecewise continuous odd function. This case is
sufficient for our purposes. The proof of the fully general case
--- which goes through induction on
$l$  --- is more tedious and it is left as a fun exercise for the
reader.

Assuming
(\ref{eq:bodd}) we have
\[
\sqrt{l}
\langle b\,,\,\xi\rangle_l
=
l^{-1/2}
\sum_{j=0}^{[l/2]} b(j/l) (\xi_j-\xi_{l-j})
\]
and hence
\begin{eqnarray*}
&&
\hskip-15mm
\expect
\Big(
\exp
\big\{
\gamma \sqrt{l}
\langle b\,,\,\xi\rangle_l
\big\}
\,\Big|\,
\overline{\vzeta}_l=\vx
\Big)
\\[5pt]
&&
\hskip1cm
=
\expect
\Big(
\expect
\Big(
\exp
\big\{
\gamma \sqrt{l}
\langle b\,,\,\xi\rangle_l
\big\}
\,\Big|\,
\vzeta_j+\vzeta_{l-j}:\,\,j=0,\dots, l
\Big)
\,\Big|\,
\overline{\vzeta}_l=\vx
\Big)
\\[5pt]
&&
\hskip1cm
=
\expect
\Big(
\prod_{j=0}^{[l/2]}
\expect
\Big(
\exp
\big\{
\gamma l^{-1/2} b(j/l) (\xi_j-\xi_{l-j})
\big\}
\,\Big|\,
\vzeta_j+\vzeta_{l-j}
\Big)
\,\Big|\,
\overline{\vzeta}_l=\vx
\Big)
\\[5pt]
&&
\hskip1cm
\le
\exp
\big\{
C \gamma^2
\sum_{j=1}^{[l/2]}
l^{-1}  b(j/l)^2
\big\}
\\[5pt]
&&
\hskip1cm
=
\exp
\big\{
C \gamma^2 (V(b) + \Ordo(1/l))
\big\}.
\end{eqnarray*}
In the second step we use the fact that the pairs
$\big(\xi_j,\xi_{l-j}\big)$,
$j=0,\dots,[l/2]$ are independent, given
$\vzeta_j+\vzeta_{l-j}$,
$j=0,\dots,[l/2]$. In the third  step we note that the variables
$\xi_j$ are bounded and
$\expect\big(\xi_j-\xi_{l-j}\big|\vzeta_j+\vzeta_{l-j}\big)=0$.
\\
(ii)
Beside
$\Xi(\vx)$ we also introduce the functions
\[
\Xi_l:\big({\rm Ran}(\vzeta)+\cdots+{\rm Ran}(\vzeta)\big)/l\to\R,
\qquad
\Xi_l(\vx):=
\expect\big(\xi_1\big|\overline{\vzeta}_l=\vx\big).
\]
We shall exploit the following  facts
\\
(1)
The functions
$\Xi(\vx)$  and
$\Xi_l(\vx)$  are uniformly bounded. This follows from the
boundedness of
$\xi_j$.
\\
(2)
The function
$\vx\mapsto\Xi(\vx)$  is smooth with bounded first two
derivatives. This follows from direct computations.
\\
(3)
There exists a finite constant
$C$, such that
\[
\abs{\Xi_l(\vx)-\Xi(\vx)}
\le C l^{-1}.
\]
This follows from the so-called equivalence of ensembles (see
e.g.~Appendix 2 of \cite{kipnislandim}).

We write
\begin{eqnarray}
\label{eq:dekomp}
{\langle b\,,\,\xi\rangle}_l
-
\Xi(\langle b\,,\,\vzeta\rangle_l)
&=&
\phantom{+}
\left(
{\langle b\,,\,\xi\rangle}_l
-
\overline{\xi}_l
\right)
+
\left(
\overline{\xi}_l
-
\Xi_l(\overline{\vzeta}_l)
\right)
\\[5pt]
\notag
&&
+
\left(
\Xi_l(\overline{\vzeta}_l)
-
\Xi(\overline{\vzeta}_l)
\right)
+
\left(
\Xi(\overline{\vzeta}_l)
-
\Xi(\langle b\,,\,\vzeta\rangle_l)
\right).
\end{eqnarray}
By applying Jensen's inequality we conclude that we have to bound
the exponential moments of type
(\ref{eq:epl12}), separately for the four terms.

Bounding the first and last terms reduces directly to
(\ref{eq:epl11}),  the third term is uniformly
$\Ordo(l^{-1})$, so we only have to bound the exponential moments of
the second term in
(\ref{eq:dekomp}). This is done by induction on
$l$. Let
$C(l)$ be the best constant such that for any
$\gamma\in\R$
\[
\expect
\Big(
\exp
\big\{
\gamma\sqrt{l}
\big(
\overline{\xi}_l
-
\Xi_l(\overline{\vzeta}_l)
\big)
\big\}
\,\Big|\,
\overline{\vzeta}_l=\vx
\Big)
\le
\exp\{C(l) \gamma^2\}.
\]
We prove that
$C(l)$ stays bounded as
$l\to\infty$.

The following identity holds
\begin{eqnarray*}
\sqrt{l+1}
\big(
\overline{\xi}_{l+1}
-
\Xi_{l+1}(\overline{\vzeta}_{l+1})
\big)
&=&
\phantom{+}
\frac{l}{\sqrt{l+1}}
\big(
\overline{\xi}_l
-
\Xi_l(\overline{\vzeta}_l)
\big)
\\[5pt]
&&
+
\frac{1}{\sqrt{l+1}}
\big(
\xi_{l+1}
-
\Xi_1(\vzeta_{l+1})
\big)
\\[5pt]
&&
+
\frac{l}{\sqrt{l+1}}
\big(
\Xi_l(\overline{\vzeta}_l)-\Xi_{l+1}(\overline{\vzeta}_{l+1})
\big)
\\[5pt]
&&
+
\frac{1}{\sqrt{l+1}}
\big(
\Xi_1({\vzeta}_{l+1})-\Xi_{l+1}(\overline{\vzeta}_{l+1})
\big)
\end{eqnarray*}
Thus
\begin{eqnarray*}
\expect
\Big(
\exp
\big\{
\gamma\sqrt{l+1}
\big(
\overline{\xi}_{l+1}
-
\Xi_{l+1}(\overline{\vzeta}_{l+1})
\big)
\big\}
\,\Big|\,
\overline{\vzeta}_{l+1}=\vx
\Big)
=
\hskip3cm
\\[5pt]
\expect
\Big(
\expect
\Big(
\exp
\big\{
\gamma\sqrt{l+1}
\big(
\overline{\xi}_{l+1}
-
\Xi_{l+1}(\overline{\vzeta}_{l+1})
\big)
\big\}
\,\Big|\,
\overline{\vzeta}_{l}\,,\,\vzeta_{l+1}
\Big)
\,\Big|\,
\overline{\vzeta}_{l+1}=\vx
\Big)
=
\\[5pt]
\expect
\Bigg(
\hskip2cm
\expect
\Big(
\exp
\big\{
\frac{ \gamma l }{\sqrt{l+1}}
\big(
\overline{\xi}_{l}
-
\Xi_{l}(\overline{\vzeta}_{l})
\big)
\big\}
\,\Big|\,
\overline{\vzeta}_{l}
\Big)
\times
\phantom{
\,\Big|\,
\overline{\vzeta}_l=\vx
\Big)M
}
\\[5pt]
\expect
\Big(
\exp
\big\{
\frac{\gamma}{\sqrt{l+1}}
\big(
\overline{\xi}_{l+1}
-
\Xi_{1}(\overline{\vzeta}_{l+1})
\big)
\big\}
\,\Big|\,
{\vzeta}_{l+1}
\Big)
\times
\phantom{
\,\Big|\,
\overline{\vzeta}_l=\vx
\Big)M
}
\\[5pt]
\exp
\big\{
\frac{\gamma l}{\sqrt{l+1}}
\big(
\Xi_l(\overline{\vzeta}_l)-\Xi_{l+1}(\overline{\vzeta}_{l+1})
\big)
\big\}
\times
\phantom{
\,\Big|\,
\overline{\vzeta}_l=\vx
\Big)M
}
\\[5pt]
\exp
\big\{
\frac{\gamma}{\sqrt{l+1}}
\big(
\Xi_1({\vzeta}_{l+1})-\Xi_{l+1}(\overline{\vzeta}_{l+1})
\big)
\big\}
\,\Bigg|\,
\overline{\vzeta}_l=\vx
\Bigg).
\end{eqnarray*}
The terms
\[
\big(
\overline{\xi}_{l+1} - \Xi_{1}(\overline{\vzeta}_{l+1})
\big)
,
\quad\,\,\,
l
\big(
\Xi_l(\overline{\vzeta}_l)- \Xi_{l+1}(\overline{\vzeta}_{l+1})
\big),
\quad\,\,\,
\big(
\Xi_1({\vzeta}_{l+1})- \Xi_{l+1}(\overline{\vzeta}_{l+1})
\big)
\]
are uniformly bounded and
\begin{eqnarray*}
\expect
\big(
\overline{\xi}_{l+1}
-
\Xi_{1}(\overline{\vzeta}_{l+1})
\,\big|\,
{\vzeta}_{l+1}
\big)
&=&
0,
\\[5pt]
\expect
\big(
\Xi_l(\overline{\vzeta}_l)-\Xi_{l+1}(\overline{\vzeta}_{l+1})
\,\big|\,
\overline{\vzeta}_{l+1}
\big)
&=&
0,
\\[5pt]
\expect
\big(
\Xi_1({\vzeta}_{l+1})-\Xi_{l+1}(\overline{\vzeta}_{l+1})
\,\big|\,
\overline{\vzeta}_{l+1}
\big)
&=&
0.
\end{eqnarray*}
Using the induction hypothesis it
follows that there exists a finite constant
$B$  such that
\[
C(l+1)\le \frac{l}{l+1}C(l) + \frac{1}{l+1} B.
\]
Hence,
$\limsup_{l\to\infty}C(l)\le B$ and the lemma follows.
\end{proof}

\begin{lemma}
\label{lemma:epl}
{\rm (Microcanonical Gaussian bounds)}
\\
There exists a
$\gamma_0>0$, depending only on $m$, on the joint distribution of
$(\xi_i,\vzeta_i)$ and on the function
$b$, such that the following bounds hold uniformly in
$l\in\N$ and
$\vx\in({\rm{Ran}}(\vzeta)+\dots+{\rm{Ran}}(\vzeta))/l$:
\\
(i)
If
$M(b)=0$, then
\begin{eqnarray}
\label{eq:epl1}
\log \expect
\Big(
\exp
\big\{
\gamma_0 l
{\langle b\,,\,\xi\rangle}_l^2
\big\}
\,\Big|\,
\overline{\vzeta}_l=\vx
\Big)
\le
1.
\end{eqnarray}
(ii)
If
$M(b)=1$ then
\begin{eqnarray}
\label{eq:epl2}
\log \expect
\Big(
\exp
\big\{
\gamma_0 l
\big(
{\langle b\,,\,\xi\rangle}_l
-
\Xi({\langle b\,,\,\vzeta\rangle}_l)
\big)^2
\big\}
\,\Big|\,
\overline{\vzeta}_l=\vx
\Big)
\le
1.
\end{eqnarray}
\end{lemma}

\begin{proof}
This is actually a corollary of Lemma \ref{lemma:epl1}: The bounds
(\ref{eq:epl1}) and (\ref{eq:epl2}) follow from (\ref{eq:epl11}),
respectively, (\ref{eq:epl12}) by exponential Gaussian averaging
(as in the proof of Lemma \ref{lemma:kurschak}).
\end{proof}

\subsubsection{Proof of Proposition   \ref{propo:aprioribounds}}
\label{subsubs:aprioriboundsproof}

Now we turn to the proof of Proposition
\ref{propo:aprioribounds}.

\begin{proof}
(i)
In order to prove
(\ref{eq:plainobrepl}) first note that by simple numerical
approximation (no probability bounds involved)
\begin{eqnarray*}
\abs{
\int_{\T}
\big|
\big\{
\wih{\xi}
-
\Xi(\wih{\eta},\wih{\zeta})
\big\}
(x)
\big|^2
\,dx
-
\frac1n\sum_{j\in\Tn}
\big|
\big\{
\wih{\xi}
-
\Xi(\wih{\eta},\wih{\zeta})
\big\}
(\frac{j}{n})
\big|^2
}
\phantom{MMMM}
&&
\\[5pt]
\le
\frac{C}{l}
=
o\left(\frac{l^2}{n^{1+3\b+\d}}\right).
&&
\end{eqnarray*}
We apply Lemma
\ref{lemma:varadhan} with
\begin{eqnarray*}
{\cal V}
=
\big|\,
\big\{
\wih{\xi}
-
\Xi(\wih{\eta},\wih{\zeta})
\big\}
(0)
\,\big|^2
=
\big|\,
{\langle a\,,\,\xi \rangle}_l
-
\Xi\big(
{\langle a\,,\,\eta \rangle}_l,
{\langle a\,,\,\zeta \rangle}_l
\big)
\,\big|^2
\end{eqnarray*}
We use  the bound
(\ref{eq:epl2}) of Lemma
\ref{lemma:epl} with the function
$b=a$. Note that
$\gamma=\gamma_0l$ can be chosen in
(\ref{eq:varadhan}). This yields the bound
(\ref{eq:plainobrepl}).
\\
(ii)
In order to prove
(\ref{eq:plaingrad}) we start again with numerical approximation:
\begin{eqnarray*}
\abs{
\int_{\T}
\big|
\px\wih{\xi}(x)
\big|^2
\,dx
-
\frac1n\sum_{j\in\Tn}
\big|
\px\wih{\xi}(\frac{j}{n})
\big|^2
}
\le
C\frac{n^2}{l^3}
=
o(n^{1-3\b-\d}).
\end{eqnarray*}
We apply Lemma
\ref{lemma:varadhan} with
\begin{eqnarray*}
{\cal V}
=
\big|
\px\wih{\xi}(0)
\big|^2
=
\frac{n^2}{l^2}
\,
\big|\,
{\langle a'\,,\,\xi \rangle}_l
\,\big|^2.
\end{eqnarray*}
We use now the bound
(\ref{eq:epl1}) of Lemma
\ref{lemma:epl} with the function
$b=a'$. Now we can choose
$\gamma=\gamma_0l^3/n^2$ and this will yield the bound
(\ref{eq:plaingrad}).
\\
(iii)
Next we prove
(\ref{eq:turbograd}). We apply Lemma
\ref{lemma:varadhan}  with
\begin{eqnarray*}
{\cal V}
&=&
\frac{|\,\px\wih{\xi}(0)\,|^2}{\wih{\eta}(0)}
=
\frac{n^{2}}{l^3}\,
\frac{\big|\, \sum_k a'(k/l) \xi_k \,\big|^2 }
     {\sum_k a(k/l) \eta_k }
\\[5pt]
&=&
\frac{n^{2}}{2l^3}\,
\frac{\big|\, \sum_k a'(k/l) (\xi_k-\xi_{-k}) \,\big|^2 }
     {\sum_k a(k/l) (\eta_k+\eta_{-k}) },
\end{eqnarray*}
where in the last equality we use the fact that the weighting function
$x\mapsto a(x)$  is
\emph{even}.  We compute the exponential moment
$\expect^{2l+1}_{N,Z}\big(\exp\{\gamma {\cal V}\}\big)$.  Let
$X$
be a standard Gaussian random variable, which is independent of all
other random variables appearing in this paper and denote by
$\langle\dots\rangle$  averaging with respect to it. We have
\begin{eqnarray*}
&&
\expect^{2l+1}_{N,Z}
\Big(
\exp\big\{\gamma {\cal V}\big\}
\Big)
\\[5pt]
&&
=
\expect^{2l+1}_{N,Z}
\Big(
\exp\big\{\gamma
\frac{n^{2}}{2l^3}\,
\frac{\big|\, \sum_k a'(k/l) (\xi_k-\xi_{-k}) \,\big|^2 }
     {\sum_k a(k/l) (\eta_k+\eta_{-k}) }
\big\}
\Big)
\\[5pt]
&&
=
\Big\langle
\expect^{2l+1}_{N,Z}
\Big(
\exp\big\{
X
\sqrt{\gamma}
\frac{n}{l^{3/2}}\,
\frac{ \sum_k a'(k/l) (\xi_k-\xi_{-k})  }
     {\sqrt{\sum_k a(k/l) (\eta_k+\eta_{-k})} }
\big\}
\Big)
\Big\rangle
\\[5pt]
&&
=
\Big\langle
\expect^{2l+1}_{N,Z}
\Big(
\expect^{2l+1}_{N,Z}
\Big(
\exp\big\{
X
\sqrt{\gamma}
\frac{n}{l^{3/2}}\,
\frac{ \sum_k a'(k/l) (\xi_k-\xi_{-k})  }
     {\sqrt{\sum_k a(k/l) (\eta_k+\eta_{-k})} }
\big\}
\,\Big|\,
{\big\{\eta_k+\eta_{-k}\big\}}_{k=0}^l
\Big)
\Big)
\Big\rangle
\\[5pt]
&&
\le
\Big\langle
\expect^{2l+1}_{N,Z}
\Big(
\exp\big\{
C
X^2
\gamma
\frac{n^{2}}{l^{3}}\,
\frac{ \sum_k a'(k/l)^2 (\eta_k+\eta_{-k})  }
     {\sum_k a(k/l) (\eta_k+\eta_{-k}) }
\big\}
\Big)
\Big\rangle
\\[5pt]
&&
\le
\Big\langle
\exp\big\{
C
X^2
\gamma
\frac{n^{2}}{l^{3}}
\big\}
\Big\rangle,
\end{eqnarray*}
where we used the facts that the random variables $\eta_k$
are non-negative, $\Om$ is finite and $\eta(\om)=0$ implies 
$\xi(\om)=0$. In the last step we used the inequality
\[
a'(x)^2\le C a(x),
\]
which follows from the conditions on $a(x)$, see subsection
\ref{subs:blocks}.
 
From this bound it follows that in Lemma
\ref{lemma:varadhan}  we can choose $\gamma=\gamma_0 l^3/n^{2}$,
with  a small but fixed $\gamma_0$, and hence the second bound in
(\ref{eq:turbograd}) follows.

\end{proof}

\section{Appendix: Some details about the PDE (\ref{eq:pde})}
\label{section:pdedetails}

\noindent
\underline{\emph{Hyperbolicity:}}
One has to analyze Jacobian the matrix
\beqs
D:=
\left(
\begin{array}{cc}
\big(\rho u\big)_\rho & \big(\rho u\big)_u
\\[8pt]
\big(\rho + \gamma u^2)_\rho &\big(\rho + \gamma u^2)_u
\end{array}
\right)
=
\left(
\begin{array}{cc}
u  & \rho
\\[8pt]
1  & 2\gamma u
\end{array}
\right).
\eeqs
The eigenvalues with the corresponding right and left  eigenvectors
are:
\beqs
Dr=\lambda r,
\quad
Ds=\mu s,
\quad
l^{\dagger}D=\lambda l^{\dagger},
\quad
m^{\dagger}D=\mu m^{\dagger},
\eeqs
($v^\dagger$ stands for  the transpose of the column 2-vector $v$).
The eigenvalues and eigenvectors are
\beqs
\left. 
\begin{array}{c}
\lambda\\\mu
\end{array}
\right\}
=
\pm
\frac12
\left\{
\sqrt{(2\gamma-1)^2u^2 + 4\rho} \pm (2\gamma+1)u,
\right\}
\eeqs
and
\beqs
&&
\left. 
\begin{array}{c}
r^\dagger\\s^\dagger
\end{array}
\right\}
=
\left(
\frac12
\left\{
\mp\sqrt{(2\gamma-1)^2u^2 + 4\rho} - (2\gamma-1)u
\right\},
1
\right),
\\[8pt]
&&
\left. 
\begin{array}{c}
l^\dagger\\m^\dagger
\end{array}
\right\}
=
\left( 1, - \frac12 \left\{ \pm\sqrt{(2\gamma-1)^2u^2 +
4\rho} - (2\gamma-1)u \right\} \right).
\eeqs
So, we can conclude
that the pde (\ref{eq:pde}) is (strictly) hyperbolic in the domain
\beqs
\gamma\not=1/2: && \quad \left\{(\rho,u)\in{\R_+\times\R}:
  (\rho,u)\not=(0,0)\right\},
\\
\gamma=1/2:&&
\quad
\left\{(\rho,u)\in{\R_+\times\R}:
  \rho\not=0\right\}.
\eeqs

\medskip
\noindent \underline{\emph{Riemann invariants:}} The Riemann
invariants $w=w(\rho,u)$, $z=z(\rho,u)$ of the pde are given by
the relations
\beqs
(w_\rho,w_u)\cdot s= 0 = (z_\rho,z_u)\cdot r.
\eeqs
That is, the level lines $w=\text{const.}$, respectively
$z=\text{const.}$ are determined by the ordinary differential
equations
\beqs
\left. 
\begin{array}{c}
w=\text{const.}\\z=\text{const.}
\end{array}
\right\}
:
\quad
\frac{d\rho}{du}=
\mp
\frac12
\left\{
\sqrt{(2\gamma-1)^2u^2 + 4\rho} \pm (2\gamma-1)u
\right\}
\eeqs
Actually  only the level lines of the functions
$w(\rho,u)$, respectively, $z(\rho,u)$ are determined.
In our case the Riemann invariants can be
found explicitly. For $\gamma\not=3/4$ we get
\beqs
\left. 
\begin{array}{c}
w(\rho,u)\\z(\rho,u)
\end{array}
\right\}
=
F\Big\{
\!\!
\left(\!
\sqrt{ (2\gamma-1)^2u^2 + 4\rho } \pm (2\gamma-1)u
\right)^{\frac{2\gamma-1}{2\gamma-2}}
\!\!
\left(\!
\sqrt{ (2\gamma-1)^2u^2 + 4\rho } \mp (2\gamma-2)u
\right)
\!\!
\Big\}
\eeqs
Where $F:\R\to\R$ is an appropriately chosen bijection (mind, that
only the level sets of the Riemann invariants are determined).

\begin{figure}[htp]
\centering \subfigure[$3/4<\g<1$]{
\includegraphics[width=(\textwidth-10mm)/2]{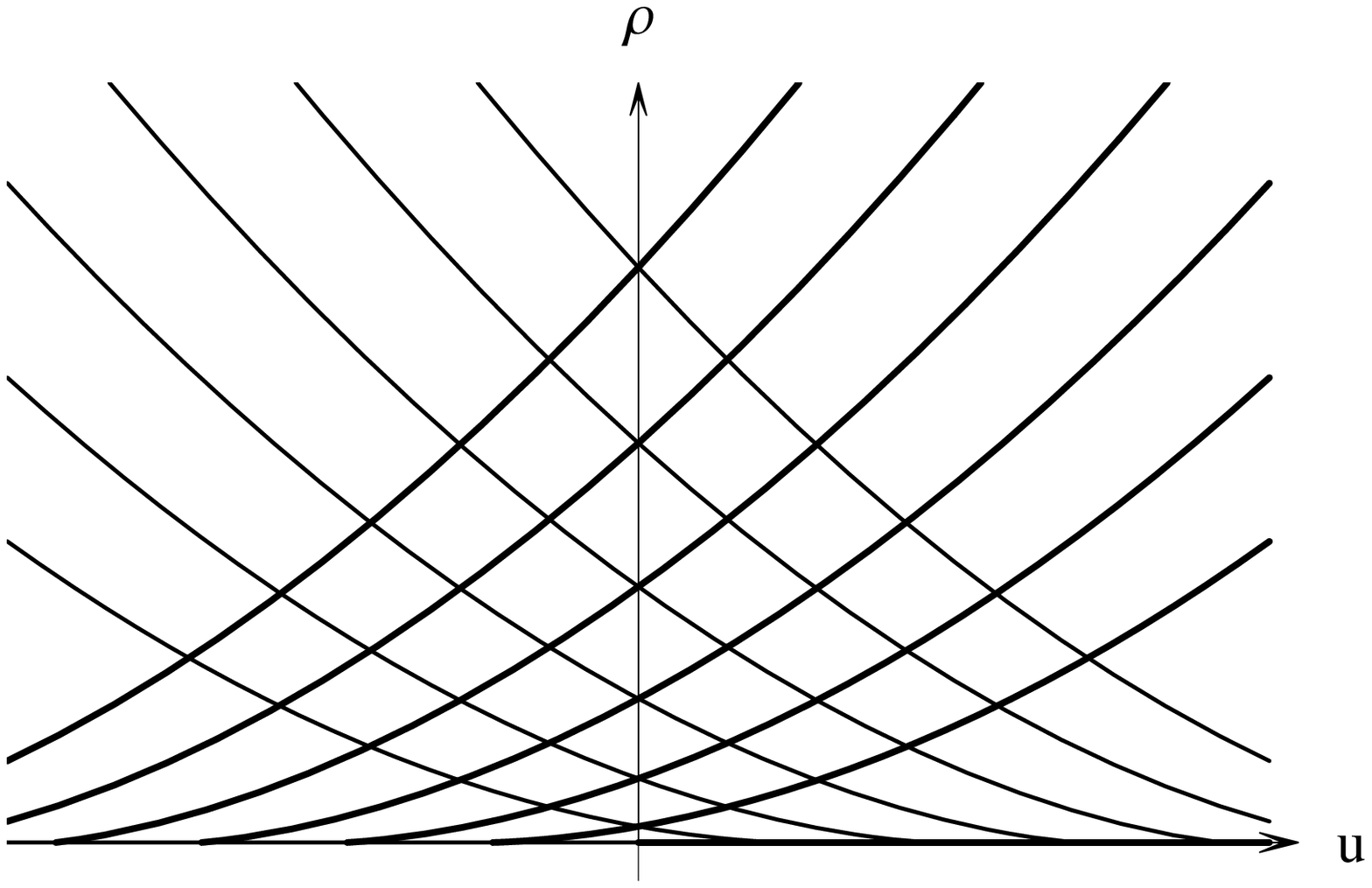}
} \subfigure[ $1<\g$]{
\includegraphics[width=(\textwidth-10mm)/2]{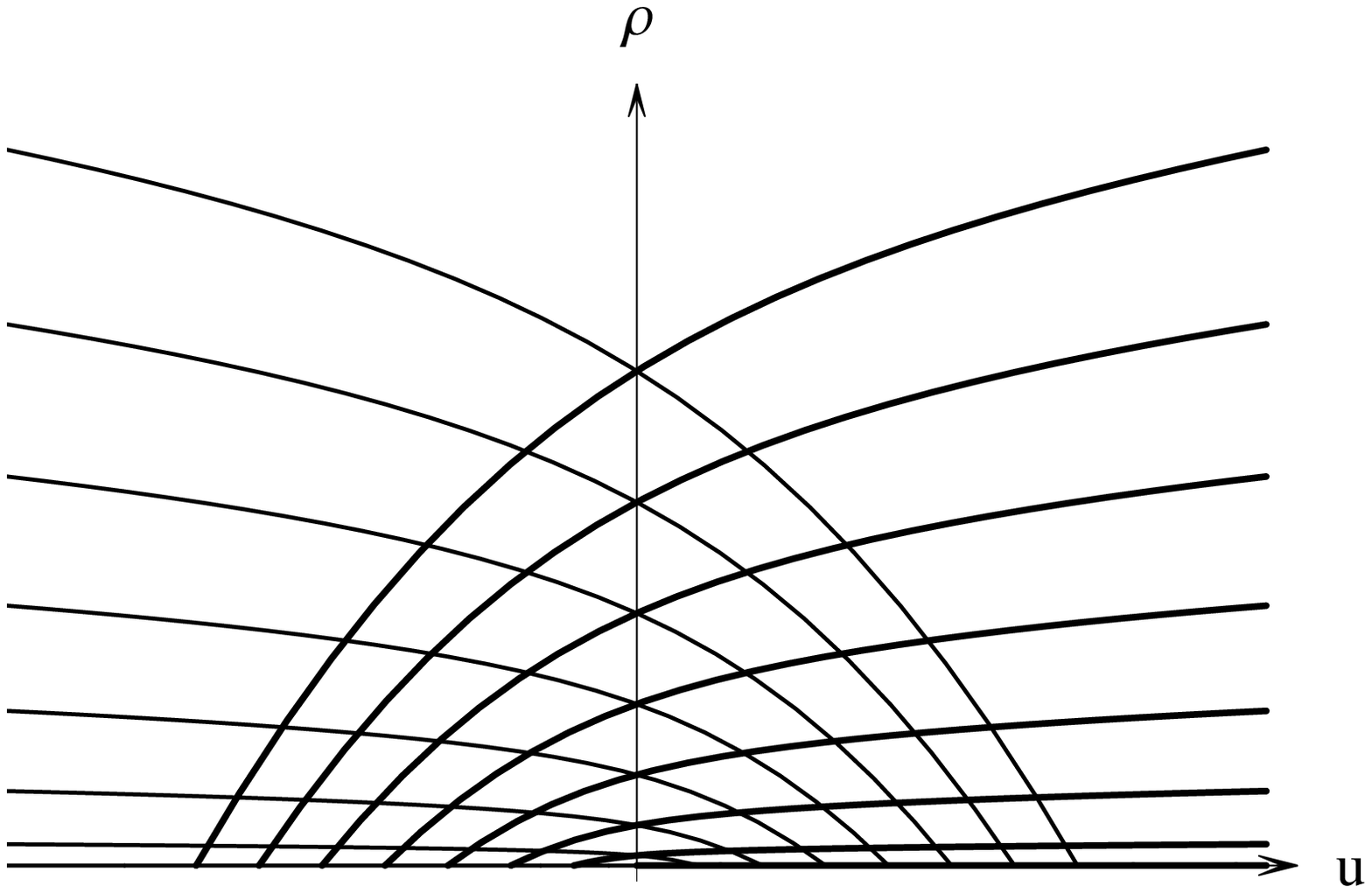}
}\caption{Level lines of Riemann-invariants} \label{fig:level}
\end{figure}

Note that due to the changes of sign of $2\gamma-1$ and $2\gamma
-2$, the above expression gives rise to \emph{qualitatively
different} behavior of the Riemann invariants. The picture changes
qualitatively at the critical values $\gamma=1/2$, $\gamma=3/4$
and $\gamma=1$. In Figure \ref{fig:level} we present the
qualitative picture of the level lines of $w(\rho,u)$  and
$z(\rho,u)$  for $3/4<\gamma<1$,
and $\gamma>1$,
respectively. (For economy reasons we omit the graphical
representation of the other 
cases, but encourage the reader
to sketch it.)
In all cases the Riemann invariants satisfy the convexity
conditions
\beq
\label{eq:cri}
\begin{array}{l}
w_{\rho\rho}w_u^2 - 2w_{\rho u}w_\rho w_u + w_{uu}w_\rho^2\ge0,
\\[5pt]
z_{\rho\rho}z_u^2 - 2z_{\rho u}z_\rho z_u + z_{uu}z_\rho^2\ge0,
\end{array}
\eeq
in ${\R_+\times\R}$ for all $\gamma$. (The sign of the function
$F(\cdot)$  is so chosen, that these expressions be non-negative.)
The inequalities are strict in  the interior of ${\R_+\times\R}$,
except
for the $\gamma=1$ case, when these expressions identically
vanish. These conditions are equivalent to saying that  the level
sets $\{(\rho,u)\in[0,\infty)\times(-\infty,\infty):
w(\rho,u)<c\}$ and $\{(\rho,u)\in[0,\infty)\times(-\infty,\infty):
z(\rho,u)<c\}$ be convex.  See \cite{lax1}, \cite{lax2} or
\cite{serre} for the importance of these convexity conditions.

It is of crucial importance for our  problem that the level
curves $w(\rho,u)=w=\text{const.}$ expressed as $u\mapsto\rho(u,w)$
are convex for $\gamma<1$, linear for $\gamma=1$ and concave for
$\gamma>1$.

\medskip
\noindent
\underline{\emph{Genuine nonlinearity:}}
Genuine nonlinearity holds if and only if
\beqs%
(\lambda_\rho,\lambda_u)\cdot r
\not=0\not=
(\mu_\rho,\mu_u)\cdot s .
\eeqs
in the interior of the domain ${\R_+\times\R}$.
Elementary computations show that
\beq
\label{eq:gnl}
\left. 
\begin{array}{c}
(\lambda_\rho,\lambda_u)\cdot r = 0
\\
(\mu_\rho,\mu_u)\cdot s = 0
\end{array}
\right\}
\quad
\Leftrightarrow
\quad
\rho=-\frac{4\gamma (2\gamma-1)^2}{(\gamma+1)^2} u^2
\text{ and }
\left\{
\begin{array}{c}
u\le0\\ u\ge0
\end{array}
\right. 
.
\eeq
Thus, for $\gamma\ge0$, $\gamma\not=0, 1/2$ the
system is genuinely nonlinear on the closed domain ${\R_+\times\R}$;
for
$\gamma=0,1/2$ it is genuinely nonlinear in the interior of
${\R_+\times\R}$ (with genuine nonlinearity marginally lost on the
boundary,
$\rho=0$). For $\gamma<0$ genuine nonlinearity is  lost in the
interior of ${\R_+\times\R}$.

\medskip
\noindent \underline{\emph{Lax entropies and entropy solutions:}}
Lax entropies of the pde (\ref{eq:pde}) are solutions of the
linear hyperbolic partial differential equation \beqs \rho
S_{\rho\rho} + (2\gamma-1)u S_{\rho u} - S_{uu} =0. \eeqs It turns
out that the system is sufficiently rich in Lax entropies. In
particular a Lax entropy globally convex in $\R_+\times\R$
is
\beq
\label{eq:cle}
S(\rho,u)=\rho\log\rho+\frac{u^2}{2}.
\eeq
Construction of other Lax entropies with particular features (e.g.
possessing scale similarity, or polynomial in $\sqrt\rho$ and $u$,
etc.) is a very instructive exercise.

\medskip
\noindent
\underline{\emph{The Maximum Principle and positively invariant
    domains:}}
For $\gamma\ge0$ our systems satisfy the conditions of the
Lax's Maximum
Principle  proved in \cite{lax1}.  Namely:
(i)
they do posses a globally strictly convex Lax entropy bounded from
below, see (\ref{eq:cle});
(ii)
the Riemann invariants $w(\rho,u)$ and $z(\rho,u)$ satisfy the
convexity condition (\ref{eq:cri});
(iii)
they are genuinely nonlinear in the interior of ${\cal D}$, see
    (\ref{eq:gnl}).

Hence it follows that \emph{convex domains bounded by level curves of
$w(\rho,u)$ and $z(\rho,u)$ are positively invariant for entropy
solutions.}

First we conclude, that ${\cal D}$ itself is positively invariant
domain, as it should be.

Second: a very essential difference between the cases $\gamma<1$,
$\gamma=1$ and $\gamma>1$ follows, which is of crucial importance for
the  main result of the present paper. In the case $\gamma<1$ all
convex domains bounded by
level curves of the Riemann invariants are \emph{unbounded
  (non-compact)}  and thus there is no a priori bound on the
solutions. Even starting with smooth initial data with compact support
nothing
prevents the entropy solutions to blow up indefinitely after
appearence of the shocks. On the other
hand, if $\gamma\ge1$ any compact subset of ${\cal D}$ is contained in
a compact convex domain bounded by level sets of the Riemann
invariants, which fact yields a priori bounds on the entropy
solutions, given bounded initial data.
A microscopic consequence of this fact is that the proof of our main
theorem is valid only for $\gamma>0$.


\bigskip

\noindent {\bf\large Acknowledgement:}
It is our pleasure to thank J\'ozsef Fritz for the many
discussions  on the content of this paper, his permanent interest
and  encouragement.
We also thank Peter Lax  for a very inspirative consultation on
hyperbolic conservation laws.
\\
The kind hospitality of Institut Henri
Poincar\'e (Paris) and that  of the Isaac
Newton Institute (Cambridge),  where parts of this work were
completed,  is gratefully acknowledged.
\\
The research work of the authors is partially supported by the
Hungarian  Scientific Research Fund (OTKA) grant no. T037685.


\vskip1cm

\hbox{\sc
\vbox{\noindent
\hsize66mm
B\'alint T\'oth\\
Institute of Mathematics\\
Technical University Budapest\\
Egry J\'ozsef u. 1.\\
H-1111 Budapest, Hungary\\
{\tt balint{@}math.bme.hu}
}
\hskip5mm
\vbox{\noindent
\hsize66mm
Benedek Valk\'o\\
Institute of Mathematics\\
Technical University Budapest\\
Egry J\'ozsef u. 1.\\
H-1111 Budapest, Hungary\\
{\tt valko{@}math.bme.hu}
}
}


\begin{thebibliography}{99}


\bibitem{balazs}
M. Bal\'azs:
Growth fluctuations in interface models.
{\sl Annales de l'Institut Henri Poincar\'e --- Probabilit\'ees et
  Statistiques}
{\bf 39}: 639-685 (2003)

\bibitem{cocozza}
C. Cocozza:
Processus des misanthropes.
{\sl Zeitschrift f\"ur Wahrscheinlichkeitstheorie und verwandte
  Gebiete}
{\bf 70}: 509-523 (1985)

\bibitem{evans}
L.C. Evans:
{\sl Partial Differential Equations.}
Graduate Studies in Mathematics
{\bf 19},  AMS, Providence RI, 1998


\bibitem{fritz1}
J. Fritz:
{\sl An Introduction to the Theory of Hydrodynamic Limits.}
Lectures in Mathematical Sciences
{\bf 18}. Graduate School of Mathematics, Univ. Tokyo, 2001.

\bibitem{fritz2}
J. Fritz:
Entropy pairs and compensated compactness for weakly asymmetric
systems.
{\sl Advanced Studies in Pure Mathematics}
(2003) (to appear), www.math.bme.hu/~jofri.



\bibitem{fritztoth}
J. Fritz, B. T\'oth: Derivation of the Leroux system as the
hydrodynamic limit of a two-component lattice gas.
to appear in
{\sl Communications in Mathematical Physics} (2004)
\\
{\tt http://arxiv.org/abs/math.PR/0304481}





\bibitem{garabedian}
P.R. Garabedian: {\sl Partial Differential Equations.}  AMS
Chelsea, Providence RI, 1998



\bibitem{john}
F. John:
{\sl Partial Differential Equations.}
Applied Mathematical Sciences, vol. 1,
Springer, New York-Heidelberg-Berlin, 1971.

\bibitem{kardarparisizhang}
M. Kardar, G. Parisi, Y.-C. Zhang:
Dynamic scaling of growing interfaces.
{\sl Physical Reviews Letters\/}
{\bf 56}: 889-892 (1986)

\bibitem{kipnislandim}
C. Kipnis, C. Landim:
{\sl Scaling Limits of Interacting Particle Systems.\/}
Springer, 1999.

\bibitem{lax1}
P. Lax:
Shock waves and entropy.
In: {\sl Contributions to   Nonlinear Functional Analysis},
ed.: E.A. Zarantonello.
Academic Press, 1971, pp. 603-634

\bibitem{lax2}
P. Lax:
{\sl Systems of Conservation Laws and the Mathematical Theory of Shock
  Waves.}
SIAM, CBMS-NSF 11, 1973.

\bibitem{leveque}
R.J. Leveque:
{\sl Numerical Methods in Conservation Laws.}
Lectures In Mathematics, ETH Z\"urich, Birkh\"auser Verlag Basel, 1990

\bibitem{levinesleeman}
H. A. Levine, B. D. Sleeman:
A system of reaction diffusion equations arising in the theory of
reinforced random walks.
{\sl SIAM Journal of Applied Mathematics} {\bf 57} 683-730 (1997)


\bibitem{othmerstevens}
H. G. Othmer, A. Stevens:
Aggregation, blowup, and collapse: the abc's of taxis in reinforced
random walks.
{\sl SIAM Journal of Applied Mathematics} {\bf 57}: 1044-1081 (1997)


\bibitem{popkovschutz}
V. Popkov, G.M. Sch\"utz:
Shocks and excitation dynamics in driven diffusive two channel
systems.
{\sl Journal of Statistical Physics} {\bf 112}: 523-540 (2003)


\bibitem{rascle}
M. Rascle:
On some ``viscous'' perurbations of quasi-linear first order
hyperbolic systems arising in biology.
{\sl Contemporary Mathematics} {\bf 17}: 133-142 (1983)

\bibitem{rezakhanlou}
F. Rezakhanlou:
Microscopic structure of shocks in one conservation laws.
{\sl Annales de l'Institut Henri Poincar\'e --- Analyse Non Lineaire}
{\bf 12}: 119-153 (1995)

\bibitem{serre}
D. Serre:
{\sl Systems of Conservation Laws.} Vol 1-2.
Cambridge University Press, 2000

\bibitem{smoller}
J. Smoller:
{\sl Shock Waves and Reaction Diffusion Equations}, Second Edition,
Springer, 1994.

\bibitem{tothvalko1}
B. T\'oth, B. Valk\'o:
Between equilibrium fluctuations and Eulerian scaling. Perturbation of
equilibrium for a class of deposition models.
{\sl Journal of   Statistical Physics}
{\bf 109}: 177-205 (2002)

\bibitem{tothvalko2}
B. T\'oth, B. Valk\'o:
Onsager relations and Eulerian hydrodynamic limit for systems with
several conservation laws.
{\sl Journal of Statistical Physics}
{\bf 112}: 497-521 (2003)

\bibitem{tothwerner1}
B. Tóth, W. Werner:
The true self-repelling motion.
{\sl Probability Theory and Related Fields}
{\bf 111}: 375-452 (1998)

\bibitem{tothwerner2}
B. Tóth, W. Werner:
Hydrodynamic equation for a deposition model.
In: {\sl In and out of equilibrium. Probability with a physics
  flavor,} V. Sidoravicius Ed., Progress in Probability
{\bf 51}, Birkhäuser, 227-248 (2002)

\bibitem{varadhan}
S.R.S. Varadhan:
Nonlinear diffusion limit for a system with nearest neighbor
interactions II. In:
{\sl Asymptotic Problems in Probability Theory, Sanda/Kyoto 1990}
75--128. Longman, Harlow 1993.

\bibitem{yau1}
H.T. Yau: Relative entropy and hydrodynamics of Ginzburg-Landau
models.
{\sl Letters in  Mathematical Physics\/}
{\bf 22}: 63-80 (1991)

\bibitem{yau2}
H.T. Yau:
Logarithmic Sobolev inequality for generalized simple exclusion
processes.
{\sl Probability Theory and Related Fields}
{\bf 109}: 507-538 (1997)

\bibitem{yau3}
H.T. Yau:
Scaling limit of particle systems, incompressible Navier-Stokes
equations and Boltzmann equation.
In: {\sl Proceedings of the International Congress of Mathematics,
  Berlin 1998}, vol 3, pp 193-205, Birkh\"auser  (1999)

\end{thebibliography}
\end{document}